\documentclass[12pt]{report}     
\usepackage{graphicx}               

\usepackage[centertags]{amsmath}
\usepackage{amsfonts}
\usepackage{latexsym}
\usepackage{amssymb}
\usepackage{amsthm}
\usepackage{verbatim}
\usepackage{amsthm}
\usepackage{bbm}
\theoremstyle{plain}
\newtheorem{theorem}{Theorem}[chapter]
\newtheorem{algorithm}[theorem]{Algorithm}

\newtheorem{corollary}[theorem]{Corollary}
\newtheorem{lemma}[theorem]{Lemma}

\newtheorem{proposition}[theorem]{Proposition}

\theoremstyle{definition}
\newtheorem{definition}[theorem]{Definition}
\newtheorem{example}[theorem]{Example}
\newtheorem*{remark}{Remark}

\numberwithin{equation}{section}

\newcommand{\CC}{\mathbb{C}}
\newcommand{\PP}{\mathbb{P}}
\newcommand{\ZZ}{\mathbb{Z}}
\newcommand{\RR}{\mathbb{R}}
\newcommand{\QQ}{\mathbb{Q}}
\newcommand{\KK}{\mathbb{K}}
\newcommand{\Aff}{\mathbb{A}}
\newcommand{\HH}{\mathbb{H}}
\newcommand{\Oh}{\mathcal{O}}

\newcommand{\Spec}{\textnormal{Spec}}
\newcommand{\Proj}{\textnormal{Proj}}

\newenvironment{Eg}[1]
{
    \begin{flushleft}
    \bfseries{Maximal Group #1}
    \addcontentsline{toc}{subsection}{\ \ \ \ \ \ Symmetry Group #1} 
    \end{flushleft}
}
{
}
\newenvironment{subEg}[1]
{
    \begin{flushleft}
    \bfseries{$\bullet$ Subgroup #1}
    \end{flushleft}
}
{

}

\newenvironment{Family}[2]
{
    \begin{flushleft}
    \bfseries{Family #1}\rm: Invariant under #2
    \end{flushleft}
}
{

}
\newenvironment{Eqn}[1]
{
    \begin{flushleft}
    \bfseries{Differential Equation #1}
    \end{flushleft}
}
{

}

\begin{document}
\pagestyle{empty}

\begin{center}
\begin{bfseries}
\begin{Large}
Picard-Fuchs Differential Equations for Families of K3 Surfaces

\end{Large}
\end{bfseries}
\bigskip
\begin{normalsize}
by
\end{normalsize}
\bigskip
\begin{bfseries}
\begin{Large}

James Patrick Smith
\vspace{12cm}
\end{Large}
\end{bfseries}

\begin{normalsize}
Submitted to the University of Warwick\\
for the degree of
\end{normalsize}

\begin{bfseries}
Doctor of Philosophy
\bigskip

Mathematics Institute,\\
University of Warwick\\
July 2006
\end{bfseries}
\end{center}

\pagenumbering{roman}
\tableofcontents                   
\setcounter{page}{1}
\listoftables                    
\listoffigures                   

\chapter*{Acknowledgments}
I am indebted to my Ph.D. supervisor, Dr. Katrin Wendland, for her advice and support throughout the writing of this thesis. I have also benefited from discussions with Miles Reid, Gavin Brown, Jaros\l aw Buczy\'nski, Weronika Buczy\'nska, Micha\l\ Kapustka, Grzegorz Kapustka and many of the participants of the Calf seminars. The University of North Carolina at Chapel Hill welcomed me for a two month stay.

Most importantly, I am very grateful to Helen for marrying me.

\chapter*{Declaration}
I declare that the contents of this thesis is original, except where explicitly stated. This work has not been submitted in any form for a degree at any other university.

\newpage{\pagestyle{empty}\cleardoublepage} 

\pagenumbering{arabic}
\setcounter{page}{1}
\pagestyle{plain}
\setcounter{chapter}{-1}
\chapter{Introduction}
We study the Picard--Fuchs differential equation of families of K3 surfaces. Specifically, we focus on one parameter families of K3 surfaces with generic Picard number 19. The Picard--Fuchs differential equation is the linear ODE satisfied by the periods of a family of Calabi--Yau manifolds.

In our specific case, the Picard--Fuchs equation has order three and its solutions satisfy a quadratic relationship and so are parametrised by a second order ODE. It is this, the so--called symmetric square root differential equation that we study. We are interested in understanding what differential equations can occur. Our approach is concrete and focuses on a number of specific examples of families. The first problem is to identify a source of good examples.

Chapter \ref{autChapter} examines K3 surfaces with a non-trivial finite group of symplectic  automorphisms. After some initial definitions, results are drawn together that allow us to show that families of symmetric K3 surfaces are lattice polarised. This means that their Picard lattices contain a fixed primitive sublattice of signature $(1,k)$. The finite groups that can act symplectically on a K3 surface have been classified in \cite{Mukai}. We compile a list of groups that provide a rank 19 lattice polarisation. We find projective representations of these groups and write down families of K3 surfaces defined by polynomial invariants of these representations.

By considering our families of K3 surfaces as lattice polarised families, we are able to prove in theorem \ref{squares} that the monodromy group of a Picard--Fuchs equation is integral and so the traces of the monodromies of its symmetric square root are square roots of integers.

In chapter 2 we describe the Griffiths--Dwork method for computing the Picard--Fuchs differential equation from the polynomials defining a family of Calabi--Yau hypersurfaces in weighted projective space. We turn this method into an efficient algorithm implemented in Macaulay2 and included as an appendix. The Picard--Fuchs equation for each of the examples of chapter \ref{autChapter} is then calculated. We find that the square root differential equations are, for our examples, either hypergeometric or generalised Lam\'e differential equations.

When we quotient a family of K3 surfaces by a group of symplectic automorphisms and minimally resolve the resulting singular points, we obtain a new family of K3 surfaces covered by the first. We prove at the end of chapter 2 that these two families have the same Picard--Fuchs differential equation. We find that a number of our examples of families have the same Picard--Fuchs equation and we find geometric relationships between them to explain these coincidences.

In chapter 3 we find the monodromy representations for the Picard--Fuchs differential equations. For hypergeometric ODEs, it is known that the local monodromies uniquely determine the global monodromy group. We use this rigidity to calculate the monodromy group and we classify all the hypergeometric ODEs with Fuchsian monodromy group satisfying the restrictions of theorem \ref{squares}.

For generalised Lam\'e equations, the second type of differential equation we have found in our examples, the local monodromy group does not determine the global monodromy representation. However, we prove that the local monodromies together with the values of the square of the trace of pairs of monodromies uniquely determine the global monodromy group.  The restrictions of theorem \ref{squares} show that the additional numbers we need to find are integers. We are able to numerically approximate these integers easily using an algorithm given in appendix \ref{maple}. This technique allows us to find the monodromy group with certainty and avoids the need to numerically approximate the monodromy matrices.

Chapter \ref{chapter:quots} looks at families of K3 surfaces that occur as quotients of symmetric surfaces. This can be viewed as a generalisation of a K3 surface occurring as the quotient of another K3 surface by a group of symplectic automorphisms, as studied in chapter \ref{autChapter}. We develop a simple combinatorial method to find examples of such quotients.

	\newpage

\chapter{K3 Surfaces with Automorphisms}\label{autChapter}
\section{Basic Definitions and Theorems}\label{SymA}
We begin by providing some definitions and clarifying some important points. Much of the well--known theory behind K3 surfaces can be found in \cite{BPV} and is not included here. Although some of the definitions and results in this thesis apply more generally, the emphasis will be strictly on \emph{complex} K3 surfaces.
\begin{definition}
By a K3 surface, we mean a simply--connected compact complex surface with trivial canonical bundle.
\end{definition}
We are particularly interested in the group of symplectic automorphisms of a K3 surface.
\begin{definition}
An automorphism of a K3 surface $X$ is said to be \emph{symplectic} whenever the induced automorphism on
$$H_{DR}^2(X)\otimes{\CC}\cong H^{2,0}(X)\oplus H^{1,1}(X)\oplus H^{0,2}(X)$$
pointwise fixes $H^{2,0}(X)$.
\end{definition}
In this section, $X$ shall always denote a K3 surface and $G$ a finite group of symplectic automorphisms of $X$. Writing $S_X = H^{1,1}(X)\cap{H^2(X,\ZZ)}$ for the Picard lattice of $X$ and $T_X = S_X^\perp\cap H^2(X,\ZZ)$ for the transcendental lattice, we introduce two further lattices due to Nikulin \cite{NikAut}:
$$T_{X,G} := H^2(X,\ZZ)^G$$
and
$$S_{X,G} := T_{X,G}^\perp\cap H^2(X,\ZZ).$$
\begin{proposition}\label{MukaiNumber}
\ \newline
i) $S_{X,G}\subset{S_X}$ and $T_{X,G}\supset{T_X}$ are primitive sublattices.\\
ii) $S_{X,G}$ is negative definite and contains no classes of self-intersection $-2$.\\
iii) $\textnormal{rank}(S_{X,G}) = 24-\mu(G)$ where
$$
\mu(G):=\frac{1}{|G|}\sum_{g\in G}\mu(\textnormal{ord}(g))
$$
with
$$
\mu(n):=\frac{24}{n\prod_{p|n}(1+1/p)},
$$
the product being taken over primes $p$.
\end{proposition}
\begin{proof}
The proof of i) and ii) can be found in \cite{NikAut}. iii) is due to Mukai \cite{Mukai}. Mukai's approach to understanding symplectic automorphisms is to consider the faithful representation of $G$ naturally induced on the vector space $\textnormal{H}^*(X,\ZZ)\otimes\QQ$. Given that $\textnormal{H}^0(X,\ZZ)$ and $\textnormal{H}^4(X,\ZZ)$ are fixed by $G$, it is seen that 
\begin{eqnarray*}
\textnormal{rank}(S_{X,G}) & = & \textnormal{rank}(\textnormal{H}^*(X,\QQ)) - \textnormal{rank}(\textnormal{H}^*(X,\QQ)^G)\\
& = & 24 -  \frac{1}{|G|} \sum_{g \in G}\textnormal{tr}_{\textnormal{H}^*(X,\QQ)}(g).
\end{eqnarray*}
The proof is completed by the fact that $\textnormal{tr}_{\textnormal{H}^*(X,\QQ)}(g) = |\textnormal{Fix}(g)| = \mu(g)$ as demonstrated in \cite{Mukai}.
\end{proof}
\begin{corollary}
The sublattice $S_{X,G}$ provides a lower bound on the rank of the Picard lattice:
$$\textnormal{rank}(S_{X,G})\geq{24-\mu(G)}.$$
Furthermore, since $S_{X,G}$ is negative definite, if $X$ is algebraic, the lower bound can be increased by $1$ since $S_X$ also contains a positive generator coming from a hyperplane section.
\end{corollary}
\begin{definition}
We shall call the number $\mu(G)$ occurring in proposition \ref{MukaiNumber} the \emph{Mukai number} of the group $G$.
\end{definition}
In many cases, our interest in symplectic automorphisms is due to the following fact:
\begin{proposition}
A finite group of automorphisms, $G$, on a K3 surface $X$ is symplectic whenever $X/G$ is a K3 surface with Du Val singularities.
\end{proposition}
This is because any non--trivial symplectic automorphism has isolated fixed points. If $p$ is a fixed point of $g\in{G}$, then the induced action of the stabiliser of $p$, $G_p$, on the tangent space $T_p$ at $p$
$$G_p \hookrightarrow \textnormal{Gl}(T_p)$$
is in fact contained in $\textnormal{Sl}(T_p)$, ensuring that the quotient is locally of the form $\CC^2/H$ for $H\subset\textnormal{Sl}(2,\CC)$ as required.

There is a complete classification of the finite groups that act symplectically on some K3 surface.
\begin{theorem}\label{thm::classification}[\cite{Mukai}]
If $G$ is a finite group acting faithfully as symplectic automorphisms of a K3 surface, then $G$ is a subgroup of one of the following 11 maximal groups.
$$
T_{48},\ N_{72},\ M_9,\ \mathfrak{S}_5,\ L_2(7),\ H_{192},\ T_{192},\ \mathfrak{A}_{4,4},\ \mathfrak{A}_6,\ F_{384},\ M_{20}.
$$
\end{theorem}
The notation for these groups is the same as that of \cite{Mukai} and \cite{Xiao}. The first reference contains a description of the maximal groups. In total, the 11 maximal groups have 79 subgroups, all of which are listed in \cite{Xiao}. In the next section, we shall investigate a few of these groups.

\section{The Finite Groups of Symplectic Automorphisms}
Our aim is to create a list of examples of one parameter families of K3 surfaces with generic Picard number $19$. To do this, we make use of the Mukai number of a group (see proposition \ref{MukaiNumber}). Let $G$ be a group acting as symplectic automorphisms of an algebraic K3 surface $X$. According to proposition \ref{MukaiNumber}, the Picard lattice, $S_X$, contains the sublattice $S_{X,G}$ and a very ample divisor class. Together, this divisor and sublattice generate a primitively embedded lattice, $M\subset{S_X}$, with $\textnormal{rank}(M) = 25 - \mu(G)$ and we may view $X$ as a lattice polarised K3 surface in the sense of \cite{Dol}. The moduli space of these $M$-polarised K3 surfaces has dimension $\mu(G) - 5$. Thus, for one dimensional moduli spaces, we are interested in groups with $\mu(G)=6$.

We are going to give examples of families of quartic hypersurfaces in $\PP^3$ and double covers of $\PP^2$ branched over a sextic (expressed as a sextic in the weighted projective space $\PP(1,1,1,3)$). We break the problem into the following pieces:
\begin{enumerate}
	\item Look for a subgroup, $G$, of one of the $11$ maximal groups of \cite{Mukai} with $\mu(G) = 6$.
	\item Find a faithful projective representation $\rho:G\hookrightarrow\PP\textnormal{Sl}(\CC^{k+1})$ ($k = 2$ or $3$). Calculate the polynomial invariants of $\rho(G)$ of homogeneous degree $6$ resp. $4$ (for $k = 2$ resp. $3$). There should be two such invariants, $p$ and $q$, with the additional property that the hypersurfaces $(\lambda_0{p} + \lambda_1{q} = 0)\subset\PP^k$ are nonsingular for general $(\lambda_0,\lambda_1)\in\PP^1$.
\end{enumerate}
Our method will produce families of K3 surfaces, either of the form $(\lambda_0{p} + \lambda_1{q} = t^2)\subset\PP(1,1,1,3)$ or $(\lambda_0{p} + \lambda_1{q} = 0)\subset\PP^3$, invariant under $G$. By the corollary to proposition \ref{MukaiNumber} and given that the families are not isotrivial, these K3 surfaces will have generic Picard number $19$. We shall deal with part (2) in the separate cases of $\PP^3$ and $\PP(1,1,1,3)$ in sections \ref{p1113} and \ref{p3}. First we tackle part (1).

According to part 1 of our method, we must find those groups that have a symplectic action on some K3 surface and have Mukai number $6$. These groups are all listed in \cite{Xiao}, but our problem is really to specify these groups in a way that can be handled by Magma. Abstract group names are not enough, we need concrete generators and relations. Furthermore, part 2 of the method requires us to find a projective representation of $G$. To do this, we will find it useful to restrict a projective representation of one of the $11$ maximal groups to its subgroups. With this in mind, we shall take each of the maximal groups and list their subgroups.

To list these groups, we make use of the \verb|SmallGroups| database in Magma. This is a database of the isomorphism types of all groups of order less than $2000$ (excluding $1024$). Since the largest group with a symplectic action has order 960, all the groups we are interested in are contained in this database. First, we identify the 11 maximal groups of theorem \ref{thm::classification} within the \verb|SmallGroups| database. This is detailed below. Typically, we find a maximal group within the \verb|SmallGroups| database by looking at all groups of the correct order. Since the maximal groups all have Mukai number 5, we eliminate those with the wrong Mukai number. Sometimes this specifies the group uniquely. Otherwise, we look at the orders of its elements or of its conjugator subgroup to pin down the group.

(i)\ $T_{48}.$\\
According to the SmallGroups database, there are 52 groups of order $|T_{48}| = 48$. Of these groups, only SmallGroup($48, n$) for $n = 14, 29, 30, 35, 37, 49$ have Mukai number 5 as required. The table in \cite{Xiao} tells us that $T_{48}$ has an element of order 8. Out of the six remaining groups, only number 29 has this property. Hence $T_{48} \cong \textnormal{SmallGroup}(48,29)$.

To give a feel for the Magma code required, we include it in this case only:
\begin{verbatim}
> load"MukaiNumber";
Loading "MukaiNumber"
> NumberOfSmallGroups(48);
52
> sg48 := SmallGroups(48);
> List := [ n : n in [1..#sg48] | Mukai(sg48[n]) eq 5 ];List;
[ 14, 29, 30, 35, 37, 49 ]
> [ n : n in List | 8 in { Order(g) : g in SmallGroup(48,n) } ];
[ 29 ]
> T48 := SmallGroup(48,29);
\end{verbatim}
The first line loads the file ``MukaiNumber'' which contains a function \verb|Mukai()| that calculates the Mukai number of a finite group.

(ii) and (iii) $N_{72}$ and $M_9.$\\
Both these groups have order 72. There are 50 groups of order 72 and of these, numbers 35, 40, 41, and 44 have Mukai number 5. Number 35 can be discounted as it has subgroups with non--integer Mukai number and so can't act symplectically on any K3 surface. The table in \cite{Xiao} specifies the orders of the elements of all the maximal groups (and their subgroups). From this, we find that $N_{72}\cong\textnormal{SmallGroup}(72,40)$ and $M_9\cong\textnormal{SmallGroup}(72,41)$.

(iv) $\mathfrak{S}_5.$\\
In Magma, this symmetric group is specified as \verb|SymmetricGroup(5)|. For the record, we have
\begin{verbatim}
> IdentifyGroup(SymmetricGroup(5));
<120, 34>
\end{verbatim}
so that $\mathfrak{S}_5\cong\textnormal{SmallGroup}(120, 34)$.

(v) $L_2(7).$\\
SmallGroup($168, 42$) is the only group of order 168 with Mukai number 5.

(vi) and (vii)  $H_{192}$ and $T_{192}.$\\
Magma shows that there are only two groups of order 192 with Mukai number 5 all of whose subgroups have integral Mukai number. The groups $H_{192}$ and $T_{192}$ have the same order structure. However, they can be distinguished by the order of their commutator subgroups as given in \cite{Xiao}. From this, we find $H_{192} \cong \textnormal{SmallGroup}(192,955)$ and 
$T_{192} \cong \textnormal{SmallGroup}(192,1493).$

(viii) $\mathfrak{A}_{4,4}.$\\
SmallGroup($288, 1026$) is the only group of order 288 with Mukai number 5 and subgroups with integral Mukai number.

(ix) $\mathfrak{A}_6.$\\
$\mathfrak{A}_6\cong\textnormal{AlternatingGroup}(6)\cong\textnormal{SmallGroup}(360,118)$

(x) $F_{384}.$\\
SmallGroup($384, 18135$) is the only group of order 384 with Mukai number 5.

(xi) $M_{20}.$\\
SmallGroup($960, 11357$) is the only group of order 960 with Mukai number 5.

We summarise these identifications in the following table:
\begin{table}[h]
  \newcommand\T{\rule{0pt}{2.3ex}}
  \newcommand\B{\rule[-1.0ex]{0pt}{0pt}}
  	\centering
		\begin{tabular}{|c|c|c|}
	\hline \T\B Group & Order & SmallGroup number\\\hline\hline
\T\B $T_{48}$ & 48 & 29\\\hline
\T\B $N_{72}$ & 72 & 40\\\hline
\T\B $M_9$ & 72 & 41\\\hline
\T\B $\mathfrak{S}_5$ & 120 & 34\\\hline
\T\B $L_2(7)$ & 168 & 42\\\hline
\T\B $H_{192}$ & 192 & 955\\\hline
\T\B $T_{192}$ & 192 & 1493\\\hline
\T\B $\mathfrak{A}_{4,4}$ & 288 & 1026\\\hline
\T\B $\mathfrak{A}_6$ & 360 & 118\\\hline
\T\B $F_{384}$ & 384 & 18135\\\hline
\T\B $M_{20}$ & 960 & 11357\\\hline
		\end{tabular}
	\caption{The Maximal Groups of Symplectic Automorphisms}
	\label{tab:MaximalGroups}
\end{table}

\begin{figure}
	
\begin{picture}(50,100)(10,-40)

  \put(25,40){\makebox(0,0){$T_{48}$}}
  \put(25,10){\makebox(0,0){$SD_{16}$}}
  \put(25,-20){\makebox(0,0){$C_8$}}
  \put(25,16){\line(0,1){18}}
  \put(25,-12){\line(0,1){17}}

\end{picture}
\begin{picture}(50,100)(0,-40)

  \put(25,40){\makebox(0,0){$N_{72}$}}
  \put(10,10){\makebox(0,0){$3^2C_4$}}
  \put(40,10){\makebox(0,0){$\mathfrak{S}_{3,3}$}}
  \put(40,-20){\makebox(0,0){$C_3\!\times\!D_6$}}
  \put(17,17){\line(1,3){5}}
  \put(33,17){\line(-1,3){5}}
  \put(40,-15){\line(0,1){18}}

\end{picture}
\begin{picture}(50,70)(-10,-40)

  \put(25,40){\makebox(0,0){$M_9$}}
  \put(25,10){\makebox(0,0){$3^2C_4$}}
  \put(25,16){\line(0,1){18}}
 
\end{picture}
\begin{picture}(80,70)(-10,-40)

  \put(25,40){\makebox(0,0){$\mathfrak{S}_5$}}
  \put(10,10){\makebox(0,0){$\mathfrak{A}_5$}}
  \put(45,10){\makebox(0,0){Hol$(C_5)$}}
  \put(17,17){\line(1,3){5}}
  \put(33,17){\line(-1,3){5}}

\end{picture}
\begin{picture}(50,100)(0,-40)

  \put(25,40){\makebox(0,0){$L_2(7)$}}
  \put(25,10){\makebox(0,0){$C_7\!\rtimes\!{C_3}$}}
  \put(25,-20){\makebox(0,0){$C_7$}}
  \put(25,16){\line(0,1){18}}
  \put(25,-12){\line(0,1){17}}

\end{picture}
\begin{picture}(90,100)(-55,-90)

  \put(0,-10){\makebox(0,0){$H_{192}$}}
  \put(-30,-40){\makebox(0,0){$\Gamma_{\!\!25}a_1$}}
  \put(0,-40){\makebox(0,0){$2^4\!D_6$}}
  \put(30,-40){\makebox(0,0){$C_2\!\!\times\!\!\mathfrak{S}_4$}}
  \put(0,-16){\line(0,-1){17}}
  \put(-5,-16){\line(-1,-1){17}}
  \put(5,-16){\line(1,-1){17}}

  \put(-30,-46){\line(0,-1){18}}
  \put(30,-46){\line(0,-1){18}}
  \put(-30,-70){\makebox(0,0){$\Gamma_{\!7}a_1$}}
  \put(30,-70){\makebox(0,0){$C_2\!\!\times\!\!\mathfrak{A}_4$}}

\end{picture}
\begin{picture}(60,100)(-0,-71)

  \put(25,40){\makebox(0,0){$T_{192}$}}
  \put(10,10){\makebox(0,0){$\Gamma_{\!\!25}a_1$}}
  \put(40,10){\makebox(0,0){$C_2\!\!\times\!\!\mathfrak{S}_4$}}
  \put(17,17){\line(1,3){5}}
  \put(33,17){\line(-1,3){5}}

  \put(10,4){\line(0,-1){18}}
  \put(40,4){\line(0,-1){18}}
  \put(10,-20){\makebox(0,0){$\Gamma_{\!7}a_1$}}
  \put(40,-20){\makebox(0,0){$C_2\!\!\times\!\!\mathfrak{A}_4$}}
\end{picture}
\begin{picture}(140,150)(-80,-120)

  \put(0,-10){\makebox(0,0){$\mathfrak{A}_{4,4}$}}
  \put(-30,-40){\makebox(0,0){$\mathfrak{A}_{4,3}$}}
  \put(0,-40){\makebox(0,0){$2^4\!D_6$}}
  \put(30,-40){\makebox(0,0){$C_2\!\!\times\!\!\mathfrak{S}_4$}}
  \put(0,-16){\line(0,-1){17}}
  \put(-5,-16){\line(-1,-1){17}}
  \put(5,-16){\line(1,-1){17}}

  \put(30,-46){\line(0,-1){18}}
  \put(30,-70){\makebox(0,0){$C_2\!\!\times\!\!\mathfrak{A}_4$}}

  \put(-27,-46){\line(1,-3){5}}
  \put(-33,-46){\line(-1,-3){5}}
  \put(-50,-70){\makebox(0,0){$C_3\!\!\rtimes\!\!D_8$}}
  \put(-10,-70){\makebox(0,0){$C_3\!\!\times\!\!\mathfrak{A}_4$}}

  \put(-27,-94){\line(1,3){5}}
  \put(-33,-94){\line(-1,3){5}}
  \put(-30,-100){\makebox(0,0){$C_2\!\!\times\!\!C_6$}}
  
  \put(-60,-94){\line(1,3){5}}
  \put(-60,-100){\makebox(0,0){$Q_{12}$}}

\end{picture}
\begin{picture}(60,70)(-15,-70)

  \put(25,40){\makebox(0,0){$\mathfrak{A}_6$}}
  \put(10,10){\makebox(0,0){$\mathfrak{A}_5$}}
  \put(40,10){\makebox(0,0){$3^2C_4$}}
  \put(17,17){\line(1,3){5}}
  \put(33,17){\line(-1,3){5}}

\end{picture}\\
\begin{picture}(150,170)(-30,-160)

  \put(45,-10){\makebox(0,0){$F_{384}$}}

  \put(47,-16){\line(2,-3){10}}
  \put(52,-16){\line(2,-1){30}}
  \put(43,-16){\line(-2,-3){10}}
  \put(38,-16){\line(-2,-1){30}}

  \put(30,-40){\makebox(0,0){$2^4\!D_6$}}
  \put(60,-40){\makebox(0,0){$\Gamma_{\!\!25}a_1$}}
  \put(60,-46){\line(0,-1){18}}
  \put(60,-70){\makebox(0,0){$\Gamma_{\!7}a_1$}}

  \put(90,-40){\makebox(0,0){$SD_{16}$}}
  \put(90,-45){\line(0,-1){17}}
  \put(90,-70){\makebox(0,0){$C_8$}}

  \put(0,-40){\makebox(0,0){$4^2\!\mathfrak{A}_4$}}
  \put(-15,-70){\makebox(0,0){$\Gamma_{\!\!13}a_1$}}
  \put(15,-70){\makebox(0,0){$4^2C_3$}}
  \put(4,-46){\line(1,-3){5}}
  \put(-4,-46){\line(-1,-3){5}}
  \put(-15,-76){\line(0,-1){18}}
  \put(-15,-100){\makebox(0,0){$\Gamma_4c_2$}}

  \put(-19,-106){\line(-1,-3){5}}
  \put(-11,-106){\line(1,-3){5}}
  \put(-30,-130){\makebox(0,0){$C_2\!\!\times\!\!{Q_8}$}}
  \put(0,-130){\makebox(0,0){$C_4^2$}}
  
  \put(0,-121){\line(1,5){9}}

\end{picture}
\begin{picture}(80,170)(-50,-160)

  \put(0,-10){\makebox(0,0){$M_{20}$}}
  \put(-30,-40){\makebox(0,0){$\mathfrak{A}_5$}}
  \put(0,-40){\makebox(0,0){$4^2\mathfrak{A}_4$}}
  \put(30,-40){\makebox(0,0){$2^4D_6$}}
  \put(0,-16){\line(0,-1){17}}
  \put(-5,-16){\line(-1,-1){17}}
  \put(5,-16){\line(1,-1){17}}
  \put(-15,-70){\makebox(0,0){$\Gamma_{\!\!13}a_1$}}
  \put(15,-70){\makebox(0,0){$4^2C_3$}}
  \put(4,-46){\line(1,-3){5}}
  \put(-4,-46){\line(-1,-3){5}}
  \put(-15,-76){\line(0,-1){18}}
  \put(-15,-100){\makebox(0,0){$\Gamma_4c_2$}}

  \put(-19,-106){\line(-1,-3){5}}
  \put(-11,-106){\line(1,-3){5}}
  \put(0,-130){\makebox(0,0){$C_2\!\!\times\!\!{Q_8}$}}
  \put(-30,-130){\makebox(0,0){$C_4^2$}}

\end{picture}

	\caption{The 11 Maximal Groups and their Subgroups of Mukai Number 6.}
	\label{fig:The11MaximalGroupsAndTheirSubgroupsWithMukaiNumber6}
\end{figure}
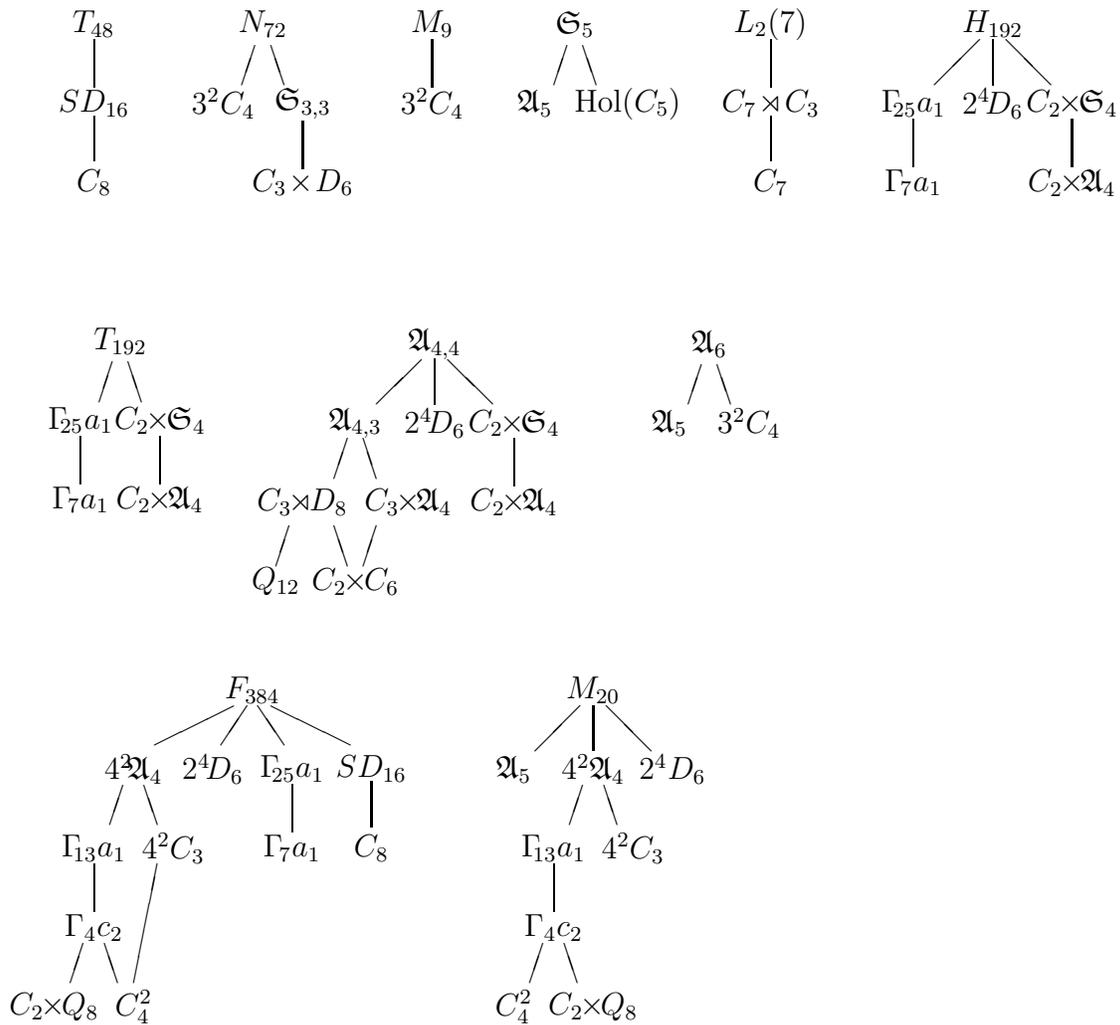

Next, we take each of the 11 maximal groups in turn and list their subgroups via the Magma command \verb|SubgroupLattice()|. This function lists the conjugacy classes of subgroups of a finite group and states which subgroups are contained in which. We then restrict our attention to those subgroups with Mukai number 6 and identify these groups with the groups named in \cite{Xiao}. The results are displayed in figure \ref{fig:The11MaximalGroupsAndTheirSubgroupsWithMukaiNumber6}. Each diagram shows one of the 11 maximal groups and its subgroups with Mukai number 6. The lines denote inclusion with subgroups written below.

\section{Representations and Invariant Theory in Magma}

In order to find examples of symmetric K3 surfaces of Picard number 19 in $\PP(1,1,1,3)$ and $\PP^3$, we would like an exhaustive list of actions of the groups of figure \ref{fig:The11MaximalGroupsAndTheirSubgroupsWithMukaiNumber6} on $\PP^2$ and $\PP^3$.

We are interested in projective representations $G\subset\PP\textnormal{Sl}(n,\CC)$ induced from a representation of the lift $\overline{G}\subset\textnormal{Sl}(n,\CC)$. If $Z$ denotes the centre of $\textnormal{Sl}(n,\CC)$, then $G$ and $\overline{G}$ are related by
$$
\overline{G}/(Z\cap\overline{G}) \cong G
$$
and so we are required to find a representation of an extension $\overline{G}$ of our group $G$.

To find projective representations of the groups with Mukai number 6, we look through lists of finite subgroups of $\textnormal{Sl}(3,\CC)$ and $\textnormal{Sl}(4,\CC)$ for groups $\overline{G}$ where $\overline{G}/(Z\cap\overline{G})$ is one of the groups in figure \ref{fig:The11MaximalGroupsAndTheirSubgroupsWithMukaiNumber6}. With one exception, the projective representations we shall find are restrictions of projective representations of one of the 11 maximal groups. Sections \ref{p1113} and \ref{p3} list some projective representations for the maximal groups $T_{48}, M_9, L_2(7)$ and $\mathfrak{A}_6$ in $\PP\textnormal{Sl}(3,\CC)$ and $L_2(7), T_{192}, F_{384}, M_{20}$ and $\mathfrak{S}_5$ in $\PP\textnormal{Sl}(4,\CC)$. Also, we find a projective representation of the non--maximal group $\mathfrak{A}_{4,3}$ in $\PP\textnormal{Sl}(3,\CC)$ that is not induced from a representation of a larger group.

Without resorting to searching through lists of finite subgroups of $\PP\textnormal{Sl}(n,\CC)$, there is another method at our disposal that is valid in most cases. If $G$ is a soluble finite group, it is possible to build up a representation of $G$ step by step from representations of the subgroups in its composition series. This is implemented in Magma as the command \verb|IrreducibleRepresentationsSchur()|. For example, the maximal group $N_{72}$ is soluble and if it has a projective representation $N_{72}\subset\PP\textnormal{Sl}(3,\CC)$ lifting to a representation $\overline{N}_{72}\subset\textnormal{Sl}(3,\CC)$, then the lift will have order $72$ or $3\times{72}$ and will also be soluble. Magma can list the possible extensions $\overline{N}_{72}$ and their irreducible representations. However, this method does not uncover any more examples of low dimensional projective representations other than those already found by more ad--hoc methods.

Magma has a good collection of algorithms for dealing with invariant theory and once we have found a projective representation of one of our groups, we are able to find its polynomial invariants. This is achieved by the commands \verb|FundamentalInvariants()| that list generators of the ring of invariants or \verb|InvariantsofDegree()| that lists independent invariants of a specified degree.

\subsection{K3 Hypersurfaces in $\PP(1,1,1,3)$}\label{p1113}
Throughout this section, we shall let $x_0,x_1,x_2$ denote the weight 1 coordinates in $\PP(1,1,1,3)$ and $t$ shall denote the weight 3 coordinate.

If $f_6(x_0,x_1,x_2,t)=0$ is any hypersurface of weighted degree 6 in $\PP(1,1,1,3)$, then after a weighted projective transformation of the form
$$
t\mapsto\alpha t+\beta g_3(\underline{x})
$$
for $\alpha, \beta \in \CC$ and $g_3$ a weight 3 polynomial in $x_0, x_1, x_2$, it may be assumed that
\begin{equation}\label{p113form}
f_6(\underline{x},t) = g_6(\underline{x})-t^2.
\end{equation}
The hypersurface is then seen to be a double cover of $\PP^2$ (with coordinates $x_0, x_1, x_2$) branched over the sextic curve $g_6(\underline{x}) = 0$. We shall write all K3 hypersurfaces in this section in the form of \ref{p113form}.

Our present aim is to find examples of one--parameter families of K3 surfaces in $\PP(1,1,1,3)$ with generic Picard number 19. As discussed earlier, this is achieved by finding degree 6 invariants, $p$ and $q$, of finite subgroups of $\textnormal{Aut}(\PP(1,1,1,3))$ appearing in Mukai's classification. The family $t^2 = p+\lambda q$ will be a family of K3 surfaces if it is generically nonsingular.

If a transformation in $\textnormal{Aut}(\PP(1,1,1,3))$ is to act symplectically on a K3 surface of the form $t^2 = g_6(\underline{x})$, then it must be either of the form
\begin{eqnarray*}
t & \mapsto & t\\
\underline{x} & \mapsto & A.\underline{x}
\end{eqnarray*}
for $A\in\textnormal{Sl}(3,\CC)$, or
\begin{eqnarray*}
t & \mapsto & -t\\
\underline{x} & \mapsto & B.\underline{x}
\end{eqnarray*}
for $B\in\textnormal{Gl}(3,\CC)$ with $\textnormal{Det}(A) = -1$. However, the second type of automorphism may be composed with the identity transformation of $\PP(1,1,1,3)$ that maps $t\mapsto -t$ and $\underline{x}\mapsto - \underline{x}$ to put it in the first form. Hence, we are interested in finite subgroups of $\textnormal{Sl}(3,\CC)$ that act on K3 surfaces according to Mukai's classification.

The finite subgroups of $\textnormal{Sl}(3,\CC)$ are classified and listed in, for example, \cite{YauYu}. These subgroups fall into four infinite families and a small number of exceptional cases:
\renewcommand{\labelenumi}{(\Alph{enumi}).}
\begin{enumerate}
	\item Diagonal Abelian groups.
	\item Subgroups isomorphic to transitive subgroups of $\textnormal{Gl}(2,\CC)$ under the isomorphism
	$$\begin{pmatrix}
	a&b\\c&d
	\end{pmatrix}\mapsto\begin{pmatrix}(a d - b c)^{-1}&0&0\\0&a&b\\0&c&d\end{pmatrix}.
	$$
	\item Trihedral groups. That is, a group of type (A) with the transformation $$\begin{pmatrix}0&1&0\\0&0&1\\1&0&0\end{pmatrix}$$
	adjoined.
	\item A group of type (C) with transformations of the form
	$$\begin{pmatrix}a&0&0\\0&0&b\\0&c&0\end{pmatrix}$$
	adjoined where $abc=-1$.
	\item The group of order $36\times{3}$ generated by
	$$
	\begin{pmatrix}
	1&0&0\\0&\omega&0\\0&0&\omega^2
	\end{pmatrix},
	\begin{pmatrix}
	0&1&0\\0&0&1\\1&0&0
	\end{pmatrix},
	\frac{1}{\omega-\omega^2}\begin{pmatrix}
	1&1&1\\1&\omega&\omega^2\\1&\omega^2&\omega
	\end{pmatrix}
	$$
	where $\omega = e^\frac{2\pi i}{3}$. Modulo the centre, this group is isomorphic to $3^2C_4$.
	\item The group of order $72\times{3}$ generated by (E) together with
	$$	\frac{1}{\omega-\omega^2}\begin{pmatrix}
	1&\omega&\omega\\\omega^2&\omega&\omega^2\\\omega^2&\omega^2&\omega
	\end{pmatrix}.
	$$
	Modulo the centre, this group is isomorphic to $M_9$.
	\item A group of order $216\times{3}$.
	\item The group of order 60 isomorphic to the alternating group $\mathfrak{A}_5$ and generated by
	$$
	\begin{pmatrix}
	1&0&0\\0&\varepsilon^4&0\\0&0&\varepsilon
	\end{pmatrix},
	\begin{pmatrix}-1&0&0\\0&0&-1\\0&-1&0\end{pmatrix},
	\frac{1}{2 t + 1}\begin{pmatrix}
	1&1&1\\2&s&t\\2&t&s\end{pmatrix}
	$$
	where $\varepsilon = e^\frac{2 \pi i}{5}, s = \varepsilon^2 + \varepsilon^3$ and $t = \varepsilon + \varepsilon^4$.
	\item The group of order 168 isomorphic to the simple group $L_2(7)$ and generated by
$$
\begin{pmatrix} 0&1&0\\0&0&1\\1&0&0 \end{pmatrix}
,\ 
\begin{pmatrix}\zeta & 0 & 0 \\0 & \zeta^2 & 0 \\0 & 0 & \zeta^4 \end{pmatrix}
,$$
$$
\frac{1}{k} \begin{pmatrix} -\zeta^2+\zeta^5 & -\zeta+\zeta^6 & \zeta^3-\zeta^4 \\ -\zeta+\zeta^6 & \zeta^3-\zeta^4 & -\zeta^2+\zeta^5 \\ \zeta^3-\zeta^4 & -\zeta^2+\zeta^5 & -\zeta+\zeta^6 \end{pmatrix}
$$
where $\zeta = e^\frac{2\pi i}{7}$ and $k = \zeta+\zeta^2+\zeta^4-\zeta^3-\zeta^5-\zeta^6$.
	\item The group generated by (H) and the centre of $\textnormal{Sl}(3,\CC)$.
	\item The group generated by (I) and the centre of $\textnormal{Sl}(3,\CC)$.
	\item The group of order $360\times{3}$ generated by (H) together with
	$$\frac{1}{\sqrt{2 t + 1}}\begin{pmatrix}
	1&\lambda_1&\lambda_1\\
	2\lambda_2&s&t\\
	2\lambda_2&t&s
	\end{pmatrix}
	$$
	where $\lambda_1 = t\omega+\frac{1}{2}(t+\omega)$ and $\lambda_2 = -t\omega -\frac{1}{2}(t+\omega+1)$. Modulo the centre, this group is isomorphic to the alternating group $\mathfrak{A}_6$.
\end{enumerate}

We now search for examples of families of K3 surfaces in $\PP(1,1,1,3)$ with a symplectic automorphism group of Mukai number 6 induced from $\textnormal{Sl}(3,\CC)$. The exceptional groups of the classification already provide a few potential examples. For example, group (L) contains (H) which, modulo the centre, corresponds to $\mathfrak{A_6}$ containing $\mathfrak{A}_5$ in $\PP\textnormal{Sl}(3,\CC)$.

\begin{Eg}{$T_{48}$}
$T_{48} = \textnormal{SmallGroup}(48,29)$ in Magma.

\begin{center}
\begin{picture}(50,80)(0,-30)

  \put(25,40){\makebox(0,0){$T_{48}$}}
  \put(25,10){\makebox(0,0){$SD_{16}$}}
  \put(25,-20){\makebox(0,0){$C_8$}}
  \put(25,16){\line(0,1){18}}
  \put(25,-12){\line(0,1){17}}

\end{picture}
\end{center}

The matrices
$$
A =  \begin{pmatrix}-1 & 0 & 0 \\0 & 0 & \frac{-1+i}{\sqrt{2}} \\0 & -\frac{1+i}{\sqrt{2}} & 0 \end{pmatrix}
,\ 
B = \begin{pmatrix}1 & 0 & 0 \\0 & 0 & -i \\0 & -i & 0 \end{pmatrix}
,\ 
C = \begin{pmatrix}1 & 0 & 0 \\0 & -\frac{1+i}{2} & \frac{-1+i}{2} \\0 & \frac{1+i}{2} & \frac{-1+i}{2} \end{pmatrix}
$$
generate $T_{48}$ in $\textnormal{PSl}(3,\CC)$ containing the subgroup $\left<A,B\right>\cong{SD_{16}}.$

In terms of the classification of finite subgroups of $\textnormal{Sl}(3,\CC)$, $T_{48}$ is a group of type (B), ie., a group isomorphic to a transitive subgroup of $\textnormal{Gl}(2,\CC)$.

\begin{subEg}{$SD_{16}$}
The $SD_{16}$ invariant polynomials of degree 6 are  $x_0^6,\  x_0^2x_1^2x_2^2,\  x_1^5x_2 + x_1x_2^5$ giving the invariant family of K3 surfaces
$$X_\lambda : t^2 = x_0^6 + x_1^5x_2 + x_1x_2^5 - 3 \lambda x_0^2x_1^2x_2^2$$
in $\PP(1,1,1,3)$. This family is singular at $\lambda = \infty$ and $\lambda^3 = 1$ and the symplectic automorphism group jumps to $T_{48}$ at $\lambda = 0$. If $\omega$ is a primitive third root of unity, then the projective transformation
$$(x_0,x_1,x_2)\mapsto(\sqrt{\omega}x_0,x_1,x_2)$$
provides an isomorphism $X_\lambda \cong X_{\omega\lambda}$ so that $\lambda^3$ is the natural parameter for this family.

\end{subEg}
\end{Eg}
 
\begin{Eg}{$M_9$}
$M_9 = \textnormal{SmallGroup}(72,41)$ in Magma.

\begin{center}
\begin{picture}(50,50)

  \put(25,40){\makebox(0,0){$M_9$}}
  \put(25,10){\makebox(0,0){$3^2C_4$}}
  \put(25,16){\line(0,1){18}}
 
\end{picture}
\end{center}

In the classification of finite subgroups of $\PP\textnormal{Sl}(3,\CC)$, $M_9$ is group (F). It contains the subgroup $3^2C_4$ of Mukai number 6 occurring as (E) in the classification.

\begin{subEg}{$3^2C_4$}
The $3^2C_4$ invariant polynomials of degree 6 are
\begin{eqnarray*}
p & := & x_0^6 + x_1^6 + x_2^6 - 10(x_0^3x_1^3 + x_0^3x_2^3 + x_1^3x_2^3),\\ 
q & := & x_0x_1x_2(x_0^3 + x_1^3 + x_2^3) - 2(x_0^3x_1^3 + x_0^3x_2^3 + x_1^3x_2^3) + 3x_0^2x_1^2x_2^2
\end{eqnarray*}
        giving the invariant family of K3 surfaces
$$X_\lambda : t^2 = p + 3 \lambda q$$
in $\PP(1,1,1,3)$. This family is singular at $\lambda = \infty, -1, -2$ and where $\lambda^2 + 16\lambda + 16 = 0$ and the symplectic automorphism group jumps to $M_9$ at $\lambda = 0$.

\end{subEg}
\end{Eg}

\begin{Eg}{$L_2(7)$}
$L_2(7) = \textnormal{SmallGroup}(168,42)$ in Magma.

\begin{center}
\begin{picture}(50,80)(0,-30)

  \put(25,40){\makebox(0,0){$L_2(7)$}}
  \put(25,10){\makebox(0,0){$C_7\!\rtimes\!{C_3}$}}
  \put(25,-20){\makebox(0,0){$C_7$}}
  \put(25,16){\line(0,1){18}}
  \put(25,-12){\line(0,1){17}}

\end{picture}\end{center}
In the classification of finite subgroups of $\PP\textnormal{Sl}(3,\CC)$, this is group (I). Labeling the generators as
$$
A =  \begin{pmatrix} 0&1&0\\0&0&1\\1&0&0 \end{pmatrix}
,\ 
B = \begin{pmatrix}\zeta & 0 & 0 \\0 & \zeta^2 & 0 \\0 & 0 & \zeta^4 \end{pmatrix}
,\ 
C = \frac{1}{k} \begin{pmatrix} -\zeta^2+\zeta^5 & -\zeta+\zeta^6 & \zeta^3-\zeta^4 \\ -\zeta+\zeta^6 & \zeta^3-\zeta^4 & -\zeta^2+\zeta^5 \\ \zeta^3-\zeta^4 & -\zeta^2+\zeta^5 & -\zeta+\zeta^6 \end{pmatrix}
$$
where $\zeta^7=1$ is primitive and $k$ is as defined for group (I), the elements $A$ and $B$ generate a subgroup isomorphic to $C_7\rtimes{C_3}$

\begin{subEg}{$C_7\rtimes{C_3}$}
The $C_7\rtimes{C_3}$ invariant polynomials of degree 6 are $x_0^5x_1  + x_1^5x_2 + x_2^5x_0$ and $x_0^2x_1^2x_2^2$ giving the invariant family of K3 surfaces
$$X_\lambda : t^2 = x_0^5x_1  + x_1^5x_2 + x_2^5x_0 - 3\lambda x_0^2x_1^2x_2^2$$
in $\PP(1,1,1,3)$. This family is singular at $\lambda = \infty$ and $\lambda^3 = 1$ and the symplectic automorphism group jumps to $L_2(7)$ at $\lambda = \frac{5}{3}$. If $\omega$ is a primitive third root of unity, then the projective transformation
$$(x_0,x_1,x_2)\mapsto(x_0,\omega x_1,\omega^2x_2)$$
provides an isomorphism $X_\lambda \cong X_{\omega\lambda}$ so that $\lambda^3$ is the natural parameter for this family.
\end{subEg}
\end{Eg}
\begin{Eg}{$\mathfrak{A}_{4,4}$}
$\mathfrak{A}_{4,4}=\textnormal{SmallGroup}(288,1026)$ in Magma.

\begin{center}
\begin{picture}(0,100)(-20,-100)

  \put(0,-10){\makebox(0,0){$\mathfrak{A}_{4,4}$}}
  \put(-30,-40){\makebox(0,0){$\mathfrak{A}_{4,3}$}}
  \put(0,-40){\makebox(0,0){$2^4\!D_6$}}
  \put(30,-40){\makebox(0,0){$C_2\!\!\times\!\!\mathfrak{S}_4$}}
  \put(0,-16){\line(0,-1){17}}
  \put(-5,-16){\line(-1,-1){17}}
  \put(5,-16){\line(1,-1){17}}

  \put(30,-46){\line(0,-1){18}}
  \put(30,-70){\makebox(0,0){$C_2\!\!\times\!\!\mathfrak{A}_4$}}

  \put(-27,-46){\line(1,-3){5}}
  \put(-33,-46){\line(-1,-3){5}}
  \put(-50,-70){\makebox(0,0){$C_3\!\!\rtimes\!\!D_8$}}
  \put(-10,-70){\makebox(0,0){$C_3\!\!\times\!\!\mathfrak{A}_4$}}

  \put(-27,-94){\line(1,3){5}}
  \put(-33,-94){\line(-1,3){5}}
  \put(-30,-100){\makebox(0,0){$C_2\!\!\times\!\!C_6$}}
  
  \put(-60,-94){\line(1,3){5}}
  \put(-60,-100){\makebox(0,0){$Q_{12}$}}

\end{picture}
\end{center}

Although $\mathfrak{A}_{4,4}$ does have a representation $\mathfrak{A}_{4,4}\hookrightarrow\textnormal{Sl}(3,\CC)$, this has no smooth invariants of degree 6.
\begin{subEg}{$\mathfrak{A}_{4,3}$}
However, modulo the centre of $\textnormal{Sl}(3,\CC)$, the matrices
$$
\begin{pmatrix}
1&0&0\\0&\eta&0\\0&0&\eta^5
\end{pmatrix},
\begin{pmatrix}
0&1&0\\0&0&1\\1&0&0
\end{pmatrix},
\begin{pmatrix}
\eta&0&0\\0&0&\eta\\0&\eta&0
\end{pmatrix}.
$$
with $\eta = e^\frac{2 \pi i}{6}$ generate a group isomorphic to the subgroup $\mathfrak{A}_{4,3}\subset\mathfrak{A}_{4,4}$. This subgroup of Mukai number 6 has invariants $x_0^6+x_1^6+x_2^6$ and $x_0^2x_1^2x_2^2$ providing the family of K3 surfaces
$$t^2 = x_0^6+x_1^6+x_2^6-3\lambda x_0^2x_1^2x_2^2 \hspace{10pt}\subset \PP(1,1,1,3).$$
If $\omega$ is a primitive third root of unity, then the projective transformation
$$
(x_0,x_1,x_2)\mapsto(\omega x_0,x_1,x_2)
$$
induces an isomorphism $X_\lambda\cong X_{\omega\lambda}$ so that $\lambda^3$ is a natural parameter for this family.

\begin{remark}
The projective representation above is not induced by the projective representation of any group $G$ of symplectic automorphisms of a K3 surface with $\mathfrak{A}_{4,4}\supseteq G\supsetneq\mathfrak{A}_{4,3}$. Any such a group would have $\mu(G) = 5$ and so could have only one nonsingular invariant of degree 6. It is seen that none of the exceptional groups (E)--(L) contain $\mathfrak{A}_{4,3}$. However, the families of groups of types (A), (C) and (D) each have the singular invariant $x_0^2x_1^2x_2^2$, and so this could be the only degree 6 invariant of $G$. Thus this representation of $G$ would not correspond to an action on a K3 surface. This leaves the possibility that $G$ is isomorphic to a transitive subgroup of $\textnormal{Gl}(n,\CC)$ (a group of type (B)). This is not possible because $\mathfrak{A}_{4,3}$ is not isomorphic to any of the transitive subgroups of $\PP\textnormal{Gl}(n,\CC)$.
\end{remark}
\end{subEg}
\end{Eg}

\begin{Eg}{$\mathfrak{A}_6$}
$\mathfrak{A}_6 = \textnormal{SmallGroup}(360,118)$ in Magma.

\begin{center}
\begin{picture}(60,70)(-5,0)

  \put(25,40){\makebox(0,0){$\mathfrak{A}_6$}}
  \put(10,10){\makebox(0,0){$\mathfrak{A}_5$}}
  \put(40,10){\makebox(0,0){$3^2C_4$}}
  \put(17,17){\line(1,3){5}}
  \put(33,17){\line(-1,3){5}}

\end{picture}
\end{center}

In the classification of finite subgroups of $\PP\textnormal{Sl}(3,\CC)$, $\mathfrak{A}_6$ is group (L). It contains the subgroups $3^2C_4$ and $\mathfrak{A}_5$ of Mukai number 6 occurring as (E) and (H) respectively in the classification.

\begin{subEg}{$3^2C_4$}
$\mathfrak{A}_6$ contains one conjugacy class of subgroups isomorphic to $3^2C_4$. It can be verified that these subgroups are in turn conjugate to the subgroup isomorphic to $3^2C_4$ embedded in $M_9$ and covered in the previous example.
\end{subEg}
\begin{subEg}{$\mathfrak{A}_5$}
The subgroup (H) of the classification is isomorphic to $\mathfrak{A}_5$ and has the following degree 6 invariants:
$$p := 8x_0^6 + 30x_0^2x_1^2x_2^2 + 3x_0(x_1^5 + x_2^5) + 5x_1^3x_2^3
$$
and
$$q^3 :=(x_0^2 + x_1x_2)^3$$
providing us with the family of K3 surfaces
$$t^2 = p+\lambda q^3$$ in $\PP(1,1,1,3)$.
This family is nonsingular except at $\lambda = \infty, -5, -8$, and $-\frac{40}{9}$.
\end{subEg}
\end{Eg}

To summarise, we have found 5 distinct families of K3 hypersurfaces in $\PP(1,1,1,3)$ invariant under the groups $SD_{16}$, $3^2C_4$, $C_7\rtimes{C_3}$, $\mathfrak{A}_5$ and $\mathfrak{A}_{4,3}$. These are shown in table \ref{tab:K3HypersurfacesInPP1113} in a rearranged order together with labels I to V for later reference. In four examples, the representation of $G$ is induced from the representation of one of the 11 maximal groups. This is not the case for example II.
\begin{table*}
	\newcommand\T{\rule{0pt}{2.6ex}}
  \newcommand\B{\rule[-1.2ex]{0pt}{0pt}}
	\centering
		\begin{tabular}[h]{|c|c|c|}\hline
			\T\B & Symmetry Group & Equation\\
			\hline\hline
			\T\B I & $SD_{16}\subset T_{48}$ & $t^2 = x_0^6 + x_1^5x_2 + x_1x_2^5 - 3 \mu^{1/3} x_0^2x_1^2x_2^2$\\\hline
			\T\B II & $\mathfrak{A}_{4,3}$ & $t^2 = x_0^6+x_1^6+x_2^6-3\mu^{1/3} x_0^2x_1^2x_2^2$\\\hline
			\T\B III & $C_7\rtimes{C_3} \subset L_2(7)$ & $t^2 = x_0^5x_1  + x_1^5x_2 + x_2^5x_0 - 3\mu^{1/3} x_0^2x_1^2x_2^2$\\\hline
			\T\B IV & $3^2C_4\subset M_9, \mathfrak{A}_6$ & $t^2=p+3\lambda q$\\
			\T\B &&$p := x_0^6 + x_1^6 + x_2^6 - 10(x_0^3x_1^3 + x_0^3x_2^3 + x_1^3x_2^3)$\\
			\T\B &&$q :=  x_0x_1x_2(x_0^3 + x_1^3 + x_2^3) - 2(x_0^3x_1^3 + x_0^3x_2^3 + x_1^3x_2^3)$ \\
			\T\B &&$+ 3x_0^2x_1^2x_2^2$\\\hline
			\T\B V&$\mathfrak{A}_5\subset \mathfrak{A}_6$ & $t^2 = p+\lambda q^3$\\
			\T\B && $p := 8x_0^6 + 30x_0^2x_1^2x_2^2 + 3x_0(x_1^5 + x_2^5) + 5x_1^3x_2^3$\\
			\T\B && $q := x_0^2 + x_1x_2$\\\hline
		\end{tabular}
	\caption{K3 Hypersurfaces in $\PP(1,1,1,3)$}
	\label{tab:K3HypersurfacesInPP1113}
\end{table*}

\subsection{K3 Hypersurfaces in $\PP^3$}\label{p3}

We now turn our attention to finding examples of symmetric K3 hypersurfaces in $\PP^3$. This amounts to finding degree 4 invariants of finite subgroups of $\textnormal{Sl}(4,\CC)$. Again, we restrict our attention to those groups with Mukai number 6, although we find them as subgroups of groups with Mukai number 5.

\begin{Eg}{$\mathfrak{S}_{5}$}
$\mathfrak{S}_5 = \textnormal{SmallGroup}(120,34)$ in Magma.

\begin{center}
\begin{picture}(80,40)(0,10)

  \put(25,40){\makebox(0,0){$\mathfrak{S}_5$}}
  \put(10,10){\makebox(0,0){$\mathfrak{A}_5$}}
  \put(45,10){\makebox(0,0){Hol$(C_5)$}}
  \put(17,17){\line(1,3){5}}
  \put(33,17){\line(-1,3){5}}

\end{picture}
\end{center}
The matrices
$$
A =  \begin{pmatrix}1&0&0&0\\0&1&0&0\\0&0&\omega&0\\0&0&0&\omega^2\end{pmatrix}
,\ 
B = \begin{pmatrix}-1&0&0&0\\0&1&0&0\\0&0&-1&0\\0&0&0&1\end{pmatrix},
$$
$$
C = \begin{pmatrix}0&1&0&0\\1&0&0&0\\0&0&0&1\\0&0&1&0\end{pmatrix}
,\ 
D = \frac{1}{2\sqrt{3}}\begin{pmatrix}1&\sqrt{3}&\sqrt{2}&-\sqrt{6}\\\sqrt{3}&-1&-\sqrt{6}&-\sqrt{2}\\\sqrt{2}&-\sqrt{6}&2&0\\-\sqrt{6}&-\sqrt{2}&0&-2\end{pmatrix}
$$
where $\omega = e^\frac{2\pi i}{3}$, generate $\mathfrak{S}_{5}$ in $\textnormal{PSl}(4,\CC)$ with $A\mapsto(2,4,5)$, $B\mapsto(1,3)$, $C\mapsto(1,3)(4,5)$ and $D\mapsto(1,2)(4,5)$. This contains the subgroups $\left<A,C,D\right>\cong{\mathfrak{A}_{5}}.$ and $\left<(2,4,3,5), (1,4,2,3,5)\right>\cong\textnormal{Hol}(C_5)$.

This projective representation of $\mathfrak{S}_5$ is conjugate to the group labeled $I_*$ in \cite{HH}. Conjugate generators have been chosen to minimise the length of the invariant polynomials.
\begin{subEg}{$\mathfrak{A}_5$}
The $\mathfrak{A}_{5}$ invariant polynomials of degree 4 are
$$p := 5x_0^4+ 5x_1^4 - 6x_0^2x_1^2 + 12x_2^2x_3^2 - 8\sqrt{2}(x_0x_2^3  + x_1x_3^3) - 48x_0x_1x_2x_3,$$
and
$$q := (2x_0x_1 + 3x_2x_3)(x_0^2 + x_1^2) - 2\sqrt{2}(x_0x_3^3 +x_1x_2^3)$$
giving the invariant family of K3 surfaces
$$X_\lambda : p+\lambda{q}=0$$
in $\PP^3$. This family is singular at $\lambda^2 = -80$ and $\lambda^2 = 1$ and the symplectic automorphism group jumps to $\mathfrak{S}_{5}$ at $\lambda = 0$. The projective transformation
$$(x_0,x_1,x_2,x_3)\mapsto(-x_0,x_1,-x_2,x_3)$$
provides an isomorphism $X_\lambda \cong X_{-\lambda}$ so that $\lambda^2$ is the natural parameter for this family.
\end{subEg}
\begin{subEg}{$\textnormal{Hol}(C_5)$}
This subgroup is conjugate to the group generated by the matrices
$$
\begin{pmatrix}
0&0&-1&0\\
1&0&0&0\\
0&0&0&1\\
0&1&0&0
\end{pmatrix},\ \ \ \ \ 
\begin{pmatrix}
\zeta&0&0&0\\
0&\zeta^2&0&0\\
0&0&\zeta^3&0\\
0&0&0&\zeta^4
\end{pmatrix}
$$
where $\zeta^5=1$. The $\textnormal{Hol}(C_5)$ invariant polynomials of degree 4 are $
x_0^3x_1 + x_1^3x_3 + x_3^3x_2 - x_2^3x_0$ and $x_0^2x_3^2 + x_1^2x_2^2$ giving the invariant family of K3 surfaces
$$
X_\lambda := x_0^3x_1 + x_1^3x_3 + x_3^3x_2 - x_2^3x_0 + \left(\frac{1+i}{2}\right)\lambda(x_0^2x_3^2 + x_1^2x_2^2)=0
$$
This family is singular at $\infty$ and where $\lambda^4=1$. The projective transformation
$$
(x_0,x_1,x_2,x_3)\mapsto (-x_0,ix_1,-ix_2,x_3)
$$
provides an isomorphism $X_\lambda\cong X_{i\lambda}$ so that $\lambda^4$ is the natural parameter for the family.
\end{subEg}
\end{Eg}

\begin{Eg}{$L_2(7)$}
$L_2(7) = \textnormal{SmallGroup}(168,42)$ in Magma.

\begin{center}
\begin{picture}(50,80)(0,-30)

  \put(25,40){\makebox(0,0){$L_2(7)$}}
  \put(25,10){\makebox(0,0){$C_7\!\rtimes\!{C_3}$}}
  \put(25,-20){\makebox(0,0){$C_7$}}
  \put(25,16){\line(0,1){18}}
  \put(25,-12){\line(0,1){17}}

\end{picture}\end{center}

If $\rho:L_2(7)\rightarrow\textnormal{Sl}(3,\CC)$ is the representation leading to the example of symmetric surfaces in $\PP(1,1,1,3)$, then we obtain a reducible representation
\begin{eqnarray*}
\rho^\prime:L_2(7) & \rightarrow & \textnormal{Sl}(4,\CC)\\
g & \mapsto & \begin{pmatrix} \rho(g) & 0 \\ 0 & 1 \end{pmatrix}.
\end{eqnarray*}
\begin{subEg}{$C_7\rtimes{C_3}$}
Under this representation, the subgroup $C_7\rtimes{C_3}$ has degree 4 invariants $x_0x_1^3 + x_1x_2^3 + x_2x_0^3,\ x_3^4$, and $x_0x_1x_2x_3$ and an invariant family of K3 surfaces
$$
X_\lambda : x_0x_1^3 + x_1x_2^3 + x_2x_0^3 + x_3^4 + 4 \lambda x_0x_1x_2x_3 = 0
$$
in $\PP^3$. This is singular at $\lambda = \infty$ and where $\lambda^4 = 1$ and the symplectic automorphism group jumps to $L_2(7)$ at $\lambda = 0$. The projective transformation
$$(x_0,x_1,x_2,x_3)\mapsto(x_0,x_1,x_2,ix_3)$$
provides an isomorphism $X_\lambda \cong X_{i\lambda}$ so that $\lambda^4$ is the natural parameter for this family.
\end{subEg}
\end{Eg}

\begin{Eg}{$T_{192}$}
$T_{192} = \textnormal{SmallGroup}(192,1493)$ in Magma.

\begin{center}
\begin{picture}(60,80)(-10,-30)

  \put(25,40){\makebox(0,0){$T_{192}$}}
  \put(10,10){\makebox(0,0){$\Gamma_{\!\!25}a_1$}}
  \put(40,10){\makebox(0,0){$C_2\!\!\times\!\!\mathfrak{S}_4$}}
  \put(17,17){\line(1,3){5}}
  \put(33,17){\line(-1,3){5}}

  \put(10,4){\line(0,-1){18}}
  \put(40,4){\line(0,-1){18}}
  \put(10,-20){\makebox(0,0){$\Gamma_{\!7}a_1$}}
  \put(40,-20){\makebox(0,0){$C_2\!\!\times\!\!\mathfrak{A}_4$}}
\end{picture}
\end{center}

If we define the matrices
\begin{eqnarray*}
I & = & \begin{pmatrix}0&1\\-1&0\end{pmatrix}\\
J & = &\begin{pmatrix}i&0\\0&-i\end{pmatrix}
\end{eqnarray*}
then these generate the quaternion group $Q_8\subset\textnormal{Sl}(2,\CC)$ and the matrices
$$
I_l =  \begin{pmatrix}I&0\\0&id\end{pmatrix},
\ 
J_l =  \begin{pmatrix}J&0\\0&id\end{pmatrix},
\ 
I_r =  \begin{pmatrix}id&0\\0&I\end{pmatrix},
\ 
J_r =  \begin{pmatrix}id&0\\0&J\end{pmatrix},
$$
$$
A = \begin{pmatrix}0&0&1&0\\0&0&0&1\\1&0&0&0\\0&1&0&0\end{pmatrix},
\ 
B = \frac{1+i}{2}\begin{pmatrix}\omega&\omega&0&0\\i\omega&-i\omega&0&0\\0&0&-i\omega^2&-\omega^2\\0&0&-i\omega^2&\omega^2\end{pmatrix}
$$
with $\omega=e^\frac{2\pi i}{3}$, generate $T_{192}$ in $\textnormal{PSl}(4,\CC)$ containing the subgroups $\left<I_l,J_l,I_r,J_r,A\right>\cong{\Gamma_{25}a_1}$ and $\left<I_lJ_r,A,B\right>\cong{C_2\times{\mathfrak{S}_4}}$. This representation of $T_{192}$ also occurrs in \cite{Mukai}.

\begin{subEg}{$\Gamma_{25}a_1$}\label{Gamma25a1}
The $\Gamma_{25}a_1$ invariant polynomials of degree 4 are  
$x_0^4+x_1^4+x_2^4+x_3^4$, and $x_0^2x_1^2 + x_2^4x_3^2$ giving the invariant family of K3 surfaces
$$X_\lambda : x_0^4+x_1^4+x_2^4+x_3^4 + 2 \lambda (x_0^2x_1^2 + x_2^2x_3^2)=0$$
in $\PP^3$. This family is singular at $\lambda^2 = 1$ and the symplectic automorphism group jumps to $T_{192}$ at $\lambda = \pm\sqrt{3}i$. The projective transformation
$$(x_0,x_1,x_2,x_3)\mapsto(ix_0,x_1,ix_2,x_3)$$
provides an isomorphism $X_\lambda \cong X_{-\lambda}$ so that $\lambda^2$ is the natural parameter for this family.
\end{subEg}
\begin{subEg}{$C_2\times{\mathfrak{S}_4}$}\label{C2A4}
The $C_2\times{\mathfrak{S}_4}$ invariant polynomials of degree 4 are
\begin{eqnarray*}
p & = & x_0^4 + x_1^4 + x_2^4 + x_3^4 - 2i\sqrt{3}(x_0^2x_1^2 + x_2^2x_3^2)\\
q & = & (x_0x_2+x_1x_3-i x_1x_2-i x_0x_3)^2
\end{eqnarray*}

giving the invariant family of K3 surfaces
$$X_\lambda : p + 2\omega\lambda q$$
in $\PP^3$ where $\omega = \frac{\sqrt{3}+i}{2}$ satisfies $\omega^{12}=1$. This family is singular where $\lambda^2 = \frac{1}{4}$ and $\lambda^2 = \frac{1}{3}$ and the symplectic automorphism group jumps to $T_{192}$ at $\lambda = 0$. The projective transformation
$$(x_0,x_1,x_2,x_3)\mapsto(ix_0,ix_1,x_2,x_3)$$
provides an isomorphism $X_\lambda \cong X_{-\lambda}$ so that $\lambda^2$ is the natural parameter for this family.
\end{subEg}
\end{Eg}
\begin{Eg}{$F_{384}$}
$F_{384} = \textnormal{SmallGroup}(384,18135)$ in Magma.

\begin{center}
\begin{picture}(150,150)(-40,-150)

  \put(45,-10){\makebox(0,0){$F_{384}$}}

  \put(47,-16){\line(2,-3){10}}
  \put(52,-16){\line(2,-1){30}}
  \put(43,-16){\line(-2,-3){10}}
  \put(38,-16){\line(-2,-1){30}}

  \put(30,-40){\makebox(0,0){$2^4\!D_6$}}
  \put(60,-40){\makebox(0,0){$\Gamma_{\!\!25}a_1$}}
  \put(60,-46){\line(0,-1){18}}
  \put(60,-70){\makebox(0,0){$\Gamma_{\!7}a_1$}}

  \put(90,-40){\makebox(0,0){$SD_{16}$}}
  \put(90,-46){\line(0,-1){18}}
  \put(90,-70){\makebox(0,0){$C_8$}}

  \put(0,-40){\makebox(0,0){$4^2\!\mathfrak{A}_4$}}
  \put(-15,-70){\makebox(0,0){$\Gamma_{\!\!13}a_1$}}
  \put(15,-70){\makebox(0,0){$4^2C_3$}}
  \put(4,-46){\line(1,-3){5}}
  \put(-4,-46){\line(-1,-3){5}}
  \put(-15,-76){\line(0,-1){18}}
  \put(-15,-100){\makebox(0,0){$\Gamma_4c_2$}}

	\put(0,-121){\line(1,4){11}}
	
  \put(-19,-106){\line(-1,-3){5}}
  \put(-11,-106){\line(1,-3){5}}
  \put(-30,-130){\makebox(0,0){$C_2\!\!\times\!\!{Q_8}$}}
  \put(0,-130){\makebox(0,0){$C_4^2$}}

\end{picture}\end{center}
The matrices
$$
A =  \begin{pmatrix}i&0&0&0\\0&-i&0&0\\0&0&1&0\\0&0&0&1 \end{pmatrix}
,\ 
B = \begin{pmatrix}1&0&0&0\\0&i&0&0\\0&0&-i&0\\0&0&0&1 \end{pmatrix}
,\ 
C = \begin{pmatrix}1&0&0&0\\0&0&1&0\\0&0&0&1\\0&1&0&0\end{pmatrix},
$$
$$
D = \begin{pmatrix}0&1&0&0\\1&0&0&0\\0&0&1&0\\0&0&0&-1\end{pmatrix}
,\ 
E = \begin{pmatrix}0&1&0&0\\1&0&0&0\\0&0&0&1\\0&0&1&0\end{pmatrix}
$$
generate $F_{384}$ in $\textnormal{PSl}(4,\CC)$ containing the subgroups $\left<A,B,C,E\right>\cong{4^2\mathfrak{A}_4}$,\\
$\left<D,EB\right>\cong{SD_{16}}$, $\left<C,D,E\right>\cong{2^4D_6}$ and $\left<A,C^{-1}BC,D,C^{-1}EC\right>\cong{\Gamma_{25}a_1}$. This projective representation of $F_{384}$ is from \cite{Mukai}.

\begin{subEg}{$4^2\mathfrak{A}_4$}\label{42A4}
The $4^2\mathfrak{A}_4$ invariant polynomials of degree 4 are  $x_0^4+x_1^4+x_2^4+x_3^4,\ x_0x_1x_2x_3$ giving the invariant family of K3 surfaces
$$X_\lambda :  x_0^4+x_1^4+x_2^4+x_3^4 + 4\lambda x_0x_1x_2x_3=0$$
in $\PP^3$. This family is singular at $\lambda = \infty$ and $\lambda^4 = 1$ and the symplectic automorphism group jumps to $F_{384}$ at $\lambda = 0$. The projective transformation
$$(x_0,x_1,x_2,x_3)\mapsto(ix_0,x_1,x_2,x_3)$$
provides an isomorphism $X_\lambda \cong X_{i\lambda}$ so that $\lambda^4$ is the natural parameter for this family.
\end{subEg}
\begin{subEg}{$SD_{16}$}
The $SD_{16}$ invariant polynomials of degree 4 are  $x_0^4+x_1^4+x_2^4+x_3^4,\ x_0x_1(x_2^2+ix_3^2)$ giving the invariant family of K3 surfaces
$$X_\lambda :  x_0^4+x_1^4+x_2^4+x_3^4 + 2\sqrt{2}\lambda x_0x_1(x_2^2+ix_3^2)=0$$
in $\PP^3$. This family is singular at $\lambda = \infty$ and $\lambda^4 = 1$ and the symplectic automorphism group jumps to $F_{384}$ at $\lambda = 0$. The projective transformation
$$(x_0,x_1,x_2,x_3)\mapsto(ix_0,x_1,x_2,x_3)$$
provides an isomorphism $X_\lambda \cong X_{i\lambda}$ so that $\lambda^4$ is the natural parameter for this family.
\end{subEg}
\begin{subEg}{$2^4D_6$}
The $2^4D_6$ invariant polynomials of degree 4 are  $x_0^4+x_1^4+x_2^4+x_3^4,\ (x_0^2+x_1^2+x_2^2+x_3^2)^2$ giving the invariant family of K3 surfaces
$$X_\lambda :  x_0^4+x_1^4+x_2^4+x_3^4 + \lambda(x_0^2+x_1^2+x_2^2+x_3^2)^2=0$$
in $\PP^3$. This family is singular at $\lambda = \infty$ and $\lambda = -1,-\frac{1}{2},-\frac{1}{3},-\frac{1}{4}$ and the symplectic automorphism group jumps to $F_{384}$ at $\lambda = 0$ and to $M_{20}$ at $\lambda = \frac{3}{4}$.
\end{subEg}
\begin{subEg}{$\Gamma_{25}a_1$}
This subgroup coincides exactly with the subgroup $\Gamma_{25}a_1\subset T_{192}$ considered above.
\end{subEg}
\end{Eg}

\begin{Eg}{$M_{20}$}
$M_{20} = \textnormal{SmallGroup}(960,11357)$ in Magma.

\begin{center}
\begin{picture}(80,170)(-40,-160)

  \put(0,-10){\makebox(0,0){$M_{20}$}}
  \put(-30,-40){\makebox(0,0){$\mathfrak{A}_5$}}
  \put(0,-40){\makebox(0,0){$4^2\mathfrak{A}_4$}}
  \put(30,-40){\makebox(0,0){$2^4D_6$}}
  \put(0,-16){\line(0,-1){17}}
  \put(-5,-16){\line(-1,-1){17}}
  \put(5,-16){\line(1,-1){17}}
  \put(-15,-70){\makebox(0,0){$\Gamma_{\!\!13}a_1$}}
  \put(15,-70){\makebox(0,0){$4^2C_3$}}
  \put(4,-46){\line(1,-3){5}}
  \put(-4,-46){\line(-1,-3){5}}
  \put(-15,-76){\line(0,-1){18}}
  \put(-15,-100){\makebox(0,0){$\Gamma_4c_2$}}

  \put(-19,-106){\line(-1,-3){5}}
  \put(-11,-106){\line(1,-3){5}}
  \put(0,-130){\makebox(0,0){$C_2\!\!\times\!\!{Q_8}$}}
  \put(-30,-130){\makebox(0,0){$C_4^2$}}

\end{picture}
\end{center}
The matrices
$$
A =  \begin{pmatrix}i&0&0&0\\0&-i&0&0\\0&0&1&0\\0&0&0&1 \end{pmatrix}
,\ 
B = \begin{pmatrix}1&0&0&0\\0&i&0&0\\0&0&-i&0\\0&0&0&1 \end{pmatrix}
,\ 
C = \begin{pmatrix}1&0&0&0\\0&0&1&0\\0&0&0&1\\0&1&0&0\end{pmatrix},
$$
$$
E = \begin{pmatrix}0&1&0&0\\1&0&0&0\\0&0&0&1\\0&0&1&0\end{pmatrix}
,\ 
F = \frac{1}{2}\begin{pmatrix}-i&-i&-1&\phantom{-}1\\-i&\phantom{-}i&-1&-1\\\phantom{-}i&\phantom{-}i&-1&\phantom{-}1\\\phantom{-}i&-i&-1&-1\end{pmatrix}
$$
generate $M_{20}$ in $\textnormal{PSl}(4,\CC)$ containing the subgroups $\left<B,E,AFB\right>\cong{2^4D_6}$,
$\left<A,B,C,E\right>\cong{4^2\mathfrak{A}_4}$ and $\left<A^{-1}CA,F\right>\cong{\mathfrak{A}_5}$
\begin{subEg}{$2^4D_6$}
The $2^4D_6$ invariant polynomials of degree 4 are  $x_0^4+x_1^4+x_2^4+x_3^4+12x_0x_1x_2x_3$ and $x_0^2x_3^2 +  x_1^2x_2^2 + 2x_0x_1x_2x_3$ giving the invariant family of K3 surfaces
$$X_\lambda :  x_0^4+x_1^4+x_2^4+x_3^4 + 12x_0x_1x_2x_3+4\lambda (x_0^2x_3^2 +  x_1^2x_2^2 + 2x_0x_1x_2x_3)=0$$
in $\PP^3$. This family is singular at $\lambda = \infty, -1, -\frac{5}{6}, -\frac{1}{2}, \frac{1}{2}$ and the symplectic automorphism group jumps to $M_{20}$ at $\lambda = 0$ and to $F_{384}$ at $\lambda = -\frac{3}{2}$. The projective transformation
$$
\frac{1}{\sqrt{2}}\begin{pmatrix}
1&0&0&i\\
0&1&i&0\\
0&i&1&0\\
i&0&0&1
\end{pmatrix}
$$
transforms this family to the family $x_0^4+x_1^4+x_2^4+x_3^4 + \lambda(x_0^2+x_1^2+x_2^2+x_3^2)^2=0$ considered earlier.
\end{subEg}
\begin{subEg}{$4^2\mathfrak{A}_4$}
This subgroup coincides exactly with the subgroup $4^2\mathfrak{A}_4\subset F_{384}$ considered earlier.
\end{subEg}
\begin{subEg}{$\mathfrak{A}_5$}
The nondegenerate quadric
$$
q := x_0^2 - x_1^2 + x_2^2 + x_3^2 + (i + 1)(- x_0x_1 + ix_0x_2 + ix_0x_3 + x_1x_2 + x_1x_3 - ix_2x_3)
$$
is invariant under this action of $\mathfrak{A}_5$. The degree 4 invariants are
$$
p := x_0^4 + x_1^4 + x_2^4 + x_3^4 + 12x_0x_1x_2x_3
$$
and
$$
q^2.
$$
We choose to define the family of quartic K3 surfaces
$$
X_\lambda :  \left( ((5+15i)p+(1-3i)q^2) + \lambda q^2 = 0 \right) \subset \PP^3.
$$
This family degenerates at $\infty$, $-\frac{16}{3}, -5, 0$ and $3$. The strange choice of parameter is made so that the degenerate points lie on the real line.
\end{subEg}
\end{Eg}

To summarise, we have found 9 distinct families of K3 hypersurfaces in $\PP^3$ invariant under the groups $C_7\rtimes{C_3}$, $\Gamma_{25}a_1$, $C_2\times\mathfrak{S}_4$, $4^2\mathfrak{A}_4$, $2^4D_6$, $SD_{16}$, $\mathfrak{A}_5$ (2 examples) and $\textnormal{Hol}(C_5)$. These are listed in table \ref{tab:K3HypersurfacesInPP3}. In all these examples, the representation is induced from the representation of one of the 11 maximal groups.

\begin{table*}
	\newcommand\T{\rule{0pt}{2.6ex}}
  \newcommand\B{\rule[-1.2ex]{0pt}{0pt}}
  \centering
		\begin{tabular}[h]{|c|c|c|}\hline
			\T\B & Symmetry Group & Equation\\
			\hline\hline
			\T\B VI&$4^2\mathfrak{A}_4\subset F_{384}, M_{20}$ & $x_0^4+x_1^4+x_2^4+x_3^4 + 4\lambda x_0x_1x_2x_3=0$\\\hline
			\T\B VII & $C_7\rtimes{C_3}\subset L_2(7)$ & $x_0x_1^3 + x_1x_2^3 + x_2x_0^3 + x_3^4 + 4 \lambda x_0x_1x_2x_3 = 0$\\\hline
			\T\B VIII & $SD_{16}\subset F_{384}$ & $x_0^4+x_1^4+x_2^4+x_3^4 + 2\sqrt{2}\lambda x_0x_1(x_2^2+ix_3^2)=0$\\\hline
			\T\B IX & $\textnormal{Hol}(C_5)\subset \mathfrak{S}_5$ & $x_0^3x_1 + x_1^3x_3 + x_3^3x_2 - x_2^3x_0 + \left(\frac{1+i}{2}\right)\lambda(x_0^2x_3^2 + x_1^2x_2^2)=0$\\\hline
			\T\B X & $\Gamma_{25}a_1\subset T_{192}, F_{384}$ & $x_0^4+x_1^4+x_2^4+x_3^4 + 2 \lambda (x_0^2x_1^2 + x_2^2x_3^2)=0$\\\hline
			\T\B XI & $2^4D_6\subset F_{384}, M_{20}$ & $x_0^4+x_1^4+x_2^4+x_3^4 + \lambda(x_0^2+x_1^2+x_2^2+x_3^2)^2=0$\\\hline
			\T\B XII & $C_2\times\mathfrak{S}_4\subset T_{192}$ & $p+\lambda q^2=0$\\
			\T\B && $p = x_0^4 + x_1^4 + x_2^4 + x_3^4 - 2i\sqrt{3}(x_0^2x_1^2 + x_2^2x_3^2)$\\
			\T\B && $q  =  x_0x_2+x_1x_3-i x_1x_2-i x_0x_3$\\\hline
			\T\B XIII & $\mathfrak{A}_5\subset M_{20}$ & $(5+15i)p+(\lambda+ 1 - 3i)q^2 = 0$\\
			\T\B && $p=x_0^4 + x_1^4 + x_2^4 + x_3^4 + 12x_0x_1x_2x_3$\\
			\T\B &&$q=x_0^2 - x_1^2 + x_2^2 + x_3^2$\\
			\T\B &&$+ (i + 1)(- x_0x_1 + ix_0x_2 + ix_0x_3 + x_1x_2 + x_1x_3 - ix_2x_3)$\\\hline
			\T\B XIV & $\mathfrak{A}_5\subset \mathfrak{S}_5$ & $p+\lambda q=0$\\
			\T\B &&$p := 5x_0^4+ 5x_1^4 - 6x_0^2x_1^2 + 12x_2^2x_3^2$\\
			\T\B &&$ - 8\sqrt{2}(x_0x_2^3  + x_1x_3^3) - 48x_0x_1x_2x_3$\\
			\T\B &&$q := (2x_0x_1 + 3x_2x_3)(x_0^2 + x_1^2) - 2\sqrt{2}(x_0x_3^3 +x_1x_2^3)$\\\hline
			\end{tabular}
	\caption{K3 Hypersurfaces in $\PP^3$}
	\label{tab:K3HypersurfacesInPP3}
\end{table*}

\section{Moduli of Lattice Polarised K3 Surfaces}\label{moduli}
In this section we look at a method to construct a coarse moduli space for lattice polarised K3 surfaces described in \cite{Dol}. We use this construction to derive theorem \ref{squares} that will be of great use to us later on. Recall the definition

\begin{definition}[\cite{Dol}]
Let $S$ be a lattice of signature $(1,k-1)$. An \emph{$S$--polarisation} of a K3 surface $X$ is a primitive embedding $\iota:S\hookrightarrow \textnormal{Pic}(X)$.
\end{definition}

Our interest in lattice polarisations stems from the fact that all of our families of K3 hypersurfaces in $\PP^3$ or $\PP(1,1,1,3)$ are lattice polarised. For a hypersurface, the K\"ahler class $[\kappa]$ coincides with the class of a hyperplane divisor. Our group actions are induced from linear automorphisms of the ambient projective spaces, and so we find $[\kappa]\in\textnormal{H}^2(X,\ZZ)^G$. With the lattice $S_{X,G} = \textnormal{H}^2(X,\ZZ)^{G\perp}$ as defined at the start of this chapter, the sublattice $S_{X,G}\oplus [\kappa]$ defines a lattice polarisation of the family.

A course moduli space for $S$-polarised K3 surfaces is obtained by first constructing the fine moduli space of \emph{marked} $S$-polarised K3s and then removing the choice of marking. A marking is just a choice of isomorphism $H^2(X,\ZZ)\stackrel{\phi}{\longleftarrow}(-E_8)^2\oplus U^3 =: \Lambda_{K3}$ and so after fixing a sublattice $S\subset\Lambda_{K3}$, we obtain a lattice polarisation of $X$ whenever $\phi(S)\subset\textnormal{Pic}(X)$. The period space of marked $S$-polarised K3 surfaces is the space
$$
\Omega_S := \{ \omega \in \PP(\Lambda_{K3}\otimes{\CC})\ |\ \left<\omega, \omega\right> = 0,\ \left<\omega, \bar{\omega}\right> > 0,\ \omega \in S^\perp \}.
$$
 The three conditions mimic the properties $\omega\wedge\omega=0$, $\omega\wedge\bar{\omega}>0$ satisfied by any $\omega\in\textnormal{H}^{2,0}(X)$ and the fact that $S\subset\textnormal{Pic}(X) = \textnormal{H}^{1,1}(X)\cap\textnormal{H}^2(X,\ZZ)$. This period space is a fine moduli space for marked $S$--polarised K3 surfaces. The lattice automorphisms
$$
\Gamma(S) := \{g\in \textnormal{O}(\Lambda_{K3})\ |\ g(s)=s\ \textnormal{for all}\ s \in S\}
$$
account for all the equivalent markings of the surface. Hence, the quotient
$$
\mathcal{K}_S := \Omega_S/\Gamma(S)
$$
is a (coarse) moduli space of $S$--polarised K3 surfaces.

With \cite{Dol}, we write $T = S^\perp\cap\textnormal{H}^2(X,\ZZ)$, and obtain a simplified description of $\mathcal{K}_S$ by defining
$$
D_S = \{ \omega \in \PP(T\otimes{\CC})\ |\ \left<\omega, \omega\right> = 0,\ \left<\omega, \bar{\omega}\right> > 0\}
$$
so that
$$
\mathcal{K}_S = D_S/\Gamma_S
$$
where $\Gamma_S\subset\textnormal{O}(T)$ are the automorphisms of $T$ induced by automorphisms of $\Lambda_{K3}$. in order to construct this quotient in any given example, it will be necessary to better understand the group $\Gamma_S$. Fortunately, we have the following fact.
\begin{proposition}\cite{Dol}
A lattice automorphism $g \in \textnormal{O}(T)$ induces a natural automorphism $\theta(g)$ of the discriminant group $T^*/T$ defined by
$$
v+T\mapsto v\circ g+T.
$$
The group $\Gamma_S$ consists precisely of those automorphisms that induce the identity on $T^*/T$. In other words,
$$
\Gamma_S = \textnormal{Ker}(\theta)
$$
and in particular, $\Gamma_S$ has finite index in the arithmetic group $\textnormal{O}(T)$ and so is itself arithmetic.
\end{proposition}

We now restrict our attention to the case where the lattice $S$ has rank 19 and we diverge from the treatment of \cite{Dol}. In this case, $T=S^\perp$ has rank 3 and signature $(2,1)$ and, $T\otimes\RR = R\otimes\RR$, where $R$ is any other lattice with the same signature. In particular, for reference we fix
$$
R=\begin{pmatrix}
0&0&-1\\
0&1&0\\
-1&0&0\end{pmatrix}
$$
and consider the isomorphism
$$
\varphi:\PP\textnormal{Sl}(2,\RR)\rightarrow \textnormal{SO}(R)\otimes\RR
$$
defined by
$$
\begin{pmatrix}
a&b\\c&d
\end{pmatrix}
\mapsto
\begin{pmatrix}
a^2&\sqrt{2}ab&b^2\\
\sqrt{2}ac&ad+bc&\sqrt{2}bd\\
c^2&\sqrt{2}cd&d^2
\end{pmatrix}.
$$

If $P\in\textnormal{Gl}(3,\RR)$ is a matrix satisfying
$$^tP.T.P=R$$
and $A\in\textnormal{SO}(R)\otimes\RR$, then $P.A.P^{-1}\in\textnormal{SO}(T)\otimes\RR$ and we obtain an isomorphism
\begin{equation}\label{quad}
\Phi := P\varphi{P^{-1}}:\PP\textnormal{Sl}(2,\RR)\rightarrow \textnormal{SO}(T)\otimes\RR
\end{equation}
defined by
$$
\begin{pmatrix}
a&b\\c&d
\end{pmatrix}
\mapsto
P.\begin{pmatrix}
a^2&\sqrt{2}ab&b^2\\
\sqrt{2}ac&ad+bc&\sqrt{2}bd\\
c^2&\sqrt{2}cd&d^2
\end{pmatrix}.P^{-1}.
$$
Under this second isomorphism, $\textnormal{SO}(T)$ pulls back to an arithmetic Fuchsian group containing $\Gamma_S$ as a finite index subgroup. This pullback allows us to construct the coarse moduli space of $S$--polarised K3 surfaces as a Shimura curve
$$\mathcal{K}_S = \mathcal{H}/\Phi^{-1}(\Gamma_S).$$

Here, $\mathcal{H}$ denotes the upper--half plane on which the Fuchsian group naturally acts. Geometrically, the identification $\mathcal{K}_S = \mathcal{H}/\Phi^{-1}(\Gamma_S)$ occurs since $D_S\subset\PP(T\otimes\CC)\cong\PP^2$ is defined by the two equations $\left<\omega,\omega\right>=0$ and $\left<\omega,\bar{\omega}\right> > 0$. The first equation cuts out the non--degenerate quadric $\cong\PP^1$ corresponding to the quadratic form defined by the lattice $T$. The second condition restricts to the union of two open half--planes in this quadric. As $T$ has odd rank, $\textnormal{SO}(T)=\textnormal{O}(T)/(-id)$. But $-id$ simply swaps the two copies of the half--plane that make up the period space and we may restrict to the action of $\Phi^{-1}(\Gamma_S)\subset\Phi^{-1}(\textnormal{SO}(T))$ on one of these copies of $\mathcal{H}$.

\begin{theorem}\label{squares}
If $\Gamma\subset\PP\textnormal{Sl}(2,\RR)$ is the monodromy group of a family of rank 19 lattice polarised K3 surfaces, then
$$
\textnormal{trace}(A)^2\in\ZZ\ \ \ \ \ \textnormal{for all}\ \ \ \ \ A\in\Gamma.
$$
\end{theorem}
\begin{proof}
Up to conjugation, we may express $\Gamma = \Phi^{-1}(\Gamma_S) \subset \Phi^{-1}(\textnormal{SO}(T))$ where $\Phi$ is conjugate to the reference isomorphism $\varphi$ defined above. Since trace is invariant under conjugation, the theorem follows from the observation that 
$$
\textnormal{trace}(\varphi(A)) = a^2 + bc + ad + d^2 = (a+d)^2 - ad + bc = \textnormal{trace}(A)^2-1
$$
and that $\varphi(A)\in\textnormal{SO}(T)$ has integer entries.
\end{proof}

\section{Lattice Polarisations of Symmetric Surfaces}

Recall that if $X$ is an algebraic K3 surface with a finite group of symplectic automorphisms, $G$, then $X$ is polarised by the lattice
\begin{equation}\label{kappa}
S = S_{X,G}\oplus\kappa
\end{equation}
where $S_{X,G} = \textnormal{H}^2(X,\ZZ)^{G\perp}$ and $\kappa$ is a primitive ample divisor.

The minimal resolution of $X/G$ is a K3 surface $\widetilde{X/G}$ fitting into the diagram
\begin{eqnarray*}
&&X\\
&&\:\downarrow\pi\\
\widetilde{X/G}&\stackrel{\sigma}{\longrightarrow}&X/G.
\end{eqnarray*}

This leads to a generically finite rational morphism $\tau:\sigma^{-1}\circ\pi:X\rightarrow\widetilde{X/G}$. Minimally eliminating the indeterminacies of $\tau$, there is a surface $\widetilde{X}$ with a blow-down map $\rho:\widetilde{X}\rightarrow X$ and a generically $|G|$ to 1 cover $\widetilde{\pi}:\widetilde{X}\rightarrow\widetilde{X/G}$ ramified over the exceptional curves of $\sigma:\widetilde{X/G}\rightarrow{X/G}$. Hence, we have a commutative diagram
\begin{eqnarray*}
\,\,\widetilde{X}&\stackrel{\rho}{\longrightarrow}&X\\
\downarrow\widetilde{\pi}&&\:\downarrow\pi\\
\widetilde{X/G}&\stackrel{\sigma}{\longrightarrow}&X/G.
\end{eqnarray*}
We define
$$
E\subset \textnormal{H}^2(\widetilde{X/G},\ZZ)
$$
as the lattice spanned by the exceptional curves of $\sigma$. The lattice $E$ is not necessarily primitively embedded in $\textnormal{Pic}(\widetilde{X/G})$. However, there will be a smallest lattice $E^\prime$ containing $E$ that is primitive in $\textnormal{Pic}(\widetilde{X/G})$. Lemma 2 of \cite{Xiao} shows that
$$
E^\prime/E \cong G/[G,G].
$$
For example, if $G$ is a cyclic group, then $E^\prime/E \cong G$. The $|G|$ to 1 cyclic cover $\widetilde{\pi}:\widetilde{X}\rightarrow\widetilde{X/G}$ will have a ramification divisor $R\in E$ with
$$\frac{1}{|G|}R\in\textnormal{Pic}(\widetilde{X/G})\setminus E.$$
See \cite{BPV} section I.16. In this case, $E$ and $\frac{1}{|G|}R$ span the primitive sublattice $E^\prime$ in $\textnormal{Pic}(\widetilde{X/G})$.

\begin{example}
For the group $G=\mathfrak{A}_5$ (see families V, XIII and XIV) we notice that $[G,G] = G$ and so $E$ is primitively embedded in the Picard lattice. In this example, $E$ is the lattice
$$
E\cong 2A_4\oplus 3A_2\oplus 4A_1
$$
with discriminant $5^2.3^3.2^4.$ The quotient families are polarised by the lattice generated by $E$ and the image of a hyperplane.
\end{example}

\begin{example}
Family VII is defined in $\PP^3$ by the equation
$$
X_\lambda: (x_0x_1^3 + x_1x_2^3 + x_2x_0^3 + x_3^4 + 4 \lambda x_0x_1x_2x_3 = 0).
$$
Being invariant under the group $C_7\rtimes C_3$, it is also invariant under a subgroup of order 7 generated by the transformation
$$(x_0,x_1,x_2,x_3)\mapsto(\zeta x_0,\zeta^2 x_1, \zeta^4 x_2, x_3)$$
with $\zeta = e^\frac{2\pi i}{7}$ (see also section \ref{wps}).
This transformation has three fixed points on $X_\lambda$ at $(1,0,0,0)$, $(0,1,0,0)$ and $(0,0,1,0)$ providing the quotient with three $A_6$ singularities. The resolution of the three resulting $A_6$ singularities of $X/C_7$ forms the tetrahedral configuration of figure \ref{fig:C7Labelsx_i} together with 4 further curves shown as unfilled circles in the diagram. The three curves on $X_\lambda$ defined by $(x_i=0)$ for $i = 0,1$ and $2$ are nodal rational curves and the curve defined by $(x_3=0)$ is smooth of genus 3. These four curves correspond to curves labeled $x_i$ on $\widetilde{X_\lambda/C_7}$ which can be shown to be rational curves by the Riemann Hurwitz formula.
\begin{figure}[h]
	\centering
		\includegraphics{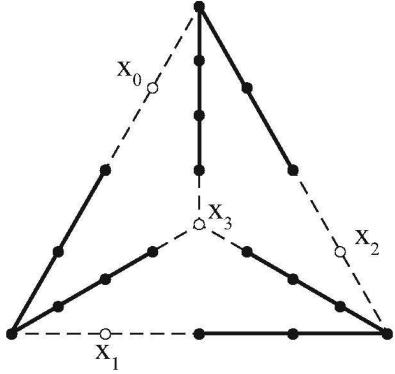}
	\caption{Exceptional Curves on $X_\lambda/C_7$.}
	\label{fig:C7Labelsx_i}
\end{figure}

Because of the ramified cyclic covering $\widetilde{\pi}:\widetilde{X_\lambda}\rightarrow\widetilde{X_\lambda/C_7}$, the exceptional lattice $E = 3A_6$ is not primitive and contains an element $R$ divisible by 7 in the Picard lattice.

Defining $M_i := E\oplus[x_i]$ for $i = 0,1,2$ and $3$, we find 
$$M_0 \cong M_1 \cong M_2 \cong A_6\oplus T_{3,4,8}$$
and
$$M_3 \cong T_{7,7,7}$$
where $T_{a,b,c}$ are the lattices from `T' shaped configurations of rational curves as featured in \cite{Bel}. Since $T_{a,b,c}$ has discriminant $abc-ab-bc-ca$, we see that for each $i$
$$\textnormal{disc}(M_i) = 4.7^2.$$
Choose one of these lattices and label is as $M$.

Let $R\in E$ denote the divisor with $\frac{1}{7}R\in\textnormal{Pic}(\widetilde{X_\lambda/C_7})\setminus E$. The overlattice $M^\prime$ spanned by $M$ and $\frac{1}{7}R$ will have discriminant
$$
\textnormal{disc}M^\prime = \frac{\textnormal{disc} M}{7^2} = 4.
$$
This rank 19 lattice will be the generic Picard lattice for this family (any lattice $N$ with $M^\prime\subset N\subset M^{\prime *}$ would have to be unimodular, even and of rank 19 -- an impossibility). 
\end{example}

We see from these examples that it will typically be possible to find the lattice that polarises the quotient surfaces in any particular case. However, in order to find the coarse moduli space of such surfaces using the method of \cite{Dol}, it is necessary to find the orthogonal complement to this polarising lattice. This is a difficult task when that polarising lattice has high discriminant.

	\newpage
	
\chapter{Picard-Fuchs Differential Equations}
\section{Differential Equations and Systems}
In this section, we draw together some techniques and theorems on ordinary differential equations that we shall require later. This material is not new and most of the definitions can be found in \cite{CL} and \cite{In} although the general approach we shall adopt is original.
\subsection{Fuchsian Differential Equations and Systems}

Let
$$
D(z)\colon \!\!\!\!=\frac{d^nf}{dz^n}+\sum_{k=1}^{n}a_{n-k}(z)\frac{d^{n-k}f}{dz^{n-k}}=0
$$
be an ordinary differential equation where the functions $a_i(z)$ are meromorphic on $\CC$. A singularity of $D$ is any $z_0\in\CC$ which is a singularity of some $a_i$. The point at infinity is singular when $D(1/y)$ is singular at $y=0$.

\begin{definition}
A singular point, $\alpha\in\CC$, of $D(z)$ is said to be a \emph{regular singular point} if
$$
\lim_{z\rightarrow\alpha}(z-\alpha)^ka_{n-k}(z)
$$
exists for $k=1,\ldots,n.$ Similarly, the point at infinity is a regular singular point if $D(1/y)$ is regular singular at $y=0$.

A differential equation whose singular points in $\PP^1$ are all regular singular is known as a \emph{Fuchsian differential equation}.
\end{definition}

In particular, since the coefficients of a Fuchsian differential equation are globally meromorphic, they must be rational functions. In this case, we are able to clear denominators and write the differential equation in the form
$$
\sum_{k=0}^{n}a_{k}(z)\frac{d^{k}f}{dz^{k}}=0
$$
where the $a_k(z)\in\CC[z]$ have no common factor.

We shall find it convenient to write any linear differential equation as a linear system of first order differential equations. For example, a second order differential equation
$$
a_2f^{\prime\prime}+a_1f^{\prime}+a_0f=0
$$
is directly equivalent to the system
\begin{align*}
f^{\prime} & =  g\\
g^{\prime} & =  -\frac{a_0}{a_2}\ f-\frac{a_1}{a_2}\ g
\end{align*}
which we may rewrite in the form
$$
\frac{d}{{dz}}
\begin{pmatrix}
    f  \\
    f^\prime \\
\end{pmatrix}
=
\begin{pmatrix}
   0 & 1  \\
   -\frac{a_0}{a_2} & -\frac{a_1}{a_2}  \\
\end{pmatrix}
\begin{pmatrix}
    f  \\
    f^\prime \\
\end{pmatrix}.
$$
We shall use the term \emph{differential system} to mean a first order matrix system of differential equations
\begin{equation}
\frac{d}{{dz}}\omega(z)= A(z)\omega(z)\label{sys}
\end{equation}
where $A(z)\in \textnormal{M}_{n}(\CC(z))$.

In agreement with the definition of regular singular point for a first order differential equation, we define
\begin{definition}
A singularity, $\alpha$, of the differential system (\ref{sys}) is said to be a \emph{regular singular point} if the matrix function $A(z)$ has a simple pole at $z=\alpha$. Also, a differential system is said to be \emph{Fuchsian} if its only singularities on $\PP^1$ are regular.
\end{definition}

As we have seen, for every differential equation, we can associate a differential system. However, if the initial differential equation is Fuchsian, the system as constructed as above will not be a Fuchsian system. We put this right later with theorem~\ref{equivSys}. The following proposition, often taken as a definition, gives an immediate way of recognising a Fuchsian system.

\begin{proposition}[\cite{CL}]
The differential system (\ref{sys}) is Fuchsian if and only if $A(z)$ has the form
\begin{eqnarray}\label{system}
A(z) = \sum_{i=1}^{k}\frac{R_i}{z-\alpha_i}
\end{eqnarray}
where $\alpha_i\in\CC$ are the finite regular singular points and $R_i\in\textnormal{M}_{n}(\CC)$ are constant matrices, called the residue matrices.

Furthermore, the system has a regular singular point at infinity with residue matrix $R_\infty = -\sum_{i=1}^{k}R_i$ (nonsingular if $\sum R_i = 0$).
\end{proposition}
\begin{proof}
The following proof is similar to that found on page 129 of \cite{CL} with the addition of the form of the residue matrix at infinity.\\
If the system $\frac{d}{dz}\omega(z) = A(z)\omega(z)$ is Fuchsian, then we may write
$$
A(z) = \sum_{i=1}^k\frac{R_i}{z-\alpha_i} + C(z)
$$
where $C(z)$ is regular at all $z\in\CC$ and possibly has a pole at infinity. Hence, $C(z)$ has polynomial entries and may be written as
$$
C(z) = \sum_{j=1}^mC_jz^j
$$
with $C_j\in\textnormal{M}_{n}(\CC)$.
Since the point at infinity must also be regular singular, setting $y=1/z$, we must have that
$$
-y^2 \frac{d}{dy}\omega  = A(y^{-1})\omega
$$
is regular singular at $y=0$. However, rearranging this gives
$$
\frac{d}{dy}\omega = -\sum_{i=1}^{k}\left(\frac{R_i}{y}+\frac{R_i}{y-1/\alpha_i}\right)-\sum_{j=1}^{m}\frac{C_j}{y^{j+2}}
$$
which can have a simple pole at $y=0$ only if $C_j=0$ for every $j$, in which case $C(z)=0$. Notice that the residue at $y=0$ is $\sum_{i=1}^{k}R_i$ as required.\\
The converse is proved similarly.
\end{proof}

\subsection{Local Solutions and Monodromy}
About any nonsingular point, $z_0$, of an n--th order differential equation, there exist $n$ linearly independent solutions, each analytic at $z_0$. The analytic continuation of some basis of solutions at $z_0$ around a closed path avoiding singular points is a new basis for the solution space at $z_0$ (see \cite{In}, page 357). This leads to a representation
$$
\rho\colon\pi_1(\PP^1\setminus{S},z_0)\rightarrow\textnormal{Gl}(n,\CC)
$$
(where $S$ is the set of singular points), called the \emph{monodromy representation}. Under this map, the image of a closed path with base point $z_0$ is defined to be the change of basis matrix corresponding to the analytic continuation of the given basis of the solution space. Strictly speaking, the monodromy representation depends on a choice of basis of solutions at the base point and so we bear in mind that the representation is only defined up to conjugacy.
The image of the monodromy representation is called the \emph{monodromy group}.

In the case of a differential system, since any solution of $\omega^\prime=A\omega$ about some point $z_0$ is a vector function, we may arrange a basis of solutions as the columns of a matrix
$$
\Omega := (\omega_1 |\, \omega_2 |\, \ldots |\, \omega_n )
$$
where $\det\Omega\neq{0}$ and $\omega_i^\prime=A\omega_i$ for $i = 1,\ldots,n$. Such a matrix is known as a fundamental matrix for the system and satisfies
$$
\Omega^\prime = A\Omega.
$$
If $M\in\textnormal{Gl}(n,\CC)$ is any constant matrix, then $\Omega.M$ is also a fundamental matrix for our system since $(\Omega.M)^\prime = \Omega^\prime.M = A(\Omega.M)$.

Analytic continuation of any solution around a closed loop $\gamma$ gives another solution and, as before, we obtain a monodromy representation
$$
\rho\colon\pi_1(\PP^1\setminus{S},z_0)\rightarrow\textnormal{Gl}(n,\CC)
$$
sending the fundamental matrix $\Omega$ to the fundamental matrix $\Omega.\rho(\gamma)$. This representation is defined up to conjugacy.
\begin{remark}
If the system
\begin{equation}
\omega^\prime=A\omega\label{sysA}
\end{equation}
has a fundamental matrix $\Omega$, then for any matrix $M\in\textnormal{Gl}(n,\CC(z))$, the system
\begin{equation}
\omega^\prime=B\omega\label{sysB}
\end{equation}
where
$$
B = M.A.M^{-1}+M^\prime.M^{-1}
$$
has a fundamental matrix $M.\Omega$. In particular, if $M\in\textnormal{Gl}(n,\CC)$, then $M^\prime=0$ and $B=M.A.M^{-1}$.
\end{remark}

We shall find it more convenient to deal with Fuchsian differential systems than with Fuchsian differential equations and so we introduce the following theorem.
\begin{theorem}[\cite{Bo}]\label{equivSys}
For any Fuchsian differential equation, there is some Fuchsian differential system with the same singular points and the same monodromy representation.
\end{theorem}
\begin{proof}
This is proved at length in \cite{Bo} pp.89-94 for Fuchsian differential equations of arbitrary order. The converse of this statement does not hold.

We shall find it useful to take the general proof of \cite{Bo} and create a specific algorithm for converting a Fuchsian equation into a Fuchsian system with the same monodromy representation in the special case of \emph{second order} equations.

Let
\begin{equation}\label{secondOrderEquation}
a_2(z)\frac{d^2f}{dz^2}+a_1(z)\frac{df}{dz}+a_0(z)f=0
\end{equation}
be a second order Fuchsian differential equation. If necessary, apply a M\"obius transformation to the parameter $z$ to assume that the equation is nonsingular at $\infty$ and has a regular singular point at 0. Also, by dividing out any repeated roots if necessary, we assume
$$
a_2(z) = \prod_{i-1}^m(z-\alpha_i).
$$
where $\{\alpha_1,\ldots,\alpha_m\}$ are the distinct singular points of the equation (and we no longer assume that $a_1$ and $a_0$ are polynomials).

The proof of theorem~\ref{equivSys} in \cite{Bo} states that there exists a Laurent polynomial
$$
Q = q_{m-2}z^{2-m}+\ldots+q_1z^{-1}+q_0\ \in\CC[z^{-1}]
$$
such that the Fuchsian equation (\ref{secondOrderEquation}) has the same singular points and monodromy representation as the differential system
$$
\frac{d}{dz}\omega = A\omega
$$
where
$$
\omega = \begin{pmatrix}f\\Qf+a_2z^{2-m}\frac{df}{dz}\end{pmatrix}
$$
and $A$ has the form
$$
A = M.B.M^{-1} + \frac{d}{dz}M.M^{-1}
$$
with
$$
B = \begin{pmatrix}
0&\frac{1}{a_2}\\
-a_0&\frac{-a_1+a_2^\prime}{a_2}
\end{pmatrix}
$$
and
$$
M = \begin{pmatrix}
1&0\\
Q&z^{2-m}
\end{pmatrix}.
$$
The only unknowns are the coefficients $q_{m-2},\ldots,q_1,q_0$. The value of $q_0$ is arbitrary as the choice will not affect the monodromy representation of the resulting system. We shall typically take $q_0 = 0$.

The remaining coefficients $q_{m-2},\ldots,q_1$ are determined by the condition that $A$ is Fuchsian at 0 so that $zA$ is nonsingular at 0. In practice, this comes down to solving some simultaneous polynomial equations in these unknowns.
\end{proof}

\begin{definition}[Hypergeometric Differential Equation]\label{hypergeometric}
The Hypergeometric differential equation $_2F_1(\alpha,\beta,\gamma)$ is the second order Fuchsian ODE
\begin{equation}\label{hypergeom}
z(1-z)\frac{d^2f}{dz^2}+(\gamma-(\alpha+\beta+1)z)\frac{df}{dz}-\alpha\beta{f} = 0.
\end{equation}
From \cite{Bo}, the system
$$
\frac{d}{dz}\omega = A \omega
$$
with
\begin{equation}\label{hypergeomSys}
A = \frac{1}{z}\begin{pmatrix}0&0\\-\alpha\beta&-\gamma\end{pmatrix}+\frac{1}{z-1}\begin{pmatrix}0&1\\0&\gamma-\alpha-\beta\end{pmatrix}
\end{equation}
has the same monodromy representation.
\end{definition}
A few examples of these hypergeometric differential equations occur later as the symmetric square--root of some Picard--Fuchs differential equations related to families of K3 surfaces. As will the following.
\begin{definition}[Lam\'e Differential Equation]\label{lame}
The differential equation
$$
L_n(p,B):= p(z)\frac{d^2f}{dz^2} + \frac{1}{2}p^\prime(z)\frac{df}{dz} - (n(n+1)z+B)f = 0
$$ with
$$
p(z) = 4(z-\alpha_1)(z-\alpha_2)(z-\alpha_3)
$$
and
$$
\alpha_1+\alpha_2+\alpha_3=0
$$
is known as Lam\'e's differential equation. It is an example of a Fuchsian equation with 4 regular singular points (it has a singularity at infinity).

The Fuchsian system from \cite{Lame}:
$$\omega^\prime = \sum_{i=1}^3\frac{R_i}{z-\alpha_i}\omega
$$
with residue matrices
$$
R_1=\begin{pmatrix}0&1\\0&\frac{1}{2}\end{pmatrix}\ \ R_2=\begin{pmatrix}0&0\\l_1&-\frac{1}{2}\end{pmatrix}\ \ R_3=\begin{pmatrix}0&0\\l_2&-\frac{1}{2}\end{pmatrix}
$$
where
$$
l_1 = \frac{\alpha_2 n(n+1)+B}{4(\alpha_2-\alpha_3)}\ \ \ \textnormal{and}\ \ \ l_1+l_2 = \frac{n(n+1)}{4}
$$
has the same monodromy representation as the Lam\'e differential equation.
\end{definition}
\begin{definition}\label{elementary}
Following \cite{In}, we say that a regular singular point of a second order differential equation is \emph{elementary} if the eigenvalues of the corresponding residue matrix differ by $\frac{1}{2}$.

The Lam\'e differential equation has three elementary singular points and an arbitrary regular singular point at infinity. The class of \emph{generalised Lam\'e differential equations} consists of those Fuchsian equations with $p\geq 4$ regular singular points, $p-1$ of which are elementary.
\end{definition}

\begin{theorem}[\cite{CL}, pages 109, 117 and 119]\label{local}
A differential system
$$
\omega^\prime=\left(\frac{R_i}{z-\alpha_i}+B(z)\right)\omega
$$
with $R_i$ constant and $B(z)$ analytic in some disc $\Delta$ centred at $\alpha_i$ has a fundamental matrix of solutions, $\Omega$, of the form
$$
\Omega = S(z)(z-\alpha_i)^P
$$
valid in the punctured disc $\Delta\setminus{\{\alpha_i\}}$. Here, $S(z)$ is a single--valued matrix function and P is a constant matrix. Furthermore, if no two \emph{distinct} eigenvalues of the residue matrix, $R_i$, differ by integers, then we may take $P=R_i$, or indeed we may take $P$ to be the canonical form of $R_i$. The local monodromy transformation about $\alpha_i$ is given by
$$
M_{\alpha_i} = \exp(2\pi i P).
$$
\end{theorem}
\begin{proof}
This well-known result is used in Morrison's paper, \cite{Mo}, and is patched together from chapter 4 of \cite{CL}, in particular theorems 1.1, 3.1 and 4.1. As is pointed out in \cite{Mo}, the case where two distinct eigenvalues of $R_i$ do differ by an integer is dealt with in \cite{CL} page 120.
\end{proof}
This result is usually used to calculate the local monodromy of a differential \emph{equation} by writing it in the form
\begin{equation}\label{transformed}
\Phi^n(f) + \sum_{i=1}^{n}a_{n-i}\Phi^{n-i}(f) = 0
\end{equation}
with $\Phi = (z-\alpha_i)\frac{d}{dz}$. This equation has the same \emph{local} monodromy as the system
$$
\Phi\Omega = A\Omega
$$
where 
$$
\Omega = 
\begin{pmatrix}
f\\\Phi(f)\\\vdots\\\Phi^{n-1}(f)
\end{pmatrix}
$$
and
$$
A = \begin{pmatrix}
0&1&\phantom{0}&\ldots&0\\
\phantom{0}&0&1&\ldots&0\\
 & &\ddots&\ddots&\vdots\\
-a_0&-a_1&-a_2&\ldots&-a_{n-1}
\end{pmatrix}.
$$
Because $A$ is guaranteed to be holomorphic at $\alpha_i$, on dividing by $(z-\alpha_i)$, this system is in the form required by theorem \ref{local}.

Since this method requires the repeated transformation of the equation into the form (\ref{transformed}) for each singular point, we find it to be inconvenient. We shall prefer to use theorem \ref{equivSys} once to find a differential system with the same monodromy representation. The residue matrices are then available without further calculation and the local monodromies are easily determined by the following corollary to theorem \ref{local}.
\begin{corollary}\label{disteigen}
Let
$$
\omega^\prime=\left(\sum_{i=1}^k\frac{R_i}{z-\alpha_i}\right)\omega
$$
be a Fuchsian system. If no two distinct eigenvalues of the residue matrix $R_i$ differ by an integer, then the local monodromy about $\alpha_i$, with respect to some unspecified basis, is
$$
\exp(2 \pi \mathbbm{i} P)
$$
where $P$ is the canonical form of residue matrix $R_i$.
\end{corollary}

Whenever we are only interested in the solutions of some differential equation up to multiplication by scalars, we will correspondingly only be interested in the \emph{projective} monodromy group. That is, the monodromy group modulo scalar matrices. Notice that about an elementary singular point, as defined in definition \ref{elementary}, the projective local monodromy transformation has order 2. Hence, the projective monodromy group of a generalised Lam\'e differential equation may be generated by elements of order 2.

\section{Picard--Fuchs Equations for K3 Surfaces}
\subsection{Periods}

We present a method for computing the Picard-Fuchs differential equation satisfied by the periods of a one--parameter family of K3 surfaces. This calculation is carried out in a number of examples.

\begin{definition}
Let $X_{\lambda_0}$ be a smooth K3 surface. A \emph{period} of $X_{\lambda_0}$ is a complex number of the form $\int_{\gamma}\omega$ where $\omega\in H^{2,0}(X_{\lambda_0})$ is a non--zero holomorphic differential 2--form and $\gamma\in H_2(X_{\lambda_0},\ZZ)$ is some cycle.
\end{definition}

If $X\rightarrow\Delta$ is a family of K3 surfaces defined over a disc $\Delta\subset\CC$, then choosing a basis $\left<\gamma_0(\lambda),\ldots,\gamma_{21}(\lambda)\right> = H_2(X_{\lambda},\ZZ)$ that varies smoothly with $\lambda$ allows us to consider the so called period point
$$
\left(\int_{\gamma_0(\lambda)}\omega,\ldots,\int_{\gamma_{21}(\lambda)}\omega\right)\in\PP^{21}.
$$
This point lies in \emph{projective} space since, on any K3 surface, $\omega$ is only uniquely defined up to a scalar multiple. Suppose our family $X\rightarrow\Delta$ is a family of K3 surfaces invariant under the action of some symplectic group $G$. Then the sublattice
$$
M_G := S_{X,G}\oplus [\kappa]
$$
is primitively embedded in the Picard lattice $S_X$ where $[\kappa]$ is a very ample divisor class.

If $\rho = \textnormal{rank}(M_G)$, then the lattice $H^2(X_\lambda,\ZZ)/M_G$ has rank $22-\rho$. Since for any algebraic class $\gamma$
$$
\int_X \omega\wedge\gamma = 0,
$$
we may define the period point
$$
\left( \int_X \omega\wedge\gamma_1,\ldots,\int_X\omega\wedge\gamma_{22-\rho} \right) \in \PP^{21-\rho}
$$
where $\{\overline{\gamma}_1,\ldots,\overline{\gamma}_{22-\rho}\}$ is a basis for $H^2(X,\ZZ)/M_G$. It is well--known that the periods are the solutions of a linear ordinary differential equation of degree $22-\rho$, called the Picard--Fuchs equation (see for example \cite{VY}).

\subsection{Determining the Picard--Fuchs Equation}\label{det}

In \cite{Mo}, the Picard--Fuchs equations for some one parameter families of Calabi--Yau threefolds are calculated. We follow this method, which uses Griffith's description of the primitive cohomology of a hypersurface (see \cite{Griff}), to find the differential equations for quasismooth families of K3 hypersurfaces in weighted projective space.

Following \cite{Mo}, let
$$
\Omega:=\sum_{j=0}^3(-1)^j\,k_jx_j\,dx_0\wedge\ldots\wedge\widehat{dx_j}\wedge\ldots\wedge{dx_3},
$$
be a differential form on the weighted projective space $\PP(k_0,k_1,k_2,k_3)$. Then all rational differential 3--forms on $\PP(k_0,k_1,k_2,k_3)$ are expressible as $A\Omega / B$ where $A$ and $B$ are weighted homogeneous polynomials with $\deg A\Omega = \deg B$.

For any hypersurface
$$
X:(Q=0)\subset\PP(k_0,k_1,k_2,k_3),
$$
the primitive cohomology in $H^2(X,\ZZ)$ is described by residues of the differential forms $P\Omega / Q^l$ on $\PP(k_0,k_1,k_2,k_3)\setminus\{Q=0\}$. The residue $\text{Res}_X(P\Omega / Q^l)\in H^2(X,\ZZ)$ is the differential form satisfying
$$
\int_\gamma\text{Res}_X\left(\frac{P\Omega}{Q^l}\right)=\frac{1}{2\pi i}\int_\Gamma \frac{P\Omega}{Q^l}
$$
for any $\gamma\in H_2(X,\ZZ)$, where $\Gamma$ is a tubular neighbourhood of $\gamma$. This can be thought of as a generalisation of the residue theorem on $\CC$.

Since integration over $\Gamma$ commutes with differentiation by $\lambda$, we see that if the periods satisfy the differential equation
$$
D\omega:=\left(c_k\frac{d^k}{d\lambda^k}+c_{k-1}\frac{d^{k-1}}{d\lambda^{k-1}}+\ldots+c_0\right)\omega=0
$$
then
$$
\int_\Gamma D \frac{\Omega}{Q}\ =\ 0
$$
so that $D\frac{\Omega}{Q}$ is an exact differential form. Therefore, finding the Picard--Fuchs equation is equivalent to finding a linear dependence relation between
$$
\frac{\Omega}{Q},  \frac{d}{d\lambda}\frac{\Omega}{Q},\ldots,\frac{d^k}{d\lambda^k}\frac{\Omega}{Q}
$$
modulo exact forms.

If we let
$$
\varphi:=\frac{1}{Q^l}\sum_{i<j}(k_ix_iA_j-k_jx_jA_i)dx_0\wedge\stackrel{i,j}{\widehat{\ldots}}\wedge dx_3
$$
where the $A_i$ are any homogeneous polynomials of weighted degree $$\deg{A_i}=\deg{Q^l}+k_i-k_0-k_1-k_2-k_3,$$
then an easy calculation shows
$$
d\varphi=\frac{l\;\sum A_j \frac{\partial{Q}}{\partial{x_j}}}{Q^{l+1}}\;\Omega\hspace{7pt}-\hspace{5pt}\frac{\sum\frac{\partial{A_j}}{\partial{x_j}}}{Q^l}\;\Omega
$$
which we shall write as
$$
\frac{l\;\sum A_j \frac{\partial{Q}}{\partial{x_j}}}{Q^{l+1}}\;\Omega\hspace{10pt}\equiv_{d\varphi}\hspace{7pt}\frac{\sum\frac{\partial{A_j}}{\partial{x_j}}}{Q^l}\;\Omega.
$$
This gives a practical way to reduce the pole modulo exact forms in order to find our dependence relation.

The full algorithm is best explained with the following example. Consider the following family of K3 surfaces:
$$
X_\lambda\ :\ \left( \sum_{i=0}^3 x_i^4 +\lambda\left(\sum_{i=0}^3 x_i^2\right)^2=0 \right) \subset\PP^3.
$$
This is example XI of section \ref{p3} and has Picard number 19 so that the Picard--Fuchs equation will have order 3.

Write $s_4:=\sum x_i^4$, $q:=(\sum x_i^2)^2$, and $Q:=s_4+\lambda q$. Then, since $\Omega$ does not depend on $\lambda$, we have
$$
\frac{d^r}{d\lambda^r}\frac{\Omega}{Q} = \frac{(-1)^r\,r!\,q^{r}\,\Omega}{Q^{r+1}}.
$$
If the differential equation to be determined is of the form
$$
c_3\frac{d^3}{d\lambda^3}\frac{\Omega}{Q}\ +\ c_{2}\frac{d^{2}}{d\lambda^{2}}\frac{\Omega}{Q}\ +\ c_{2}\frac{d}{d\lambda}\frac{\Omega}{Q}\ +\ c_0\frac{\Omega}{Q}\ \equiv_{d\varphi}\ 0,
$$
then we must have
$$
c_3\frac{d^3}{d\lambda^3}\frac{\Omega}{Q}\ =\ -6\,c_3\,q^3 \frac{\Omega}{Q^4}\ =\ 3\sum 2A_j\frac{\partial Q}{\partial x_j}\frac{\Omega}{Q^4}
$$
for some polynomials $A_j\in\CC[x_0,x_1,x_2,x_3,\lambda]$ so that $c_3\frac{d^3}{d\lambda^3}\frac{\Omega}{Q}$ can be reduced modulo exact forms and expressed as a linear combination of lower order derivatives. In other words, we must have
$$
-c_3\,q^3=\sum A_j\frac{\partial Q}{\partial x_j}\in\left\langle\frac{\partial Q}{\partial x_0},\frac{\partial Q}{\partial x_1},\frac{\partial Q}{\partial x_2},\frac{\partial Q}{\partial x_3}\right\rangle.
$$
In particular, if the family degenerates at $\lambda_0\in\CC$, then $\frac{\partial Q}{\partial x_i}=0$ for each $i$ at some point $(x_0,\ldots,x_3)$, and so we must have $c_3({\lambda_0})=0$. Since our particular example degenerates at $\lambda=-1$, $-\frac{1}{2}$, $-\frac{1}{3}$, $-\frac{1}{4}$ (and at $\infty$), we take $c_3=(1+\lambda)(1+2\lambda)(1+3\lambda)(1+4\lambda)$ (or, if this process were to fail, some power of this product).

To express $-c_3 q^3$ as an element of the ideal $J:=\left\langle\frac{\partial Q}{\partial x_0},\frac{\partial Q}{\partial x_1},\frac{\partial Q}{\partial x_2},\frac{\partial Q}{\partial x_3}\right\rangle$, a Gr\"obner basis for $J$ must be calculated. The normal form for $-c_3 q^3$ with respect to this Gr\"obner basis is obtained and, having kept track of the change of basis from $\left\{\frac{\partial Q}{\partial x_i}\right\}$ to the Gr\"obner basis, the expression $-c_3 q^3=\sum A_j \frac{\partial Q}{\partial x_j}$ is found.

The intermediate polynomials $A_i$ are differentiated to obtain
$$
c_3\frac{d^3}{d\lambda^3}\frac{\Omega}{Q}\ \equiv_{d\varphi} 2 \sum \frac{\partial A_i}{\partial x_i}\frac{\Omega}{Q^3}.
$$
Although the right hand side is not reducible further, it must be true that $\sum \frac{\partial A_i}{\partial x_i}\in\left\langle\frac{\partial Q}{\partial x_0},\frac{\partial Q}{\partial x_1},\frac{\partial Q}{\partial x_2},\frac{\partial Q}{\partial x_3},-q^2\right\rangle$ so that for some $c_2\in\CC[\lambda]$, we may reduce $c_3\frac{d^3}{d\lambda^3}\frac{\Omega}{Q}+c_2\frac{d^2}{d\lambda^2}\frac{\Omega}{Q}$ modulo an exact form.

Again, to find the expression $\sum \frac{\partial A_i}{\partial x_i} = \sum B_i\frac{\partial Q}{\partial x_i}\,-\, c_2 q^2$, we calculate a Gr\"obner basis for $\langle J, -q^2\rangle$ and find the normal form for $\sum \frac{\partial A_i}{\partial x_i}$ as before.

In our example, this yields
$$
c_2 = 144 \lambda^3+225 \lambda^2+ 105 \lambda+15
$$
and
$$
c_3\frac{d^3}{d\lambda^3}\frac{\Omega}{Q}+c_2\frac{d^2}{d\lambda^2}\frac{\Omega}{Q}\ \equiv_{d\varphi} \sum \frac{\partial B_i}{\partial x_i}\frac{\Omega}{Q^2}.
$$
We repeat the last step by finding an expression for the right hand side in terms of $\left\{\frac{\partial Q}{\partial x_i}\right\}$ and $q$ to obtain
$$
c_3\frac{d^3}{d\lambda^3}\frac{\Omega}{Q}+c_2\frac{d^2}{d\lambda^2}\frac{\Omega}{Q}+c_1\frac{d}{d\lambda}\frac{\Omega}{Q}\ \equiv_{d\varphi} \sum \frac{\partial C_i}{\partial x_i}\frac{\Omega}{Q}.
$$
It is necessarily true that $c_0:=-\sum \frac{\partial C_i}{\partial x_i}$ is a polynomial only in $\lambda$ and does not depend on the $x_i$. We have now determined the Picard--Fuchs differential equation for our family.

Appendix \ref{macaulay} contains the full algorithm to compute the Picard--Fuchs differential equation. The algorithm is implemented using the algebra system Macaulay2 since it handles Gr\"obner bases efficiently and keeps track of the change of basis matrix.

\subsection{The Symmetric Square Root}

In section \ref{moduli}, we looked at a construction of the coarse moduli space of lattice polarised K3 surfaces. This construction is straightforward once the orthogonal complement of the polarising lattice is determined (or stated as an assumption), although it can be difficult to determine this transcendental lattice in an explicit example. We shall construct this moduli space without explicitly determining the generic transcendental lattice by instead using the Picard--Fuchs differential equation. 

Recall that the coarse moduli space of K3 surfaces polarised by a lattice $S$ is given by
$$
D_S/\Gamma_S
$$
where
$$
D_S = \{\ \omega\in\PP(T\otimes\CC)\ |\ \left<\omega,\omega\right>=0,\ \left<\omega,\bar\omega\right>>0\ \}
$$
is the period space, $T = S^\perp$, and $\Gamma_S$ is a group of automorphisms if $D_S$. The first condition $\left<\omega,\omega\right> = 0$ ensures that the periods lie on the non--degenerate quadric defined by the quadratic form on $T$. Hence the solutions of the Picard--Fuchs differential equation also satisfy this non--degenerate quadratic relationship. This is noted in \cite{pet} and \cite{Dor}.

\begin{definition}
The \emph{symmetric square} of a second order ODE
\begin{equation}\label{one}
    a_2\frac{d^2\omega}{d\lambda^2}+a_1\frac{d\omega}{d\lambda}+a_0\omega=0
\end{equation}
is the third order ODE
\begin{equation}\label{two}
a_2^2\frac{d^3\omega}{d\lambda^3} + 3a_1a_2\frac{d^2\omega}{d\lambda^2} + (a_2(a_0+a_1^\prime)+a_1(2a_1-a_2^\prime))\frac{d\omega}{d\lambda} + 2(a_2a_0^\prime+a_0(2a_1-a_2^\prime))\omega=0.
\end{equation}
\end{definition}

If $\{\alpha,\beta\}$ are independent solutions of (\ref{one}), then it can be shown that $\{\alpha^2, \alpha\beta,\beta^2\}$ are linearly independent solutions of the symmetric square (\ref{two}). Conversely, if the solutions of a third order ODE satisfy a non--degenerate quadratic relationship, then it is the symmetric square of a second order ODE. The Picard--Fuchs equation of a family of rank 19 lattice polarised K3 surfaces is thus the symmetric square of some second order ODE (see \cite{Dor}).

We shall call this second order ODE the \emph{symmetric square root} of the Picard--Fuchs differential equation.

Typically, the leading coefficient of the Picard--Fuchs equation is not a square as is suggested by the form of (\ref{two}). The coefficients of (\ref{two}) must share a common factor for each square--free factor of the top coefficient. This imposes some further relationships on the coefficients of the Picard--Fuchs differential equation. For example, if the top coefficient is entirely square--free, then $a_2$ must divide each coefficient of (\ref{two}). Because $a_2,\,a_1,$ and $a_0$ share no common factor, we are forced to conclude that
$$
a_2\;|\;(2a_1-a_2^\prime)\hspace{20pt}\text{or}\hspace{20pt}2a_1-a_2^\prime=0.
$$
However, $\deg(a_2)>\deg(2a_1-a_2^\prime)$ since (\ref{one}) is Fuchsian whenever (\ref{two}) is Fuchsian so the first possibility cannot occur. Therefore we conclude that $a_1=\frac{1}{2}a_2^\prime$ whenever the top coefficient is square free. Substituting this back into (\ref{two}), we see that the \emph{typical} form for our Picard--Fuchs equation is
$$
a_2\frac{d^3\omega}{d\lambda^3}+\frac{3}{2}a_2^\prime\;\frac{d^2\omega}{d\lambda^2}+(a_0+\frac{1}{2}a_2^{\prime\prime})\frac{d\omega}{d\lambda}+2 a_0^\prime\omega=0.
$$

and the symmetric square root will \emph{typically} be of the form
\begin{equation}\label{typicalform}
a_2\frac{d^2\omega}{d\lambda^2}+\frac{1}{2}a_2^\prime\frac{d\omega}{d\lambda}+a_0\omega=0
\end{equation}

\begin{remark}
From the lattice polarised description of the period space, we know that there is some basis of solutions to the Picard--Fuchs equation coming from integral cohomology classes. With respect such a basis, the monodromy representation is integral as it is a finite index subgroup of $\textnormal{SO}(T)$. Thanks to theorem \ref{squares}, the monodromy group of the symmetric square root is a discrete subgroup of $\textnormal{PSl}(2,\RR)$ with traces in some real quadratic number field. More on this in chapter \ref{mono}.
\end{remark}

\begin{remark}
The existence of the symmetric square root of our Picard--Fuchs differential equation is also convenient in some unexpected ways. For example, in order to calculate the local monodromies about a degenerate point of a differential system, we may use corollary \ref{disteigen} so long as the distinct eigenvalues of the residue matrix do not differ by an integer. If two eigenvalues \emph{do} differ by an integer, we are forced to make a series of non-linear changes to the system to reduce it to a system with good eigenvalues. This procedure is described in \cite{In}, but is awkward and best avoided if possible. Although this situation does occur for some of the examples to appear later, it can be avoided by virtue of the fact that the eigenvalues of the symmetric square root equation do not differ by an integer in all of our examples. Hence, in practice, we immediately pass to the symmetric square root of the equation to make our life easier.
\end{remark}

\section{Picard--Fuchs Equations for Families in $\PP(1,1,1,3)$}

We calculate the Picard--Fuchs differential equation for the families of K3 surfaces found in section \ref{p1113}. Since these families are invariant under groups with Mukai number 6, they have generic Picard number 19 and so their Picard--Fuchs equations are third order. We also give the second order symmetric square roots of these ODEs and, for families I, II and III, we change the parameter for the family from $\lambda$ to $\mu=\lambda^3$ to take into account the isomorphism $X_\lambda\cong X_{\omega\lambda}$ where $\omega^3=1$.

\begin{Family}{I}{$SD_{16}$}
The family:
$$X_\lambda : t^2 = x_0^6 + x_1^5x_2 + x_1x_2^5 - 3 \lambda x_0^2x_1^2x_2^2.$$
degenerates at: $\lambda = \infty$ and $\lambda^3 = 1$

The Picard--Fuchs operator for $X_\lambda$ is
$$
(\lambda^3-1)\frac{d^3}{d\lambda^3} + \frac{9}{2}\lambda^2\frac{d^2}{d\lambda^2} + \frac{13}{4}\lambda\frac{d}{d\lambda} + \frac{1}{8}
$$
which is the symmetric square of
$$
(\lambda^3-1)\frac{d^2}{d\lambda^2} + \frac{3}{2}\lambda^2\frac{d}{d\lambda} + \frac{\lambda}{16}.
$$
After substituting $\mu = \lambda^3$, we get
$$
\mu(1-\mu)\frac{d^2}{d\mu^2} + \left(\frac{2}{3}-\frac{7}{6}\mu\right)\frac{d}{d\mu} - \frac{1}{144}
$$
which is the hypergeometric differential equation $_2F_1(\frac{1}{12},\frac{1}{12}, \frac{2}{3})$. 
\end{Family}
 
\begin{Family}{II}{$\mathfrak{A}_{4,3}$}
The family
$$t^2 = x_0^6+x_1^6+x_2^6-3\lambda x_0^2x_1^2x_2^2$$
degenerates at $\lambda=\infty$ and where $(\lambda^3-1)$

The Picard--Fuchs operator for $X_\lambda$ is
$$
(\lambda^3-1)\frac{d^3}{d\lambda^3} + \frac{9}{2}\lambda^2\frac{d^2}{d\lambda^2} + \frac{13}{4}\lambda\frac{d}{d\lambda} + \frac{1}{8}
$$
This is the same as for family I and so it follow that this is the symmetric square of
$$
(\lambda^3-1)\frac{d^2}{d\lambda^2} + \frac{3}{2}\lambda^2\frac{d}{d\lambda} + \frac{\lambda}{16}.
$$
After substituting $\mu = \lambda^3$, we again get
$$
\mu(1-\mu)\frac{d^2}{d\mu^2} + \left(\frac{2}{3}-\frac{21}{18}\mu\right)\frac{d}{d\mu} - \frac{1}{144}
$$
which is the hypergeometric differential equation $_2F_1(\frac{1}{12},\frac{1}{12}, \frac{2}{3})$.
\end{Family}

\begin{Family}{III}{$C_7\rtimes{C_3}$}
The family
$$X_\lambda : t^2 = x_0^5x_1  + x_1^5x_2 + x_2^5x_0 - 3\lambda x_0^2x_1^2x_2^2$$
degenerates at $\lambda = \infty$ and $\lambda^3 = 1$

The Picard--Fuchs operator for $X_\lambda$ is
$$
(\lambda^3-1)\frac{d^3}{d\lambda^3} + \frac{9}{2}\lambda^2\frac{d^2}{d\lambda^2} + \frac{13}{4}\lambda\frac{d}{d\lambda} + \frac{1}{8}
$$
This is the same as for families I and II and again leads to the hypergeometric differential equation $_2F_1(\frac{1}{12},\frac{1}{12}, \frac{2}{3})$.
\end{Family}

\begin{Family}{IV}{$3^2C_4$}
The family
$$t^2 = p + 3 \lambda q$$
where
\begin{eqnarray*}
p & := & x_0^6 + x_1^6 + x_2^6 - 10(x_0^3x_1^3 + x_0^3x_2^3 + x_1^3x_2^3),\\ 
q & := & x_0x_1x_2(x_0^3 + x_1^3 + x_2^3) - 2(x_0^3x_1^3 + x_0^3x_2^3 + x_1^3x_2^3) + 3x_0^2x_1^2x_2^2
\end{eqnarray*}
degenerates at $\lambda = \infty$ and where $(\lambda+1)(\lambda+2)(\lambda^2 + 16\lambda + 16) = 0$.

The Picard--Fuchs operator for $X_\lambda$ is
\begin{eqnarray*}
(\lambda+1)(\lambda+2)(\lambda^2+16\lambda+16)\frac{d^3}{d\lambda^3}&\\
+ (6\lambda^3+\frac{171}{2}\lambda^2+198\lambda+120)\frac{d^2}{d\lambda^2}&\\
+ (\frac{27}{4}\lambda^2+63\lambda+72)\frac{d}{d\lambda}&\\
+ \frac{3}{4}\lambda+3&
\end{eqnarray*}
which is the symmetric square of
$$
p(\lambda)\frac{d^2}{d\lambda^2} + \frac{1}{2}p^\prime(\lambda)\frac{d}{d\lambda} + \frac{3}{16}\lambda^2+\frac{3}{2}\lambda+\frac{3}{2}.
$$

with $p(\lambda) = (\lambda+1)(\lambda+2)(\lambda^2+16\lambda+16)$.
\end{Family}

\begin{Family}{V}{$\mathfrak{A}_5$}
The family
$$t^2 = p+\lambda q^3$$
where 
$$ p := 8x_0^6 + 30x_0^2x_1^2x_2^2 + 3x_0(x_1^5 + x_2^5) + 5x_1^3x_2^3$$
and 
$$q := x_0^2 + x_1x_2$$
degenerates at $\lambda = \infty$ and where $(\lambda+5)(\lambda+8)(9\lambda+40)=0$.

The Picard--Fuchs operator is
$$
(\lambda+5)(\lambda+8)(\lambda+\frac{40}{9})\frac{d^3}{d\lambda^3}
+
(\frac{9}{2}\lambda^2+\frac{157}{3}\lambda+\frac{440}{3})\frac{d^2}{d\lambda}
+
(\frac{29}{9}\lambda+\frac{671}{36})\frac{d}{d\lambda}
+
\frac{1}{9}.
$$

This is the symmetric square of
$$
p(\lambda)\frac{d^2}{d\lambda^2}
+
\frac{1}{2}p^\prime(\lambda)\frac{d}{d\lambda}
+
\left(\frac{1}{18}\lambda+\frac{43}{144}\right)
$$
where $p(\lambda) = (\lambda+5)(\lambda+8)(\lambda+\frac{40}{9})$

After the substitution $\lambda = \frac{\mu-157}{27}$, this differential equation is of Lam\'e type (see definition \ref{lame}) with $p(\mu)=4(\mu+59)(\mu-22)(\mu-37)$, $n=-\frac{1}{3}$, $B = \frac{95}{36}$.
\end{Family}

\section{Picard--Fuchs Equations for Families in $\PP^3$}
We continue by giving the Picard--Fuchs differential equations for the families of K3 surfaces in $\PP^3$.

\begin{Family}{VI}{$4^2\mathfrak{A}_4$}
The family
$$X_\lambda :  x_0^4+x_1^4+x_2^4+x_3^4 + 4\lambda x_0x_1x_2x_3=0$$
degenerates at $\lambda = \infty$ and where $\lambda^4 = 1$.

The Picard--Fuchs operator for $X_\lambda$ is
$$
(\lambda^4-1)\frac{d^3}{d\lambda^3}+6\lambda^3\frac{d^2}{d\lambda^2} + 7\lambda^2\frac{d}{d\lambda} + \lambda
$$
which is the symmetric square of
$$
(\lambda^4-1)\frac{d^2}{d\lambda^2} + 2\lambda^3\frac{d}{d\lambda} + \frac{\lambda^2}{4}.
$$
After substituting $\mu = \lambda^4$, we get
$$
\mu(1-\mu)\frac{d^2}{d\mu^2} + \left(\frac{3}{4}-\frac{5}{4}\mu\right)\frac{d}{d\mu} - \frac{1}{64}
$$
which is the hypergeometric differential equation $_2F_1(\frac{1}{8},\frac{1}{8}, \frac{3}{4})$.
\end{Family}

\begin{Family}{VII}{$C_7\rtimes{C_3}$}
The family
$$
X_\lambda : x_0x_1^3 + x_1x_2^3 + x_2x_0^3 + x_3^4 + 4 \lambda x_0x_1x_2x_3 = 0
$$
degenerates at $\lambda = \infty$ and where $\lambda^4 = 1$.

The Picard--Fuchs operator for $X_\lambda$ is
$$
(\lambda^4-1)\frac{d^3}{d\lambda^3} + 6\lambda^3\frac{d^2}{d\lambda^2} + 7\lambda^2\frac{d}{d\lambda} + \lambda
$$
This is the same as for family VI and so it follows that this is the symmetric square of
$$
(\lambda^4-1)\frac{d^2}{d\lambda^2} + 2\lambda^3\frac{d}{d\lambda} + \frac{\lambda^2}{4}.
$$
After substituting $\mu = \lambda^4$, we again get
$$
\mu(1-\mu)\frac{d^2}{d\mu^2} + \left(\frac{3}{4}-\frac{5}{4}\mu\right)\frac{d}{d\mu} - \frac{1}{64}
$$
which is the hypergeometric differential equation $_2F_1(\frac{1}{8},\frac{1}{8}, \frac{3}{4})$.
\end{Family}

\begin{Family}{VIII}{$SD_{16}$}
The family
$$X_\lambda :  x_0^4+x_1^4+x_2^4+x_3^4 + 2\sqrt{2}\lambda x_0x_1(x_2^2+ix_3^2)=0$$
degenerates at $\lambda = \infty$ and where $\lambda^4 = 1$.

The Picard--Fuchs operator for $X_\lambda$ is
$$
(\lambda^4-1)\frac{d^3}{d\lambda^3}+6\lambda^3\frac{d^2}{d\lambda^2} + 7\lambda^2\frac{d}{d\lambda} + \lambda
$$
This is the same as for families VI and VII. Again, this leads to the hypergeometric differential equation $_2F_1(\frac{1}{8},\frac{1}{8}, \frac{3}{4})$.
\end{Family}

\begin{Family}{IX}{$\textnormal{Hol}(C_5)$}
The family
$$X_\lambda :  x_0^3x_1 + x_1^3x_3 + x_3^3x_2 - x_2^3x_0 + \left(\frac{1+i}{2}\right)\lambda(x_0^2x_3^2 + x_1^2x_2^2)=0$$
degenerates at $\lambda = \infty$ and where $\lambda^4 = 1$.

The Picard--Fuchs operator for $X_\lambda$ is
$$
(\lambda^4-1)\frac{d^3}{d\lambda^3}+6\lambda^3\frac{d^2}{d\lambda^2} + 7\lambda^2\frac{d}{d\lambda} + \lambda
$$
This is the same as for families VI, VII and VIII. Again, this leads to the hypergeometric differential equation $_2F_1(\frac{1}{8},\frac{1}{8}, \frac{3}{4})$.
\end{Family}

\begin{Family}{X}{$\Gamma_{25}a_1$}
The family
$$X_\lambda : x_0^4+x_1^4+x_2^4+x_3^4 + 2 \lambda (x_0^2x_1^2 + x_2^2x_3^2)=0$$
degenerates where $\lambda^2 = 1.$

The Picard--Fuchs operator for $X_\lambda$ is
$$
(\lambda^2-1)^2\frac{d^3}{d\lambda^3} + 6\lambda(\lambda^2-1)\frac{d^2}{d\lambda^2} + (7\lambda^2-3)\frac{d}{d\lambda} + \lambda
$$
which is the symmetric square of
$$
(\lambda^2-1)\frac{d^2}{d\lambda^2} + 2\lambda\frac{d}{d\lambda} + \frac{1}{4}.
$$
After substituting $\mu = \lambda^2$, we get
$$
\mu(1-\mu)\frac{d^2}{d\mu^2} + \left(1-\frac{3}{2}\mu\right)\frac{d}{d\mu} - \frac{1}{16}
$$
which is the hypergeometric differential equation $_2F_1(\frac{1}{4},\frac{1}{4}, 1)$.
\end{Family}

\begin{Family}{XI}{$2^4D_6$}
The family
$$X_\lambda :  x_0^4+x_1^4+x_2^4+x_3^4 + \lambda(x_0^2+x_1^2+x_2^2+x_3^2)^2=0$$
degenerates at $\lambda = \infty$ and where $(\lambda+1)(\lambda+\frac{1}{2})(\lambda+\frac{1}{3})(\lambda+\frac{1}{4})=0$

The Picard--Fuchs operator for $X_\lambda$ is
 
\begin{eqnarray*}
(\lambda+1)(\lambda+\frac{1}{2})(\lambda+\frac{1}{3})(\lambda+\frac{1}{4})&\frac{d^3}{d\lambda^3} & + \\
\left(6 \lambda^3+\frac{75}{8}\lambda^2+\frac{35}{8}\lambda+\frac{5}{8}\right)&\frac{d^2}{d\lambda^2} & + \\
\left(\frac{27}{4}\lambda^2+7\lambda+\frac{13}{8}\right)&\frac{d}{d\lambda} & + \\
\frac{3}{4}\lambda+\frac{3}{8}
\end{eqnarray*}
which is the symmetric square of
$$
p(\lambda)\frac{d^2}{d\lambda^2} +
\frac{1}{2}p^\prime(\lambda)\frac{d}{d\lambda} +
\left(\frac{3}{16}\lambda^2+\frac{3}{16}\lambda+\frac{1}{24}\right)
$$
where $p(\lambda)=(\lambda+1)(\lambda+\frac{1}{2})(\lambda+\frac{1}{3})(\lambda+\frac{1}{4})$.

\end{Family}

\begin{Family}{XII}{$C_2\times{\mathfrak{S}_4}$}
The family
$$X_\lambda : p + 2\omega\lambda q^2=0$$
where
\begin{eqnarray*}
p & = & x_0^4 + x_1^4 + x_2^4 + x_3^4 - 2i\sqrt{3}(x_0^2x_1^2 + x_2^2x_3^2)\\
q & = & x_0x_2+x_1x_3-i x_1x_2-i x_0x_3
\end{eqnarray*}
and $\omega$ is a primitive 12th root of unity degenerates at $\lambda=\infty$ and where $(\lambda^2-\frac{1}{4})(\lambda^2-\frac{1}{3})$.

The Picard--Fuchs operator for $X_\lambda$ is
$$
(4\lambda^2-1)(3\lambda^2-1)\frac{d^3}{d\lambda^3} + 3\lambda(24\lambda^2-7)\frac{d^2}{d\lambda^2} + (81\lambda^2-9)\frac{d}{d\lambda} + 9\lambda
$$
which is the symmetric square of
$$
(4\lambda^2-1)(3\lambda^2-1)\frac{d^2}{d\lambda^2} + \lambda(24\lambda^2-7)\frac{d}{d\lambda} + \frac{1}{4}(9\lambda^2-2).
$$
After substituting $\mu = \lambda^2$, we get
$$
p(\mu)\frac{d^2}{d\mu^2} + \frac{1}{2}p^\prime(\mu)\frac{d}{d\mu} + \frac{3}{64}\mu-\frac{1}{96}
$$
with $p(\mu) = \mu\left(\mu-\frac{1}{4}\right)\left(\mu-\frac{1}{3}\right)$.

After the substitution $\mu = \frac{\nu+7}{36}$, this is the Lam\'{e} differential equation with $p(\nu) = 4(\nu+7)(\nu-2)(\nu-5)$, $n = -\frac{1}{4}$, and $B = \frac{3}{16}$.
\end{Family}

\begin{Family}{XIII}{$\mathfrak{A}_5\subset M_{20}$}
The family
$$
(5+15i)p+(1-3i)q^2 + \lambda q^2 = 0
$$
with
\begin{eqnarray*}
	p & = & x_0^4 + x_1^4 + x_2^4 + x_3^4 + 12x_0x_1x_2x_3\\
	q & = & x_0^2 - x_1^2 + x_2^2 + x_3^2 + (i + 1)(- x_0x_1 + ix_0x_2 + ix_0x_3 + x_1x_2 + x_1x_3 - ix_2x_3)
\end{eqnarray*}
degenerates at $\lambda = \infty$ and at the roots of $\lambda(\lambda-3)(\lambda+5)(3\lambda+16)$.
The Picard--Fuchs operator for this family is
 
\begin{eqnarray*}
\lambda(\lambda-3)(\lambda+5)(\lambda+\frac{16}{3})&\frac{d^3}{d\lambda^3} & + \\
\left(6\lambda^3+33\lambda^2-13\lambda-120\right)&\frac{d^2}{d\lambda^2} & + \\
\left(\frac{27}{4}\lambda^2+25\lambda-7\right)&\frac{d}{d\lambda} & + \\
\frac{3}{4}\lambda+\frac{3}{2}
\end{eqnarray*}
which is the symmetric square of
$$
p(\lambda)\frac{d^2}{d\lambda^2} +
\frac{1}{2}p^\prime(\lambda)\frac{d}{d\lambda} +
\left(\frac{3}{16}\lambda^2+\frac{3}{4}\lambda-\frac{2}{3}\right)
$$
where
$$
p(\lambda) = \lambda(\lambda-3)(\lambda+5)(\lambda+\frac{16}{3}).
$$
\end{Family}

\begin{Family}{XIV}{$\mathfrak{A}_5\subset\mathfrak{S}_5$}
The family
$$X_\lambda : p+\lambda{q}=0$$
where
$$p := 5x_0^4+ 5x_1^4 - 6x_0^2x_1^2 + 12x_2^2x_3^2 - 8\sqrt{2}(x_0x_2^3  + x_1x_3^3) - 48x_0x_1x_2x_3,$$
and
$$q := (2x_0x_1 + 3x_2x_3)(x_0^2 + x_1^2) - 2\sqrt{2}(x_0x_3^3 +x_1x_2^3)$$
degenerates where $(\lambda^2+80)(\lambda^2 - 1)=0.$

The Picard--Fuchs operator for $X_\lambda$ is
$$
(\lambda^2-1)(\lambda^2+80)\frac{d^3}{d\lambda^3} + 9\lambda(\lambda^2+26)\frac{d^2}{d\lambda^2} + \frac{3(6\lambda^4+403\lambda^2+2240)}{(\lambda^2+80)}\frac{d}{d\lambda} + \frac{3\lambda(2\lambda^2+85)}{(\lambda^2+80)}
$$
which is the symmetric square of
$$
(\lambda^2-1)(\lambda^2+80)\frac{d^2}{d\lambda^2} + 3\lambda(\lambda^2+26)\frac{d}{d\lambda} + \frac{3(\lambda^4+85\lambda^2+160)}{4(\lambda^2+80)}.
$$
After substituting $\mu = \lambda^2$, we get
$$
p(\mu)\frac{d^2}{d\mu^2} + \frac{1}{2}p^\prime(\mu)\frac{d}{d\mu} + \frac{3}{16}(\mu^2+85\mu+160)
$$
where $p(\mu) = \mu(\mu-1)(\mu+80)^2$.

\end{Family}

\section{Families With The Same Picard--Fuchs Equation}\label{wps}
Our main aim for this section is to prove the following.
\begin{theorem}\label{thm::samePF}
If $X_\lambda$, $\lambda\in\CC$, is a family of algebraic K3 surfaces and $G$ a finite group of symplectic automorphisms, then the quotient family of K3 surfaces $\widetilde{X_\lambda/G}$ will have the same Picard--Fuchs differential equation as $X_\lambda$.
\end{theorem}
In order to prove this, we recall that both of these families have lattice polarisations. The family $\widetilde{X_\lambda/G}$ is the minimal resolution of the singular surfaces $X_\lambda/G$ and the exceptional curves generate a negative definite lattice $E\subset\textnormal{Pic}(\widetilde{X_\lambda/G})$. Being a family of projective K3 surfaces, there is also a positive element $\kappa^\prime\in\textnormal{Pic}(\widetilde{X_\lambda/G})$. Because $E$ is negative definite, we may assume $\kappa^\prime\in E^\perp$. The smallest primitive lattice containing $E\oplus\kappa^\prime$ polarises the family $\widetilde{X_\lambda/G}$.

We also recall the definitions
$$T_{X_\lambda,G} = \textnormal{H}^2(X_\lambda,\ZZ)^G$$
and
$$S_{X\lambda,G} = T_{X_\lambda,G}^\perp.$$
Let $\mathfrak{X}$ denote the total space of the family $X_\lambda$ and write
\begin{eqnarray*}
&&\mathfrak{X}\\
&&\downarrow \pi\\
\widetilde{\mathfrak{X}/G} & \stackrel{\sigma}{\longrightarrow} & \mathfrak{X}/G.
\end{eqnarray*}
The map $\tau := \sigma^{-1}\circ\pi:\mathfrak{X}\rightarrow\widetilde{\mathfrak{X}/G}$ is a rational map which is a $|G|$ to 1 cover of its image wherever it is defined. Over each fibre, these maps restrict as follows:
\begin{eqnarray*}
&&X_\lambda\\
&\stackrel{\tau}{\swarrow}&\downarrow \pi\\
\widetilde{X_\lambda/G} & \stackrel{\sigma}{\longrightarrow} & X_\lambda/G.
\end{eqnarray*}
The pull-back $\tau^*:\textnormal{H}^2(\widetilde{X_\lambda/G},\ZZ)\rightarrow T_{X_\lambda,G}\subset\textnormal{H}^2(X_\lambda,\ZZ)$ satisfies
\begin{equation}\label{eqn:multiply}
\left<\tau^*(y_1),\tau^*(y_2)\right> = |G|\left<y_1,y_2\right>
\end{equation}
for all $y_1,y_2\in E^\perp$ (see for example \cite{NikAut}).

Hence $\kappa^\prime$ pulls back to a positive element $\kappa = \tau^*(\kappa^\prime)\in T_{X_\lambda,G}$. Similarly, a non-zero holomorphic 2-form $\omega_{\widetilde{X_\lambda/G}}\in\textnormal{H}^{2,0}(\widetilde{X_\lambda/G})\subset\textnormal{H}^2(\widetilde{X_\lambda/G},\ZZ)\otimes\CC$ pulls back to a 2-form $\omega_{X_\lambda} := \tau^*(\omega_{\widetilde{X_\lambda/G}})$ that must also be non-zero and holomorphic.

Since
$$
\left<\kappa,\omega_{X_\lambda}\right> = \left<\tau^*(\kappa^\prime),\tau^*(\omega_{\widetilde{X_\lambda/G}})\right> = |G|\left<\kappa^\prime,\omega_{\widetilde{X_\lambda/G}}\right> = 0
$$
and similarly $\left<\kappa,\overline{\omega}_{X_\lambda}\right> = 0$, we have $\kappa\in\textnormal{H}^2(X_\lambda,\ZZ)^G\cap\textnormal{H}^{1,1}(X_\lambda)\subset\textnormal{Pic}(X_\lambda)$ and $X_\lambda$ is polarised by the lattice $S_{X_\lambda,G}\oplus\kappa$.

By (\ref{eqn:multiply}), $\tau^*$ injects $E^\perp$ into $T_{X_\lambda,G} = \textnormal{H}^2(X_\lambda,\ZZ)^G$. Therefore, if $\{\gamma_0,\ldots,\gamma_k\}$ is a basis for the lattice $(E\oplus\kappa^\prime)^\perp$, then $\{\tau^*(\gamma_0),\ldots,\tau^*(\gamma_k)\}$ will be a basis for the $\QQ$--vector space $(S_{X_\lambda,G}\oplus\kappa)^\perp\otimes\QQ.$ There is therefore a $\QQ$--linear transformation
$$
P : (S_{X_\lambda,G}\oplus\kappa)^\perp\otimes\QQ \rightarrow (S_{X_\lambda,G}\oplus\kappa)^\perp\otimes\QQ
$$
such that $\{P(\tau^*(\gamma_i))\}$ is a basis of the lattice $(S_{X_\lambda,G}\oplus\kappa)^\perp$.

To prove theorem \ref{thm::samePF}, with respect to the 2--forms $\omega_{\widetilde{X_\lambda/G}}$ and $\omega_{X_\lambda} = \tau^*(\omega_{\widetilde{X_\lambda/G}})$, we consider the periods as the linear maps
\begin{eqnarray*}
p_{\widetilde{X_\lambda/G}} : \textnormal{H}^2(\widetilde{X_\lambda/G},\ZZ) & \rightarrow & \CC\\
\gamma & \mapsto & \left<\omega_{\widetilde{X_\lambda/G}},\gamma\right> = \int_{\widetilde{X_\lambda/G}}\omega_{\widetilde{X_\lambda/G}}\wedge\gamma
\end{eqnarray*}
and
\begin{eqnarray*}
p_{X_\lambda} : \textnormal{H}^2(X_\lambda,\ZZ) & \rightarrow & \CC\\
\delta & \mapsto & \left<\omega_{X_\lambda},\delta\right> = \int_{X_\lambda}\omega_{X_\lambda}\wedge\delta.
\end{eqnarray*}
Then by (\ref{eqn:multiply}) we have
$$
p_{X_\lambda}(\tau^*(\gamma)) = |G|.p_{\widetilde{X_\lambda/G}}(\gamma)
$$
for any $\gamma\in E^\perp$. Thus the projective periods of $X_\lambda$ and $\widetilde{X_\lambda/G}$ are related by the linear transformation $P$. The differential equation satisfied by these periods will therefore be the same.
\begin{flushright}
$\square$
\end{flushright}

Families I, II and III each have the hypergeometric differential equation $_2F_1(\frac{1}{12},\frac{1}{12},\frac{2}{3})$ as their Picard--Fuchs equation. Similarly, families VI, VII, VIII and IX have Picard--Fuchs equations $_2F_1(\frac{1}{8},\frac{1}{8},\frac{3}{4})$. In light of theorem \ref{thm::samePF}, these coincidences can be explained if we can find isomorphisms between the quotients of these families by finite groups of symplectic automorphisms. 

Before we look at this, we introduce some more examples of families of K3 surfaces. Whereas our families I---XIV are subfamilies of the 19 dimensional moduli spaces of K3 surfaces in $\PP^3$ and $\PP(1,1,1,3)$, we shall now look at weighted projective spaces that provide a 1--dimensional moduli space of K3 hypersurfaces. The algorithm of section \ref{det} and appendix \ref{macaulay} that finds the Picard--Fuchs differential equations is valid for quasismooth families of hypersurfaces in weighted projective space. Out of the 95 weighted projective spaces containing K3 hypersurfaces, only
$$
\PP(7,8,9,12),\ \ \PP(7,8,10,25)\ \ \textnormal{and}\ \ \PP(5,6,8,11)
$$
contain exactly 1--dimensional families.

\begin{example}
A hypersurface in $\PP(7,8,9,12)$ is a K3 surface with DuVal singularities whenever it is defined by a nonsingular polynomial of weighted degree $36 = 7+8+9+12$. If $x,y,z,t$ are the coordinates with weights 7, 8, 9 and 12, then the only weight 36 monomials are
$$
x^4y,\ \ y^3t,\ \ x^4,\ \ t^3\ \ \textnormal{and}\ \ xyzt.
$$
Up to a weighted projective transformation of the form
$$
(x,y,z,t)\mapsto(\mu_xx,\mu_yy,\mu_zz,\mu_tt)
$$
for $\mu_x,\mu_y,\mu_z,\mu_t\in\CC^*$, any degree 36 polynomial is of the form
$$
X_\lambda: (x^4y + y^3t + x^4 + t^3 + 4\lambda xyzt = 0)\subset\PP(7,8,9,12).
$$
This family degenerates at $\lambda=\infty$ and where $\lambda^4=1$. Indeed, the transformation
$$
(x,y,z,t)\mapsto(x,y,iz,t)
$$
provides an isomorphism $X_\lambda\cong X_{i\lambda}$. The algorithm of appendix \ref{macaulay} finds that the Picard--Fuchs differential operator for this family is
$$
(\lambda^4-1)\frac{d^3}{d\lambda^3} + 6\lambda^3\frac{d^2}{d\lambda^2} + 7\lambda^2\frac{d}{d\lambda} + \lambda.
$$
Substituting $\mu = \lambda^4$, we get the hypergeometric differential equation $_2F_1(\frac{1}{8},\frac{1}{8},\frac{3}{4})$. This has occurred before as the Picard--Fuchs equation for the families in $\PP^3$ invariant under the groups $C_7\rtimes C_3$, $4^2\mathfrak{A}_4$, $SD_{16}$ and $\textnormal{Hol}(C_5)$ (examples VI, VII, VIII and IX). These families are listed in table \ref{tab:FamiliesWithTheSamePicardFuchsEquation}.

\begin{table}
	\newcommand\T{\rule{0pt}{2.6ex}}
  \newcommand\B{\rule[-1.2ex]{0pt}{0pt}}
	\centering
		\begin{tabular}{|c|c|}
			\hline
			\T\B Example&Family\\
			\hline\hline
			\T\B $\PP(7,8,9,12)$ & $x^4y + y^3t + x^4 + t^3 + 4\lambda xyzt = 0$\\\hline
			\T\B VI: $C_7\rtimes C_3$ & $x_0x_1^3 + x_1x_2^3 + x_2x_0^3 + x_3^4 + 4 \lambda x_0x_1x_2x_3 = 0$\\\hline
			\T\B VII: $4^2\mathfrak{A}_4$ & $x_0^4 + x_1^4 + x_2^4 + x_3^4 + 4\lambda x_0x_1x_2x_3 = 0$\\\hline
			\T\B VIII: $SD_{16}$ & $x_0^4+x_1^4+x_2^4+x_3^4 + 2\sqrt{2}\lambda x_0x_1(x_2^2+ix_3^2)=0$\\\hline
			\T\B IX: $\textnormal{Hol}(C_5)$ & $x_0^3x_1 + x_1^3x_3 + x_3^3x_2 - x_2^3x_0 + \left(\frac{1+i}{2}\right)\lambda(x_0^2x_3^2 + x_1^2x_2^2)=0$\\
			\hline
		\end{tabular}
	\caption{Families with the Picard--Fuchs Equation $_2F_1(\frac{1}{8},\frac{1}{8},\frac{3}{4})$}
	\label{tab:FamiliesWithTheSamePicardFuchsEquation}
\end{table}

We expect that if two families have the same Picard--Fuchs equation, then they are likely to be related geometrically. For example, they may be isomorphic, or perhaps one family is a quotient of another by some group action. In the present case, we note that the family invariant under the group $4^2\mathfrak{A}_4$ is also invariant under the Abelian subgroup $C_4\times C_4$ generated by the projective transformations
\begin{eqnarray*}
	(x_0,x_1,x_2,x_3) & \mapsto & (ix_0,-ix_1,x_2,x_3)\\
	(x_0,x_1,x_2,x_3) & \mapsto & (x_0,ix_1,-ix_2,x_3)\\
	(x_0,x_1,x_2,x_3) & \mapsto & (x_0,x_1,ix_2,-ix_3)\\
	(x_0,x_1,x_2,x_3) & \mapsto & (-ix_0,x_1,x_2,ix_3).
\end{eqnarray*}
The  monomials $x_0^4,\ x_1^4,\ x_2^4,\ x_3^4$ and $x_0x_1x_2x_3$ generate the graded ring of invariants of this group action. The family
$$
F_\lambda: ( x_0^4 + x_1^4 + x_2^4 + x_3^4 + 4\lambda x_0x_1x_2x_3 = 0 )\subset\PP^3
$$
is invariant under $C_4^2$ and the general quotient surface $F_\lambda/{C_4^2}$ has 6 $A_3$ singularities (see for example \cite{Xiao}). On the minimal resolution, there are 18 exceptional curves and 4 other rational curves coming from the images of the curves defined by $(x_i = 0)$ for $i = 0, 1, 2$, and $3$. Each of the singular points lie in the intersection between two of these 4 curves $x_i=0$, and the curves are arranged in the tetrahedral configuration of figure \ref{fig:C_4^2}. In this diagram, exceptional curves are shown as solid vertices and the other four rational curves are the non--solid vertices.
\begin{figure}
	\centering
		\includegraphics{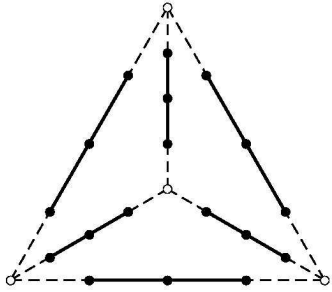}
	\caption{Configuration of Curves on $X/C_4^2$.}
	\label{fig:C_4^2}
\end{figure}

This configuration of curves can be used to indicate the existence of a birational morphism from another K3 surface. For example, it can be shown that the general degree  36 hypersurface in $\PP(7,8,9,12)$ has singular points $A_7, A_6, A_3$ and $A_2$ also arranged in a tetrahedral configuration as in figure~\ref{fig:P78912}.
\begin{figure}
	\centering
		\includegraphics{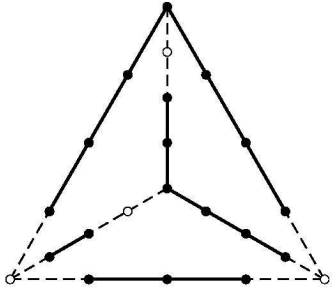}
	\caption{Configuration of Curves on $X_{36}\subset\PP(7,8,9,12)$.}
	\label{fig:P78912}
\end{figure}

To compare this weighted projective space example to the family $F_\lambda$ invariant under the group $G=C_4^2\subset 4^2\mathfrak{A}_4$, notice that since the graded ring of invariants of this group action is
\begin{eqnarray*}
    \CC[x_0,x_1,x_2,x_3]^G & = & \CC[x_0^4, x_1^4, x_2^4, x_3^4, x_0x_1x_2x_3]\\
                           & \cong & \CC[X_0,X_1,X_2,X_3,X_4]/(X_0X_1X_2X_3-X_4^4)
\end{eqnarray*}
then
$$
\PP^3/G \cong \Proj(\CC[x_0,x_1,x_2,x_3]^G) \cong (X_0X_1X_2X_3=X_4^4)\subset\PP^4.
$$
For a discussion of quotients constructed from invariants, see section~\ref{GIT}. The quotient surfaces may be described as
$$
F_\lambda / C_4^2 \cong (X_0+X_1+X_2+X_3+4\lambda X_4 = 0 ,\hspace{6pt} X_0X_1X_2X_3=X_4^4)\subset\PP^4.
$$
\begin{proposition}\label{78912isom}
The minimal resolutions of the families $F_\lambda / C_4^2$ and $X_{36}\subset\PP(7,8,9,12)$ are isomorphic.
\end{proposition}
\begin{proof}It is worth noticing that by the mirror symmetry construction of Batyrev, the family $F_\lambda / C_4^2$ is the mirror of the general quartic in $\PP^3$. As a result, its Picard lattice is known to be $(E_8)^2\oplus{U}\oplus\left<-4\right>$ and the transcendental lattice is $U\oplus \left<4\right>$. By Belcastro \cite{Bel}, the Picard and transcendental lattices of $X_{36}$ are the same. By a result of Nikulin (Theorem 1.14.4 in \cite{Nik}), the embedding of this Picard lattice into the K3 lattice is unique and we should expect that these families of K3 surfaces are isomorphic.

In fact, we shall write down an explicit birational morphism between these two families. This will determine a birational morphism between their minimal resolutions. Since any birational morphism between minimal nonruled nonsingular surfaces is in fact biregular, we obtain our isomorphism.

The morphism
\begin{eqnarray*}
    \varphi_{|-K_\PP|} : \PP(7,8,9,12) & \rightarrow & \PP^4\\
    (x,y,z,t)     & \mapsto     & (x^4y, y^3t, z^4, t^3, xyzt)
\end{eqnarray*}
defined by the anticanonical linear system on $\PP(7,8,9,12)$ restricts to $X_{36}$ to give a well defined birational morphism onto
$$\left(\sum_{i=0}^3{X_i}=4\lambda X_4 ,\hspace{6pt} X_0X_1X_2X_3=X_4^4\right)\subset\PP^4.$$
It is easily checked that this has a two--sided inverse given by
$$(X_0,X_1,X_2,X_3,X_4)\mapsto(X_0^2X_1^4X_3,\hspace{3pt} X_0^2X_1^5X_3,\hspace{3pt} X_0^2X_1^5X_3X_4,\hspace{3pt} X_0^3X_1^7X_3^2)$$
so that $\varphi_{|-K_\PP|}|_{X_{36}}$ is the required birational morphism.
\end{proof}

So these two families are related in a concrete way. Similarly, the group $C_7\rtimes C_3$ contains an element of order 7 generating the subgroup $C_7$. Under the representation given in chapter \ref{autChapter}, the order 7 element is
$$
\begin{pmatrix}
\zeta&0&0&0\\
0&\zeta^2&0&0\\
0&0&\zeta^4&0\\
0&0&0&1
\end{pmatrix}.
$$
 On an invariant quartic K3 surface, $X$, this order 7 symplectic automorphism has three fixed points leading to three $A_6$ singularities on the quotient $X/C_7$. In this example, it can be shown that the curves also form a tetrahedral configuration shown in figure~\ref{fig:C_7} and so we should expect a similar result to that of proposition~\ref{78912isom}. Indeed this is the case, although the isomorphism is a little less clear because the ring of invariants for this action of $C_7$ is not a complete intersection ring and so has numerous generators and relations. However, amongst the invariants are the monomials
$$
x_0x_1^3,\ x_1x_2^3,\ x_2x_0^3,\ x_3^4,\ x_0x_1x_2x_3 \in \CC[x_0,x_1,x_2,x_3]^{C_7}.
$$
satisfying a relationship of the form $X_0X_1X_2X_3=X_4^4$ so in this case, the quotient family projects onto the previous quotient family $F_\lambda/C_4^2$.
\begin{figure}
	\centering
		\includegraphics{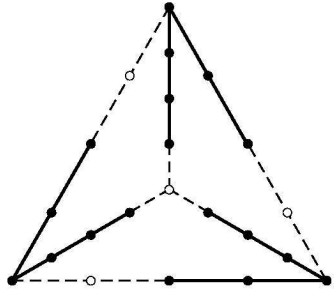}
	\caption{Configuration of Curves on $X/C_7$.}
	\label{fig:C_7}
\end{figure}

For the families XIII and IX invariant under the groups $SD_{16}$ and $\textnormal{Hol}(C_5)$ also with the same differential equation, the situation is less clear. It would be interesting to find an isomorphism as before, but if one exists it is probably disguised by an unsuitable choice of defining equations. Certainly, $SD_{16}$ has an element of order 8, and with this Abelian group, the quotient family has the suggestive configuration of curves given in figure \ref{fig:C_8}. The pattern so far has been to take the quotient of the family by its largest Abelian subgroup of symplectic automorphisms. The largest Abelian subgroup of $\textnormal{Hol}(C_5)$ is the cyclic group of order 5. This order 5 automorphism has 4 fixed points leading to 4 $A_4$ singularities on the quotient. To fit these exceptional curves in to a similar tetrahedral configuration would leave 6 curves non--exceptional to be accounted for rather than 4.
\begin{figure}
	\centering
		\includegraphics{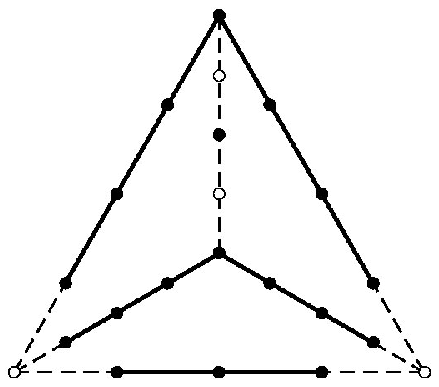}
	\caption{Configuration of Curves on $X/C_8$.}
	\label{fig:C_8}
\end{figure}
\end{example}

We look at another example of families of K3 surfaces with the same Picard--Fuchs equation as each other.

\begin{example}
In $\PP(7,8,10,25)$, the weight 50 monomials are
$$
x^6y,\ y^5z,\ z^5,\ t^2,\ xyzt\ \textnormal{and}\ x^2y^2z^2.
$$
This time we have an additional weighted projective transformation of the form
$$t\mapsto t+\mu xyz$$
and any K3 hypersurface is isomorphic to one of the form
$$
X_\lambda:(t^2 = x^6y+y^5z+z^5-3\lambda x^2y^2z^2)\subset\PP(7,8,10,25).
$$
This family degenerates at $\lambda = \infty$ and where $\lambda^3=1$ and there is an isomorphism $X_\lambda\cong X_{\omega\lambda}$ with $\omega^3=1$ that is induced by the transformation
$$
(x,y,z,t)\mapsto(\omega^2 x,y,z,t).
$$
The algorithm of appendix \ref{macaulay} shows that the Picard--Fuchs differential operator is
$$
(\lambda^3-1)\frac{d^3}{d\lambda^3} + \frac{9}{2}\lambda^2\frac{d^2}{d\lambda^2} + \frac{13}{4}\lambda\frac{d}{d\lambda} + \frac{1}{8}.
$$
This is the same as the Picard--Fuchs equation occurring in examples I, II and III and with respect to the parameter $\mu = \lambda^3$, the symmetric square root is the hypergeometric differential equation $_2F_1(\frac{1}{12},\frac{1}{12},\frac{2}{3})$.

We perform the same calculation for $\PP(5,6,8,11)$, but start by reordering the weights for later convenience as $\PP(5,8,6,11)$. Writing $X$, $Y$, $Z$ and $T$ for the weights of degrees 5, 6, 8 and 11, the degree 30 monomials are
$$
X^6,\ Y^3Z,\ Z^5,\ T^2Y,\ XYZT,\ \textnormal{and}\ X^2YZ^2.
$$
Taking into account weighted projective transformations, especially the transformation
$$T\mapsto T+\mu XZ,$$
a K3 surface in this space is a member of the family
$$
(T^2Y = X^6 + Y^3Z + Z^5 -3\lambda X^2YZ^2) \subset \PP(5,8,6,11).
$$
The calculation of the Picard--Fuchs equation for this family is the same as for the family in $\PP(7,8,10,25)$ and again we obtain the hypergeometric differential equation $_2F_1(\frac{1}{12},\frac{1}{12},\frac{2}{3})$.

\begin{table}
	\newcommand\T{\rule{0pt}{2.6ex}}
  \newcommand\B{\rule[-1.2ex]{0pt}{0pt}}
	\centering
		\begin{tabular}{|c|c|c|}
			\hline
			\T\B Example&Family&Name\\
			\hline\hline
			\T\B $\PP(7,8,10,25)$ & $t^2 = x^6y+y^5z+z^5-3\lambda x^2y^2z^2$&$\mathcal{A}$\\\hline
			\T\B $\PP(5,8,6,11)$ & $T^2Y = X^6 + Y^3Z + Z^5 -3\lambda X^2YZ^2$&$\mathcal{B}$\\\hline
			\T\B I: $SD_{16}\subset T_{48}$ & $t^2 = x_0^6 + x_1^5x_2 + x_1x_2^5 - 3 \lambda x_0^2x_1^2x_2^2$&$\mathcal{C}$\\\hline
			\T\B II: $\mathfrak{A}_{4,3}$ & $t^2 = x_0^6+x_1^6+x_2^6-3\lambda x_0^2x_1^2x_2^2$&$\mathcal{D}$\\\hline
			\T\B III: $C_7\rtimes{C_3} \subset L_2(7)$ & $t^2 = x_0^5x_1  + x_1^5x_2 + x_2^5x_0 - 3\lambda x_0^2x_1^2x_2^2$&$\mathcal{E}$\\\hline
		\end{tabular}
	\caption{Families with the Picard--Fuchs Equation $_2F_1(\frac{1}{12},\frac{1}{12},\frac{2}{3})$}
	\label{tab:MoreFamiliesWithTheSamePicardFuchsEquation}
\end{table}

With reference to the labels for the families given in table~\ref{tab:MoreFamiliesWithTheSamePicardFuchsEquation}, in these examples, there are  isomorphisms between the following minimal resolutions:
$$
\mathcal{A}\cong \mathcal{B}\cong \mathcal{C}/C_8\cong \mathcal{D}/(C_2\times{C_6}) \cong \mathcal{E}/C_7
$$
where $C_8\in SD_{16}$, $C_2\times C_6 \in \mathfrak{A}_{4,3}$ and $C_7\in C_7\rtimes C_3$ are the largest Abelian subgroups of the symmetry groups.
 
For $\PP(7,8,10,25)$, the anticanonical morphism
\begin{eqnarray*}
    \varphi_{\left|-K_1\right|} : \PP(7,8,10,25) & \rightarrow & \PP^5 \\
                (x,y,z,t)   & \mapsto     & (x^6y,y^5z,z^5,x^2y^2z^2,t^2,xyzt)\\
                            &             &=: (x_0,x_1,x_2,x_3,x_4,x_5)
\end{eqnarray*}
has image $(x_0x_1x_2 = x_3^3,\ x_3x_4 = x_5^2)\subset\PP^5$. The inverse is the map
$$(x_0,x_1,x_2,x_3,x_4,x_5)\mapsto (x_0x_1^5x_2^7x_3,\ x_1^5x_2^7x_3^4,\ x_1^6x_2^9x_3^5,\ x_1^{15}x_2^{22}x_3^{12}x_5).$$

For $\PP(5,8,6,11)$ the anticanonical morphism is
\begin{eqnarray*}
    \varphi_{\left|-K_2\right|}: \PP(5,8,6,11) & \rightarrow & \PP^5 \\
                (X,Y,Z,T)   & \mapsto     & (X^6,Y^3Z,Z^5,X^2YZ^2,YT^2,XYZT)\\
                            &             &=: (x_0,x_1,x_2,x_3,x_4,x_5)
\end{eqnarray*}
and also has image $(x_0x_1x_2 = x_3^3,\ x_3x_4 = x_5^2)\subset\PP^5$. The inverse is the map
$$(x_0,x_1,x_2,x_3,x_4,x_5)\mapsto (x_0x_1^5x_2^4,\ x_0x_1^8x_2^6x_3,\ x_0x_1^6x_2^5,\ x_0x_1^{10}x_2^8x_3^2x_5).$$

The composition $\varphi_{\left|-K_2\right|}^{-1}\circ\varphi_{\left|-K_1\right|}^{\phantom{-1}}\colon\PP(7,8,10,25)\rightarrow\PP(5,8,6,11)$ is the birational morphism
$$(x,y,z,t) \mapsto (xy,y^3,yz,yt)$$
with inverse
$$(X,Y,Z,T) \mapsto (XY^2,Y^3,Y^3Z,Y^8T).$$

This map restricts to induce an isomorphism between the minimal resolutions of the families of K3 surfaces in these two spaces.

The other isomorphisms follow in a similar manner. For example, the group $C_2\times C_6$ is generated by the projective transformations
\begin{eqnarray*}
(t,x_0,x_1,x_2) & \mapsto & (t,x_0,\upsilon x_1,\upsilon^5 x_2)\\
(t,x_0,x_1,x_2) & \mapsto & (t,\upsilon x_0,\upsilon^5 x_1,x_2)\\
(t,x_0,x_1,x_2) & \mapsto & (t,\upsilon^5 x_0,x_1,\upsilon x_2)
\end{eqnarray*}
where $\upsilon^6=1$. The ring of invariants is
\begin{eqnarray*}
\CC[t,x_0,x_1,x_2]^{C_2\times C_6}  & = & \CC[t,x_0^6,x_1^6,x_2^6,x_0x_1x_2]\\
																		& \cong & \CC[T,X_0,X_1,X_2,Y]/(X_0X_1X_2 = Y^6)
\end{eqnarray*}
so that the quotient is given by
$$
\PP(1,1,1,3)/(C_2\times C_6) \cong (X_0X_1X_2 = Y^6) \subset \PP(3,6,6,6,3).
$$
This may be straightened by the anticanonical morphism
\begin{eqnarray*}
\PP(3,6,6,6,3) & \rightarrow & \PP^5\\
(T,X_0,X_1,X_2,Y) & \mapsto & (X_0,X_1,X_2,Y^2,T^2,YT)
\end{eqnarray*}
with image $(x_0x_1x_2 = x_3^3,\ x_3x_4 = x_5^2)\subset\PP^5$ establishing another isomorphism of K3 surfaces.

Examining the invariants of $C_8$ defined by
$$
(t,x_0,x_1,x_2) \mapsto (t,\nu x_0,\nu^3 x_1,\nu^4 x_2)
$$
with $\nu^8=1$ and of $C_7$ defined by
$$
(t,x_0,x_1,x_2) \mapsto (t,\zeta x_0,\zeta^2 x_1,\zeta^4 x_2)
$$
with $\zeta^7=1$ leads to similar isomorphisms.

\begin{remark}
Since all weighted projective spaces are toric, they are all rational and in particular, any two weighted projective spaces of the same dimension are birationally equivalent. If, furthermore, a birational equivalence $f$ preserves the anticanonical linear system (i.e., $f^*\left|-K_b\right| = \left|-K_a\right|$), then the anticanonical models coincide:
\begin{eqnarray*}
\PP(a_0,\ldots,a_n) & \stackrel{f}{\longrightarrow} & \PP(b_0,\ldots,b_n) \\
     \varphi_{\left|-K_a\right|} \searrow      &             &   \swarrow \varphi_{\left|-K_b\right|}  \\
                    &  \PP^m       &
\end{eqnarray*}
Since families of K3 hypersurfaces in weighted projective space correspond to the space's anticanonical linear system, we can detect whether two families are birationally equivalent by checking whether the anticanonical models of the weighted projective spaces are isomorphic.

This method can be used in more examples than just those examined here. For examples, by examining their anticanonical models, it can be shown that the weighted projective spaces $\PP(5,6,22,33)$ and $\PP(3,5,11,14)$ both contain birational families of K3 surfaces with generic Picard number 18.
\end{remark}

\end{example}

	\newpage
	
\chapter{Monodromy}\label{mono}
\section{How to Calculate the Monodromy Representation}
We tackle the problem of calculating the monodromy representation of a given Fuchsian differential system. Corollary \ref{disteigen} makes it easy to calculate the \emph{local} monodromies about the singular points. However, each local monodromy will be calculated with respect to a different basis and to find the global monodromy representation we usually need to analytically continue a fixed basis of solutions around each singular point in turn. This analytic continuation can only be carried out by numerical approximation which introduces problems when trying to find the monodromy rigorously. There is a situation where these problems can be avoided altogether -- the case where the local system of solutions is rigid.
\subsection{Rigid Systems - The Hypergeometric Case}
A second order Fuchsian ODE with three regular singular points can always be transformed into a hypergeometric ODE and, in this case, we can use the rigidity of the local system of solutions to determine the monodromy group. First, we state precisely what we mean by rigidity.
\begin{definition}
Let $A_1,\ldots,A_k\in\textnormal{Gl}(n,\CC)$ be matrices satisfying
\begin{equation}\label{eqn:prod}
	\prod_{i=1}^kA_i=id.
\end{equation}
Let $C_j$ denote the conjugacy class of $A_j$. We say that the group $\left<A_1,\ldots,A_k\right>$ is \emph{rigid} if whenever we choose representatives $B_j\in C_j$ satisfying $\prod_{i=1}^kB_i=id$, there is a fixed $P\in\textnormal{Gl}(n,\CC)$ with $B_j = P.A_j.P^{-1}$ for $j=1,\ldots,k$.
\end{definition}
We say that the Fuchsian system (\ref{system}) is rigid if its monodromy group is rigid. To calculate the monodromy group of a rigid Fuchsian system, it is enough to determine the local monodromies and then find some conjugates of these that satisfy (\ref{eqn:prod}). Furthermore, hypergeometric differential equations are known to be rigid thanks to Levelt's theorem (see \cite{Beu}):
\begin{theorem}[Levelt]\label{Levelt}
If $A, B\in\textnormal{Gl}(n,\CC)$ have eigenvalues $a_1,\ldots,a_n$ and $b_1,\ldots,b_n$ with
\begin{equation}\label{dist}
\{a_1,\ldots,a_n\}\cap\{b_1,\ldots,b_n\}=\emptyset
\end{equation}
and if $A.B^{-1}$ is a pseudo--reflection (meaning that $A.B^{-1}-id$ has rank 1), then $\left<A,B\right>$, the group generated by $A$ and $B$, is rigid.
\end{theorem}

The hypergeometric differential system (\ref{hypergeomSys}) with the same monodromy representation as the hypergeometric differential equation (\ref{hypergeom}) has residue matrices
\begin{eqnarray*}
R_0 & = & \begin{pmatrix}
0&0\\
-\alpha\beta&-\gamma
\end{pmatrix}\\
R_1 & = & \begin{pmatrix}
0&1\\
0&\gamma-\alpha-\beta
\end{pmatrix}\\
R_\infty & = & \begin{pmatrix}
0 & -1\\
\alpha\beta & \alpha+\beta
\end{pmatrix}.
\end{eqnarray*} 
\begin{proposition}
The local monodromy matrices of $_2F_1(\alpha,\beta,\gamma)$ only depend on the values of $\alpha$, $\beta$ and $\gamma$ modulo integers.
\end{proposition}
\begin{proof}
As long as the distinct eigenvalues of a residue matrix $R=R_0$, $R_1$ or $R_\infty$ do not differ by integers, the local monodromy is given by $e^{2 \pi i R}$. This expression is not effected be adding integers to $\alpha$, $\beta$ or $\gamma$. In the case where the eigenvalues do differ by a non--zero integer, the local monodromy can be found as a limit of local monodromies of systems of the general type. In this case, the invariance under addition of integers is preserved.
\end{proof}
Let $M_0, M_1$ and $M_\infty$ be the generators for the monodromy group of the hypergeometric differential system, and write $A=M_\infty$ and $B^{-1} = M_0$. Then $A.B^{-1} = M_1^{-1}$ is a pseudo--reflection and Levelt's theorem shows that the monodromy group is rigid whenever (\ref{dist}) is additionally satisfied. This happens precisely when
\begin{equation}\label{eqn:irred}
\alpha,\beta,\gamma-\alpha\ \textnormal{and}\ \gamma-\beta\not\in\ZZ.
\end{equation}

\begin{definition}
A subgroup $G\subset\textnormal{Gl}(n,\CC)$ is said to be \emph{irreducible} if it fixes no proper linear subspace of $\CC^n$.
\end{definition}
the monodromy group of $_2F_1(\alpha,\beta,\gamma)$ is irreducible if and only if condition (\ref{eqn:irred}) is satisfied (see \cite{Beu}). So the hypergeometric differential equation is rigid whenever its monodromy representation is irreducible.

\begin{corollary}\label{integers}
If the monodromy group of $_2F_1(\alpha,\beta,\gamma)$ is irreducible, ie. if none of the values $\alpha,\beta,\gamma-\alpha$ or $\gamma-\beta$ are integers, then the monodromy group is determined by the values of $\alpha$, $\beta$ and $\gamma$ modulo $\ZZ$.
\end{corollary}
Although the monodromy representation of the hypergeometric differential equation is undoubtedly known, it is difficult to find in the literature. For this reason, we calculate the monodromy group from scratch. We would also like to find all the parameters $\alpha$, $\beta$ and $\gamma$ for which the hypergeometric ODE $_2F_1(\alpha,\beta,\gamma)$ has a monodromy group satisfying theorem \ref{squares}.

Since irrational and strictly complex values of $\alpha$, $\beta$ and $\gamma$ lead to non--discrete monodromy groups, we only need to consider rational values of the parameters. Because of corollary \ref{integers}, to find all the irreducible monodromy groups, we restrict our attention to rational parameters $\alpha$, $\beta\in(0,1)$ and $\gamma\in[0,1)$ with $\alpha$, $\beta\neq\gamma$. From now, we assume these restrictions and split up monodromy calculation into cases.

\bfseries{Case}\rm \ \ $\gamma-\alpha-\beta\neq0,-1$.

If we additionally assume that $\gamma\neq0$ and $\alpha\neq\beta$, then the canonical forms of the residue matrices are
\begin{eqnarray*}
J_0 & = & \begin{pmatrix}
0&0\\
0&-\gamma
\end{pmatrix}\\
J_1 & = & \begin{pmatrix}
0&0\\
0&\gamma-\alpha-\beta
\end{pmatrix}\\
J_\infty & = & \begin{pmatrix}
\alpha&0\\
0&\beta
\end{pmatrix}
\end{eqnarray*}
and the local monodromies are
\begin{eqnarray*}
L_0 & = & e^{-\pi i \gamma}\begin{pmatrix}
e^{\pi i \gamma}&0\\
0&e^{-\pi i \gamma}
\end{pmatrix}\\
L_1 & = & e^{\pi i (\gamma-\alpha-\beta)}\begin{pmatrix}
e^{-\pi i (\gamma-\alpha-\beta)}&0\\
0&e^{\pi i (\gamma-\alpha-\beta)}
\end{pmatrix}\\
L_\infty & = & e^{\pi i (\alpha+\beta)}\begin{pmatrix}
e^{\pi i (\alpha-\beta)}&0\\
0&e^{\pi i (\beta-\alpha)}
\end{pmatrix}.
\end{eqnarray*}
By rigidity, the monodromy representation is determined by finding the unique (up to simultaneous conjugation) matrices $M_0$, $M_1$ and $M_\infty$ conjugate to the corresponding local monodromies and satisfying $M_0.M_1.M_\infty=id$. It can be verified that the matrices
$$\begin{array}[h]{r c l}
M_1 & = & e^{\pi i (\gamma-\alpha-\beta)}\begin{pmatrix}
\phantom{-}\cos((\gamma-\alpha-\beta)\pi) & \sin((\gamma-\alpha-\beta)\pi)\\
-\sin((\gamma-\alpha-\beta)\pi) & \cos((\gamma-\alpha-\beta)\pi)
\end{pmatrix}\\
M_\infty & = & \\
\multicolumn{3}{c}{e^{\pi i (\alpha+\beta)}\begin{pmatrix}
e^{-\pi i(\alpha-\beta)} & i\frac{e^{2\pi i(2\alpha+\beta)}+e^{2\pi i(\alpha+2\beta)}+e^{2\pi i(\alpha+\gamma)}+e^{2\pi i(\beta+\gamma)}-2e^{2\pi i(\alpha+\beta+\gamma)}-2e^{2\pi i(\alpha+\beta)}}{e^{\pi i(\alpha+\beta+2\gamma)}-e^{3\pi i(\alpha+\beta)}}\\
0&e^{\pi i(\alpha-\beta)}
\end{pmatrix}}\\
M_0 & = & (M_1.M_\infty)^{-1}
\end{array}$$
satisfy $M_i = P_i.L_i.P_i^{-1}$ for $i = 0, 1, \infty$ where
\begin{eqnarray*}
P_0 & = & \begin{pmatrix}
e(\gamma)-\frac{1}{2}(e(\beta)+e(\gamma-\alpha)) & \frac{i(e(\gamma)-1)(e(\gamma)+e(\alpha+\beta)-2e(\alpha))}{e(\alpha+\beta)-e(\gamma)}\\
\frac{1}{2}i(e(\beta)-e(\gamma-\alpha)) & e(\gamma)-1
\end{pmatrix}\\
P_1 & = & \begin{pmatrix}
-\frac{i}{2} & 1\\
-\frac{1}{2} & i
\end{pmatrix}\\
P_\infty & = & \begin{pmatrix}
a & -1/b\\
b & 0
\end{pmatrix}
\end{eqnarray*}
$a=i e(\frac{\alpha-\beta}{2})(2e(\alpha+\beta)-e(2\alpha+\beta)-e(\alpha+2\beta)-e(\gamma+\alpha)-e(\gamma+\beta)+2e(\alpha+\beta+\gamma))$
$b = (1-e(\alpha-\beta))(e(\gamma+\frac{\alpha+\beta}{2})-e(\frac{3}{2}(\alpha+\beta)))$

where, for brevity, we have written $e(x) := e^{2\pi i x}$. Hence the matrices $M_1$ and $M_\infty$ above generate the monodromy group of the hypergeometric differential equation for a general choice of parameter.

We now briefly look at the possible special cases.

Still assuming $\gamma-\alpha-\beta\neq0,-1$, if $\gamma=0$, then the canonical form of the residue matrix at $0$ becomes
$$J_0 = \begin{pmatrix}0&1\\0&0\end{pmatrix}.$$
Similarly, if $\alpha = \beta$, then
$$J_\infty = \begin{pmatrix}0&1\\0&0\end{pmatrix}.$$
In each case, the corresponding local monodromy matrices are
$$
L_0 = \begin{pmatrix}
1&2\pi i\\
0&1
\end{pmatrix}$$
or
$$
L_\infty = e^{2\pi i \alpha}\begin{pmatrix}
1&2\pi i\\
0&1
\end{pmatrix}.$$
Assuming either or both of these possibilities and carrying out the calculation of the monodromy group as before yields exactly the same matrices as before (with $\gamma = 0$ or $\alpha = \beta$ substituted). Hence, the monodromies calculated above cover the case of $\gamma -\alpha-\beta\neq0,-1$.

\bfseries Case\rm \ \ $\gamma-\alpha-\beta=0.$

Since $\gamma = \alpha+\beta\in(0,2)$, the case $\gamma=0$ does not apply here and we only need the extra assumption $\alpha\neq\beta$ in order to calculate the local monodromies as
\begin{eqnarray*}
L_0 & = & e^{-\pi i (\alpha+\beta)}\begin{pmatrix}
e^{\pi i (\alpha+\beta)}&0\\
0&e^{-\pi i (\alpha+\beta)}
\end{pmatrix}\\
L_1 & = & \begin{pmatrix}
1& 2 \pi i\\
0&1
\end{pmatrix}\\
L_\infty & = & e^{\pi i (\alpha+\beta)}\begin{pmatrix}
e^{\pi i (\alpha-\beta)}&0\\
0&e^{\pi i (\beta-\alpha)}
\end{pmatrix}.
\end{eqnarray*}
The global monodromy group is
\begin{eqnarray*}
M_0 & = & e^{-\pi i (\alpha+\beta)}\begin{pmatrix}
0&-1\\
1&2\cos((\alpha+\beta)\pi)
\end{pmatrix}\\
M_1 & = & \begin{pmatrix}
1&4\sin(\alpha\pi)\sin(\beta\pi)\\
0&1
\end{pmatrix}\\
M_\infty & = & e^{\pi i (\alpha+\beta)}\begin{pmatrix}
2\cos((\alpha-\beta)\pi)&1\\
-1&0
\end{pmatrix}
\end{eqnarray*}
and, as before, it can be shown that this is still valid without the assumption $\alpha\neq\beta$.

\bfseries Case\rm \ \ $\gamma-\alpha-\beta=-1.$

By adding 1 to $\gamma$, we are able to reduce this to the previous case and we are done. 
\begin{remark}
When the monodromy group is reducible, Levelt's theorem cannot be applied to find the global monodromy group. However, since the solutions of a hypergeometric ODE vary smoothly with respect to $\alpha$, $\beta$ and $\gamma$, in the exceptional cases, the monodromy can be calculated as the limit of irreducible groups. The expressions derived above for the monodromy group are still valid when $\alpha$, $\beta$, $\gamma-\alpha$ or $\gamma-\beta$ are integers.
\end{remark}
We are interested in enumerating all the examples of hypergeometric differential equations that could possible arise as the Picard--Fuchs equation of a family of lattice polarised K3 surfaces. For this, we use theorem \ref{squares} that states that the trace of any projective monodromy matrix must be the square root of an integer. Without any exceptions, for the hypergeometric ODE, the trace of the projectivised monodromy matrices are:
\begin{eqnarray*}
\textnormal{trace}(\overline{M}_0) & = & 2\cos(\gamma\pi),\\
\textnormal{trace}(\overline{M}_1) & = & 2\cos((\gamma-\alpha-\beta)\pi),\\
\textnormal{trace}(\overline{M}_\infty) & = & 2\cos((\alpha-\beta)\pi).
\end{eqnarray*}
Since $\textnormal{trace}(A^2) = \textnormal{trace}(A)^2-2\textnormal{det}(A)$ for any $2\times 2$ matrix $A$ and because of the identity $(2\cos(x))^2-2=2\cos(2x)$, we also have
\begin{eqnarray*}
\textnormal{trace}(\overline{M}_0^2) & = & 2\cos(2\gamma\pi) = 2\cos(2(1-\gamma)\pi),\\
\textnormal{trace}(\overline{M}_1^2) & = & 2\cos(2(\gamma-\alpha-\beta)\pi),\\
\textnormal{trace}(\overline{M}_\infty^2) & = & 2\cos(2(\alpha-\beta)\pi).
\end{eqnarray*}
In order to satisfy theorem \ref{squares}, it is necessary that these traces are integers. Since $2\cos(x\pi)\in[-2,2]$ for real values of $x$, these integer traces can only be $-2, -1, 0, 1$ or $2$. Modulo integers, there are only finitely many possible values of $\alpha$, $\beta$ and $\gamma$ that can satisfy theorem \ref{squares}.

The possibilities are determined by
\begin{eqnarray*}
2\cos(2\pi x) = \phantom{-}2 &\ \ \Rightarrow\ \  & x \equiv 0\ \ (\textnormal{mod}\ \ZZ)\\
2\cos(2\pi x) = \phantom{-}1 & \Rightarrow & x \equiv \pm \frac{1}{6}\ \ (\textnormal{mod}\ \ZZ)\\
2\cos(2\pi x) = \phantom{-}0 & \Rightarrow & x \equiv \pm\frac{1}{4}\ \ (\textnormal{mod}\ \ZZ)\\
2\cos(2\pi x) = -1 & \Rightarrow & x \equiv \pm\frac{1}{3}\ \ (\textnormal{mod}\ \ZZ)\\
2\cos(2\pi x) = -2 & \Rightarrow & x \equiv \frac{1}{2}\ \ (\textnormal{mod}\ \ZZ)
\end{eqnarray*}
so that if we let $S = \{ x \in \QQ\ |\ x\equiv0,\frac{1}{6},\frac{1}{4},\frac{1}{3},\frac{1}{2},\frac{2}{3},\frac{3}{4},\frac{5}{6}\ \textnormal{mod}\ \ZZ\}$ then we require $1-\gamma$, $\gamma-\alpha-\beta$ and $\alpha-\beta\in S$.
\begin{remark}
We could expect that although the generators have the correct traces, the products of generators may not. In fact, it can be checked that this problem does not occur.
\end{remark}
It was known classically (see for example \cite{ford}) that independent solutions $f_1,f_2$ to the hypergeometric ODE $_2F_1(\alpha,\beta,\gamma)$ provide a multi--valued map
\begin{eqnarray*}
\CC & \rightarrow & \PP^1\\
z & \mapsto & f_1(z)/f_2(z)
\end{eqnarray*}

that sends the upper--half plane to a curvilinear triangle with angles $\lambda\pi$, $\mu\pi$, $\nu\pi$ at the vertices where $\lambda = |1-\gamma|$, $\mu = |\gamma-\alpha-\beta|$ and $\nu = |\alpha-\beta|$. Furthermore, there exist $\alpha^\prime$, $\beta^\prime$, $\gamma^\prime$ with $\alpha-\alpha^\prime\in\ZZ$, $\beta-\beta^\prime\in\ZZ$ and $\gamma-\gamma^\prime\in\ZZ$ so that $\lambda^\prime$, $\mu^\prime$ and $\nu^\prime$ defined as above satisfy
$$
0\leq\lambda^\prime,\mu^\prime,\nu^\prime < 1
$$
and
$$
\lambda^\prime+\mu^\prime+\nu^\prime<1+2\min(\lambda^\prime,\mu^\prime,\nu^\prime).
$$
This is explained fully in \cite{Beu}. With respect to this choice of $\alpha^\prime$, $\beta^\prime$ and $\gamma^\prime$, the projective monodromy group is Fuchsian if and only if
$$
\lambda^\prime+\mu^\prime+\nu^\prime < 1.
$$
In particular, if $\lambda^\prime$, $\mu^\prime$ or $\nu^\prime = 0$ then the monodromy group is Fuchsian.
\begin{theorem}
The projective monodromy group, $\Gamma$, of $_2F_1(\alpha,\beta,\gamma)$ is a Fuchsian group with $\textnormal{trace}(A)^2 \in \ZZ$ for all $A\in\Gamma$ if and only if
there exists $\alpha^\prime\equiv\alpha\ \textnormal{mod}\,\ZZ$, $\beta^\prime\equiv\beta\ \textnormal{mod}\,\ZZ$, $\gamma^\prime\equiv\gamma\ \textnormal{mod}\,\ZZ$ such that $\lambda^\prime = |1-\gamma^\prime|$, $\mu^\prime = |\gamma^\prime-\alpha^\prime-\beta^\prime|$ and $\nu^\prime = |\alpha^\prime-\beta^\prime|$ satisfy
$$
\lambda^\prime, \mu^\prime, \nu^\prime \in \left\{ 0,\frac{1}{6},\frac{1}{4},\frac{1}{3},\frac{1}{2},\frac{2}{3},\frac{3}{4},\frac{5}{6}\right\},
$$
and
$$
\lambda^\prime + \mu^\prime + \nu^\prime < 1.
$$
\end{theorem}
\begin{example}
For example, the projective monodromy group of $_2F_1(\frac{1}{6},\frac{5}{6},\frac{1}{2})$ is generated by
$$
\overline{M}_0 = \begin{pmatrix}
0&\omega^2\\
-\omega&0
\end{pmatrix},\ \ \overline{M}_1 = \begin{pmatrix}
0&-1\\
1&0
\end{pmatrix},\ \ \overline{M}_\infty = \begin{pmatrix}
\omega &0\\
0&\omega^2
\end{pmatrix}
$$
where $\omega=e^{2\pi i /3}$. This is a finite group, isomorphic to the dihedral group $D_{12}$, and is not conjugate to any finite subgroup of $\PP\textnormal{Sl}(2,\RR)$ because $(\lambda^\prime,\mu^\prime,\nu^\prime) = (\frac{1}{2}, \frac{1}{2}, \frac{2}{3})$ and $1\leq\frac{1}{2}+\frac{1}{2}+\frac{2}{3}<1+2\min(\frac{1}{2}, \frac{1}{2}, \frac{2}{3})$.
\end{example}
\begin{example}\label{finiteLocalMono}
The projective monodromy group of $_2F_1(\frac{1}{3},\frac{1}{12},\frac{2}{3})$ is generated by
$$
\overline{M}_1 = \frac{1}{\sqrt{2}}\begin{pmatrix}
1&1\\
-1&1
\end{pmatrix},\ \ \overline{M}_\infty = \frac{1}{\sqrt{2}}\begin{pmatrix}
1&2+\sqrt{3}\\
-2+\sqrt{3}&1
\end{pmatrix}.
$$
In this example, $\overline{M}_0^3=(\overline{M}_1.\overline{M}_\infty)^{-3}=\overline{M}_1^4=\overline{M}_\infty^4=id$. Despite having generators of finite order, this group is infinite since, for example, $\overline{M_1}.\overline{M}_\infty^{-1}$ has infinite order.
\end{example}
\begin{example}
The hypergeometric differential equation $_2F_1(\frac{1}{12},\frac{1}{12},\frac{2}{3})$ has monodromy group
\begin{eqnarray*}
M_0 & = & e^{\frac{\pi i}{3}}\begin{pmatrix}-1&-1\\1&0\end{pmatrix},\\
M_1 & = & i\begin{pmatrix}0&1\\-1&0\end{pmatrix},\\
M_\infty & = & e^{\frac{\pi i}{6}}\begin{pmatrix}1&1\\0&1\end{pmatrix}
\end{eqnarray*}
so that the projective monodromy group for $_2F_1(\frac{1}{12},\frac{1}{12},\frac{2}{3})$, being generated by $\begin{pmatrix}0&1\\-1&0\end{pmatrix}$ and $\begin{pmatrix}1&1\\0&1\end{pmatrix}$, is equal to $\textnormal{PSl}(2,\ZZ)$.
\end{example}
\begin{remark}
The finite monodromy groups of hypergeometric differential equations, as classified by Schwarz, are all finite groups of rotations of a sphere. Therefore the \emph{Fuchsian} hypergeometric monodromy groups are all infinite. Example \ref{finiteLocalMono} demonstrates that the Fuchsian groups satisfying theorem \ref{squares} include groups whose local monodromies are all of finite order. It would be interesting to know if there are any families of K3 surfaces with hypergeometric monodromy group, but without any maximally unipotent local monodromies.
\end{remark}
Using the rigidity of the local system of solutions, we are able to find the monodromy group of a hypergeometric differential equation with very little trouble. As we see next, this method can be adapted to be used in the case where the local system of solutions is not rigid.
\subsection{Non-Rigid Systems}
We have found a number of Picard--Fuchs differential equations with more than three regular singular points. The local system of solutions to the symmetric square--root of these differential equations is not rigid. To find the monodromy representation, we need to analytically continue a fixed basis of solutions around each singular point in turn. Without a closed form for the solutions of the differential equation, we are forced to perform this analytic continuation numerically. This approach has the problem that the matrices obtained will be filled with inexact decimal expansions, and the monodromy group will not be rigorously determined. In practice, however, we are able to determine the projective monodromy group with certainty by taking advantage of a special property of the differential equations found in the previous chapter.

We will later see that the non--hypergeometric differential equations that we have found are all generalised Lam\'e differential equations. This means that all but one of the projective local monodromies will have order 2 and be conjugate to the pseudo reflection  $\begin{pmatrix} 1&0\\0&-1\end{pmatrix}$.

Notice that
$$
\begin{pmatrix} 1&i\\i&1\end{pmatrix}^{-1}.\begin{pmatrix} 1&0\\0&-1\end{pmatrix}.\begin{pmatrix} 1&i\\i&1\end{pmatrix} = i.\begin{pmatrix} 0&1\\-1&0\end{pmatrix}
$$
so that, in the projective monodromy group, the transformations $\begin{pmatrix} 1&0\\0&-1\end{pmatrix}$ and $\begin{pmatrix} 0&1\\-1&0\end{pmatrix}$ are conjugate.

Although the monodromy group will not be rigid, to uniquely determine the monodromy group from the local monodromies we will see that, for a generalised Lam\'e equation, we only need to additionally know the values of $t_{i,j} := \textnormal{trace}((M_i.M_j)^2)$ for pairs of generators of the group. This missing data can only be determined numerically, but since theorem \ref{squares} tells us that $t_{i,j}\in\ZZ$, we know that we are approximating integers. We do not require the approximation to be very accurate before we are certain that we have determined the correct values.

\begin{theorem}\label{squaretraces}
Let $M\in\textnormal{Gl}(2,\CC)$ and write $N= \begin{pmatrix}1&0\\0&-1\end{pmatrix}$. If
$$
\textnormal{trace}((N.M)^2)\neq\textnormal{trace}(M^2)
$$
then the group generated by $M$ and $N$ is rigid.
\end{theorem}
\begin{proof}
Writing $\alpha_1, \alpha_2$, resp. $\beta_1, \beta_2$ for the eigenvalues of $N.M$, resp. $M$, we note that $\alpha_1\alpha_2=-\beta_1\beta_2$ and
\begin{eqnarray*}
\textnormal{trace}((N.M)^2) & = & \alpha_1^2+\alpha_2^2\\
\textnormal{trace}(M^2) & = & \beta_1^2+\beta_2^2
\end{eqnarray*}
From this, it is a simple exercise to show that $\textnormal{trace}((N.M)^2)\neq\textnormal{trace}(M^2)$ if and only if $\{\alpha_1,\alpha_2\}\cap\{\beta_1,\beta_2\}=\emptyset$. The proposition follows from Levelt's theorem (theorem \ref{Levelt}), applied in the case $A = M$ and $B = N.M$.
\end{proof}
Our generalised Lam\'e equations are second order Fuchsian differential equations with $k+1$ regular singular points at $\alpha_0,\alpha_1,\ldots,\alpha_k$. At $k$ of the singular points $\alpha_1,\ldots,\alpha_k$, the projective local monodromy matrix is conjugate to $N = \begin{pmatrix}1&0\\0&-1\end{pmatrix}$ and at $\alpha_0$, the monodromy is arbitrary. Under this assumption, every pair of projective local monodromies will include one conjugate to $N$. Subject to the conditions on the traces, each pair of monodromies generate a rigid group. From this, we deduce that the monodromy group is determined up to simultaneous conjugation by specifying suitable values for $\textnormal{trace}((M_{\alpha_i}M_{\alpha_j})^2)$.

In particular, in any given example, we need to find the values
$$
t_{i,j} := \textnormal{trace}((M_{\alpha_i}M_{\alpha_j})^2)\ \ \ \ \ \textnormal{for}\ \ \ i=0,\ldots,k\ \ \ \ \textnormal{and}\ \ \ i\neq j
$$
and we require
\begin{eqnarray*}
t_{i,0} & \neq & \textnormal{trace}(M_{\alpha_0}^2)\\
t_{i,j} & \neq & 2\ \ \ \ \textnormal{for}\ \ j\neq0.
\end{eqnarray*}
If we then find specific matrices $M_{\alpha_0}, \ldots,M_{\alpha_k}$ conjugate to the local monodromies at $\alpha_0,\ldots,\alpha_k$ satisfying $M_{\alpha_0}\ldots{M_{\alpha_k}}=id$ and with $\textnormal{trace}((M_{\alpha_i}M_{\alpha_j})^2) = t_{i,j}$, then we know we have found the monodromy representation of our ODE.
\subsection{Numerical Determination of the Monodromy Group}
We shall numerically and non-rigorously determine the monodromy group of a Picard-Fuchs differential equation in the case of more than three singular points. We know that to uniquely specify the projective monodromy group of a generalised Lam\'e equation, we need to find the values $\textnormal{trace}((M_i.M_j)^2)$ for pairs of generators $M_i$ and $M_j$. In this section, we provide an algorithm to find these numbers in any given example.

For convenience we choose to deal with a Fuchsian differential system with the same singular points and monodromy representation as the initial differential equation (see Theorem \ref{equivSys}). This is not essential, but the convenience of dealing with a Fuchsian differential system should not be underestimated. For example, the local monodromies of a differential system are readily obtained by using corollary \ref{disteigen}. Also, the algorithm we obtain in this section is far simpler to explain when using differential systems.

We take as our starting point a Fuchsian system
\begin{equation}\label{system2}
\frac{d}{dz}\omega = \sum_{i=0}^k\frac{R_i}{(z-\alpha_i)}\omega.
\end{equation}
with singular points $\alpha_0,\ldots,\alpha_k$. We assume that $k+1>3$, the point $\infty$ is not a singular point of the system and that the eigenvalues of $R_i$ differ by $\frac{1}{2}$ for $ i = 1,\ldots,k$ so that the system is of generalised Lam\'e type.

At the time of writing, there is a convenient java applet, \cite{Java}, available on the internet that can numerically approximate the monodromy representation of such a system. The applet is limited in that it does not produce an answer if the singular points are too close together and it cannot be used on a few of our examples. Instead, we shall use Maple to find the trace of a monodromy transformation of a Fuchsian system along a closed path. The final Maple procedure is given in appendix \ref{maple}.

\begin{figure}
	\setlength{\unitlength}{2cm}
	\begin{picture}(6,3)(-1,-1.5)
  	\put(0,0){\circle*{0.05}}
  	\put(1,0){\circle*{0.05}}
 	 	\put(2,0){\circle*{0.05}}
  	\put(3,0){\circle*{0.05}}
  	\put(4,0){\circle*{0.05}}
 
  	\put(2.5,1){\vector(-1,0){2}}
  	\put(0.5,1){\vector(0,-1){2}}
  	\put(0.5,-1){\vector(1,0){1}}
  	\put(1.5,-1){\vector(0,1){2}}
  	\put(1.5,1){\vector(1,0){1}}
  	\put(1,-1){\makebox(0,-0.25){$\gamma_i$}}
 
  	\put(2.5,1){\vector(1,0){1}}
  	\put(3.5,1){\vector(0,-1){2}}
  	\put(3.5,-1){\vector(1,0){1}}
  	\put(4.5,-1){\vector(0,1){2}}
  	\put(4.5,1){\vector(-1,0){2}}
  	\put(4,-1){\makebox(0,-0.25){$\gamma_j$}}
 		\put(1,0){\makebox(0.2,0.2){$\alpha_i$}}
 		\put(4,0){\makebox(0.2,0.2){$\alpha_j$}}
 
		\put(2.5,1){\makebox(0.2,0.2){$p_0$}}
  	\put(2.5,1){\circle*{0.05}}
	\end{picture}
	\caption{Monodromy Paths from a Base Point}
	\label{fig:MonodromyPathsFromABasePoint}
\end{figure}
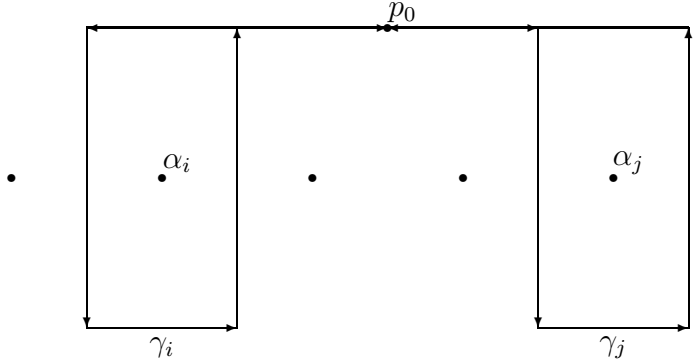

If $M_0,\ldots,M_k$ are the monodromy matrices to be determined, then our aim is to calculate the integers $\textnormal{trace}(M_i.M_j)^2$. We will then find matrices conjugate to the local monodromies and whose products have the correct traces. Subject to the conditions of theorem~\ref{squaretraces}, these matrices will generate the monodromy group. Bearing this in mind, We set up the monodromy paths corresponding to generators of the monodromy representation. We pick a base point, $p_0$, and choose one anti--clockwise loop from $p_0$ enclosing each singular point of the system. For simplicity, and without loss of generality, we shall take these loops to be piecewise linear as shown in figure~\ref{fig:MonodromyPathsFromABasePoint}. We define such a closed path in Maple by specifying the list of coordinates of its vertices, starting with $p_0$.

To calculate $\textnormal{trace}((M_i.M_j)^2)$, we fix the two singular points $\alpha_i$ and $\alpha_j$ and consider the composite path $\gamma_i\gamma_j$ around these points.

\begin{lemma}\label{ident}
There exists a fundamental matrix, $\Omega$, for the arbitrary Fuchsian system (\ref{system2}) with $\Omega(p_0) = I_n$, the identity matrix.
\end{lemma}
\begin{proof}
Let $\Omega_0$ be any fundamental matrix for (\ref{system2}). Writing $M = \Omega_0(p_0)$, we know that $\Omega = \Omega_0.M^{-1}$ is also a fundamental matrix. As required, we have $\Omega(p_0) = I_n$.
\end{proof}

If $\Omega$ is a fundamental matrix of solutions satisfying lemma~\ref{ident}, then the analytic continuation of $\Omega$ around the path $\gamma_i\gamma_j$ will be the fundamental matrix $\Omega.M_j.M_i$ satisfying $(\Omega.M_j.M_i)(p_0) = M_j.M_i$. This calculation can now be stated as an initial value problem that can be solved using Maple's \verb|dsolve()| command.

For Maple, we are required to convert our second order differential system
$$
\omega^\prime = A.\omega
$$
into a four linear differential equations
\begin{eqnarray*}
\omega_{11}^\prime & = & a_{11}\omega_{11} + a_{12}\omega_{21}\\
\omega_{21}^\prime & = & a_{21}\omega_{11} + a_{22}\omega_{21}\\
\omega_{12}^\prime & = & a_{11}\omega_{12} + a_{12}\omega_{22}\\
\omega_{22}^\prime & = & a_{21}\omega_{12} + a_{22}\omega_{22}
\end{eqnarray*}
where
$$
A = \begin{pmatrix}
a_{11}&a_{12}\\
a_{21}&a_{22}
\end{pmatrix}
$$
and $(\omega_{11},\omega_{21})^t$ and $(\omega_{21},\omega_{22})^t$ are linearly independent solutions of the system, or rather, columns of a fundamental matrix.

Maple's \verb|dsolve()| command can be used to numerically solve differential equations with initial conditions. However, it is limited to taking initial values at real points and finding the solutions at other real values by solving along a real interval. Figure~\ref{fig:MapleSDsolveCommand} gives a typical example of this use of \verb|dsolve|.
\begin{figure}[h]
	\begin{verbatim}
    eqn  := { diff(f(z),z) - z*f(z) };
    init := { f(0) = 7 }
    F := dsolve( eqn union init, f(z), type=numeric );
                     
    F(1);
              [z = 1., f(z) = 11.5411547475995189]
\end{verbatim}
	
	\caption{Maple's dsolve Command}
	\label{fig:MapleSDsolveCommand}
\end{figure}

Our closed path $\gamma_i\gamma_j$ splits into a composition of linear paths $l_1,\ldots,l_m\subset\CC$. For each $k\in\{1,\ldots,m\}$, let $T_k$ be the M\"obius transformation fixing $\infty$ and mapping the endpoints of $l_k$ to 0 and 1. These transformations send the linear sections of our path to a real interval that can be handled by Maple.

The algorithm of appendix~\ref{maple} breaks down into the following steps.
\begin{enumerate}
	\item Use $T_1$ to map $l_1$ to $(0, 1)$.
	\item With the initial conditions $\Omega(0) = I_2$, use \verb|dsolve()| to find $S:=\Omega(1)$.
	\item Repeat steps (A) and (B) for each of the next linear paths along $\gamma_i\gamma_j$, this time using  the initial conditions $\Omega(0) = S$.
\end{enumerate}
When all the path segments have been traversed, the final solution, $S$, is a numerical approximation of the monodromy matrix corresponding to the closed path. Although the matrix depends on the initial choice of basis of solutions, the trace of the matrix does not. We square this trace, record the integer that this number approximates and discard the matrix.

The numerical accuracy of Maple's \verb|dsolve| command is determined by the values of the variables \emph{abserr} and \emph{relerr}. In the algorithm, these values determine acceptable errors between steps, but do not rigorously guarantee any particular accuracy in the final answer. In practice, the default accuracy used in appendix \ref{maple} approximates the integral traces to within 4 decimal places.

\section{Arithmetic Fuchsian Groups}
The monodromy groups of the symmetric square roots of our Picard--Fuchs differential equations are Fuchsian groups defined up to conjugacy. In order to write down the monodromy group, we must choose conjugacy representatives of the generators. This seemingly arbitrary choice is unfortunate and we want to find a way of either choosing good representatives or, better still, of specifying the conjugacy class in a matrix free way. To do this, we take a detailed look at arithmetic Fuchsian groups.

By Fuchsian group, we mean a discrete subgroup of $\textnormal{PSl}(2,\RR)$. Our treatment of arithmetic Fuchsian groups follows that of \cite{AB} and \cite{Ka}. We are interested in Fuchsian groups that are derived from some quaternion algebra. To introduce these, we make some standard definitions.
\begin{definition}
Let $\KK$ be a field of characteristic $\neq{2}$. A \emph{quaternion algebra} over $\KK$ is a central simple algebra of rank 4 over $\KK$.

Equivalently, and more concretely, a quaternion algebra is any algebra isomorphic to one with basis $\{1,i,j,k\}$ and relations
$$
i^2=a,\ \ \ \ \ \ j^2=b,\ \ \ \ \ \ k=ij=-ji
$$
with $a,b\in\KK^*$. Such an algebra is denoted by the Hilbert symbol $\left(\frac{a,b}{\KK}\right)$.

Although the Hilbert symbol provides a convenient way to denote a quaternion algebra, it must be noted that the symbol is not uniquely defined for a given quaternion algebra. For example,

$$
\left(\frac{3,-5}{\QQ}\right)\cong \left(\frac{2,5}{\QQ}\right)
$$
via the isomorphism
\begin{eqnarray*}
i & \mapsto & I = \frac{4}{3}i+j+\frac{k}{3}\\
j & \mapsto & J = -\frac{5}{3}i-j+\frac{k}{3}.
\end{eqnarray*}
Wherever this ambiguity could cause problems, we shall use the normal form for quaternion algebras over $\QQ$ from \cite{AB}.
\end{definition}
\begin{definition}
For $x=x_0+x_1i+x_2j+x_3k$, we define the conjugate
$$\bar{x}:=x_0-x_1i-x_2j-x_3k,$$
the reduced norm
$$
\textnormal{N}(x) := x\bar{x}=x_0^2-ax_1^2-bx_2^2+abx_3^2
$$
and reduced trace
$$
\textnormal{Tr}(x) := x+\bar{x}=2x_0.
$$
These correspond to the determinant and trace of $x$ under any 2 dimensional representation of the algebra.
\end{definition}
\begin{definition}
By an \emph{order}, $\Oh$, in a (quaternion) algebra $A$ over $\KK$, we mean a finitely generated $\Oh_\KK$--module generating the algebra $A$ over $\KK$ which is also a subring of $A$ containing 1.

In particular, an order of a quaternion algebra $A$ over $\QQ$ is a free $\ZZ$--module of rank 4 that is a subring of $A$ containing 1.
\end{definition}
\begin{theorem}[\cite{Ka}, p 119]\label{afg}
Let $A$ be a quaternion algebra over a number field $\KK$ with $n=[\KK\colon\QQ]$. Let $\textnormal{Gal}(\KK/\QQ)=\{\sigma_i=\textnormal{id},\sigma_2,\ldots,\sigma_n\}$.

Suppose $\KK$ is totally real (ie. $\sigma_i(\KK)\subset\RR$ for $i=1,\ldots,n$) and there exist $\RR$--isomorphisms
\begin{align*}
\rho_1  \colon {A}\otimes_\KK\RR & \rightarrow \textnormal{M}_{2}(\RR)\\
\rho_i  \colon {A^{\sigma_i}}\otimes_\KK\RR & \rightarrow \HH=\left(\frac{-1,-1}{\RR}\right) \hspace{40pt}(2\leq{i}\leq{n}).
\end{align*}
If $\Oh^1$ is the group of elements of reduced norm 1 in any order $\Oh$ of $A$, then
$$
\Gamma(A,\Oh):=\rho_1(\Oh^1)/\{\pm I_2\}
$$
is a discrete subgroup of $\textnormal{PSl}(2,\RR)$.
\qed\end{theorem}
In particular, for a quaternion algebra over $\QQ$ to satisfy the conditions of theorem \ref{afg}, we only require that $A\otimes\RR\cong\textnormal{M}_2(\RR)$. In this case, we say that $A$ is indefinite.

Writing $A=\left(\frac{a,b}{\QQ}\right)$ with $a\geq{b}$, we may define a standard representation
$$
\rho:\left(\frac{a,b}{\QQ}\right)\rightarrow\textnormal{M}_2(\QQ(\sqrt{a}))
$$
by
\begin{align*}
\rho(1) =&\ I_2\\
\rho(i)=&
\begin{pmatrix}
\sqrt{a} & 0 \\
0 & -\sqrt{a} \\
\end{pmatrix}\\
\rho(j)=&
\begin{pmatrix}
0 & 1 \\
b & 0 \\
\end{pmatrix}\\
\rho(k)=&\rho(i)\rho(j) = 
\begin{pmatrix}
0 & \sqrt{a} \\
-b\sqrt{a} & 0 \\
\end{pmatrix}.
\end{align*}
If $a>0$, it is clear that
$$
A\otimes_\QQ\RR\overset{\rho\otimes\RR}{\cong}\textnormal{M}_2(\RR).
$$
If $a>0$ and $b>0$, however, we may also define the slightly less standard representation
$$
\tau:\left(\frac{a,b}{\QQ}\right)\rightarrow\textnormal{M}_2(\QQ(\sqrt{a},\sqrt{b}))\subset\textnormal{M}_2(\RR)
$$
by
\begin{align*}
\tau(1) =&\ I_2\\
\tau(i)=&
\begin{pmatrix}
\sqrt{a} & 0 \\
0 & -\sqrt{a} \\
\end{pmatrix}\\
\tau(j)=&
\begin{pmatrix}
0 & \sqrt{b} \\
\sqrt{b} & 0 \\
\end{pmatrix}\\
\tau(k)=&
\begin{pmatrix}
0 & \sqrt{ab} \\
-\sqrt{ab} & 0 \\
\end{pmatrix}.
\end{align*}

\begin{definition}
We say that a Fuchsian group is \emph{arithmetic} if it is commensurable with some $\Gamma(A,\Oh)$ as defined in theorem \ref{afg}.

We say that a Fuchsian group $\Gamma$ is \emph{derived} from an order $\Oh$ if $\Gamma$ is a finite index subgroup of $\Gamma(A,\Oh)$
\end{definition}

\begin{theorem}
A Fuchsian group, $\Gamma$, is arithmetic if and only if $\Gamma^{(2)}$, the group generated by squares of elements in $\Gamma$, is a finite subgroup of some $\Gamma(A,\Oh)$.

Also, if $G$ is a Fuchsian group with $\textnormal{Area}(\mathcal{H}/G)<\infty$, then $G$ is a finite subgroup of some $\Gamma(A,\Oh)$ where $A$ is a quaternion algebra over a number field $\mathbb{K}$ if and only if $\textnormal{Tr}(g)\in\Oh_\mathbb{K}$ for all $g \in G$ where $\Oh_\mathbb{K}$ denotes the ring of integers of $\mathbb{K}$.
\end{theorem}

By virtue of theorem \ref{squares}, we immediately get
\begin{corollary}
If $\Gamma_S$ is the monodromy group of a rank 19 lattice polarised family of K3 surfaces, then $\Gamma_S^{(2)}$ is derived from an order in a quaternion algebra over the rational numbers.
\end{corollary}

In any example, it is possible to find the quaternion algebra $A$ and the order $\Oh$ so that $\Gamma^{(2)}_S$ is a finite index subgroup of $\Gamma(A,\Oh)$. We take a look at this in the next section.

\section{Monodromy of Picard--Fuchs Equations}

Arithmetic Fuchsian groups, $\Gamma$, with $\textnormal{Area}(\mathcal{H}/\Gamma)<\infty$ fall into two distinct types: those with cusps, and those without. The Fuchsian groups without cusps have compact fundamental domains.

Similarly, either a quaternion algebra has zero divisors, or it is a division algebra.

It is well-known that these two dichotomies coincide. A cocompact arithmetic Fuchsian group is commensurable with the elements of unit norm in an order of a quaternion algebra that is a division algebra. And conversely.

We take a look at an example calculation of the monodromy group of a non-hypergeometric Picard--Fuchs differential equation. We look at family XI. This is defined in $\PP^3$ by
$$
x_0^4+x_1^4+x_2^4+x_3^4 + \lambda(x_0^2+x_1^2+x_2^3+x_4^2)^2 = 0
$$
and is invariant under an action of the group labeled $2^4D_6$. This family degenerates at the roots of $(\lambda+1)(\lambda+\frac{1}{2})(\lambda+\frac{1}{3})(\lambda+\frac{1}{4})$ and at $\infty$.

The symmetric square root of the Picard--Fuchs differential equation is given in chapter 2. We are interested in a differential system with the same singular points and monodromy as this ODE. For this, we apply the algorithm of theorem \ref{equivSys} and obtain the following Fuchsian system:
$$
\frac{d}{d\lambda}\omega = \left(\frac{R_{-1}}{\lambda+1}+ \frac{R_{-\frac{1}{2}}}{\lambda+\frac{1}{2}} + \frac{R_{-\frac{1}{3}}}{\lambda+\frac{1}{3}} + \frac{R_{-\frac{1}{4}}}{\lambda+\frac{1}{4}} \right)\omega
$$
with residue matrices
\begin{eqnarray*}
R_{-1} & = & \begin{pmatrix}
-\frac{2347}{5760} & \frac{1}{6} \\
-\frac{12267769}{5529600} & \frac{5227}{5760}
\end{pmatrix},\\
R_{-\frac{1}{2}} & = & \begin{pmatrix}
\frac{2947}{960} & -2\\
\frac{7270249}{1843200} & -\frac{2467}{960}
\end{pmatrix},\\
R_{-\frac{1}{3}} & = & \begin{pmatrix}
-\frac{5067}{640} & \frac{9}{2}\\
-\frac{3032881}{204800} & \frac{5387}{640}
\end{pmatrix},\\
R_{-\frac{1}{4}} & = & \begin{pmatrix}
\frac{6487}{1440} & -\frac{8}{3}\\
\frac{37410529}{5529600} & -\frac{5767}{1440}
\end{pmatrix}
\end{eqnarray*}
at the finite singular points and residue matrix
\begin{eqnarray*}
R_\infty & = & -(R_{-1} + R_{-\frac{1}{2}} + R_{-\frac{1}{3}} + R_{-\frac{1}{4}} )\\
& = & \begin{pmatrix}
\frac{3}{4} & 0\\
\frac{291119}{46080} & -\frac{11}{4}
\end{pmatrix}
\end{eqnarray*}
at the singular point at $\infty$. This differential system is not uniquely defined by the property of having the same monodromy representation as the initial differential equation. Looking at the residue matrices, our choice is clearly not minimal either. Nevertheless, it does the job. The four residue matrices at the finite singular points each have eigenvalues $0$ and $\frac{1}{2}$ and the matrix $R_\infty$ has eigenvalues $\frac{3}{4}$ and $-\frac{11}{4}$. Hence this differential system is of generalised Lam\'e type.

We write the eigenvalues of the residue matrices at the singular points in the so-called Riemann symbol:
$$
\left\{\newcommand\T{\rule{0pt}{2.3ex}}
\newcommand\B{\rule[-1.0ex]{0pt}{0pt}}
\centering
\begin{tabular}{c c c c c}
	\T\B $-1$ & $-\frac{1}{2}$ & $-\frac{1}{3}$ & $-\frac{1}{4}$ & $\infty$ \\
	\T\B $0$  & $0$            & $0$            & $0$            & $-\frac{11}{4}$\\
	\T\B $\frac{1}{2}$ & $\frac{1}{2}$ & $\frac{1}{2}$ & $\frac{1}{2}$ & $\frac{3}{4}$
\end{tabular}\right\}
$$
The first row states the singular points and the subsequent rows show the eigenvalues at that point.

Choosing a monodromy path encircling each finite singular point in an anti-clockwise direction from a base point, the monodromy group will be generated by elements $M_{-1}$, $M_{-\frac{1}{2}}$, $M_{-\frac{1}{3}}$, $M_{-\frac{1}{4}}$ each individually conjugate to 
$$\exp(2 \pi i \begin{pmatrix}0&0\\0&\frac{1}{2}\end{pmatrix}) = \begin{pmatrix}1&0\\0&-1\end{pmatrix}.$$
The monodromy about $\infty$ will be given by $$M_\infty = (M_{-1}M_{-\frac{1}{2}}M_{-\frac{1}{3}}M_{-\frac{1}{4}})^{-1}$$ and will itself be conjugate to 
$$\exp(2 \pi i \begin{pmatrix}\frac{3}{4}&0\\0&-\frac{11}{4}\end{pmatrix}) = -i\begin{pmatrix}1&0\\0&-1\end{pmatrix}.$$
Letting $t_{i,j} = \textnormal{trace}((M_i.M_j)^2)$ for $i,j\in\{\infty,-1,\frac{1}{2},\frac{1}{3},\frac{1}{4}\}$, the algorithm of appendix \ref{maple} finds 
$$\begin{matrix}
t_{\infty,-1} = -6 &\ \ \ t_{\infty,-\frac{1}{2}} = -22 &\ t_{\infty,-\frac{1}{3}} = -30 & \ t_{\infty,-\frac{1}{4}} = -14\\
&\ t_{-1,-\frac{1}{2}} = 10 & t_{-1,-\frac{1}{3}} = 34 & t_{-1,-\frac{1}{4}} =  30\\
& & t_{-\frac{1}{2},-\frac{1}{3}} = 10 & t_{-\frac{1}{2},-\frac{1}{4}} = 22\\
& & & t_{-\frac{1}{3},-\frac{1}{4}} = 6
\end{matrix}$$
None of these traces are equal to 2 or $-2$ as required by theorem \ref{squaretraces}, and so the monodromy group is uniquely determined by these traces. From this point, it is an easy exercise to find generating matrices conjugate to the local monodromies and with the correct traces. The problem can be turned into a set of simultaneous polynomial equations and solved using an algebra package such as Maple or Mathematica. In this example, the matrices
\begin{eqnarray*}
	M_\infty & = & \frac{1}{\sqrt{2}}\begin{pmatrix}0&-1+\sqrt{3}\\-1-\sqrt{3}&0\end{pmatrix}\\
	M_{-1} & = & -\sqrt{-1}\begin{pmatrix}1&-1+\sqrt{3}\\-1-\sqrt{3}&-1\end{pmatrix}\\
	M_{-\frac{1}{2}} & = & -\sqrt{-1}\begin{pmatrix}\sqrt{3}&2\\-2&-\sqrt{3}\end{pmatrix}\\
	M_{-\frac{1}{3}} & = & -\sqrt{-1}\begin{pmatrix}1&1+\sqrt{3}\\1-\sqrt{3}&-1\end{pmatrix}\\
	M_{-\frac{1}{4}} & = & -\frac{\sqrt{-1}}{\sqrt{2}}\begin{pmatrix}0&1+\sqrt{3}\\1-\sqrt{3}&0\end{pmatrix}\\
\end{eqnarray*}
can be checked to have the correct properties. Being Fuchsian, the projective monodromy group is a discrete subgroup of $\PP\textnormal{Sl}(2,\RR)$ and in this case is generated by
\begin{eqnarray*}
\overline{M_\infty} & = & \frac{1}{\sqrt{6}}(3j-k)\\
\overline{M_{-1}} & = & \frac{1}{\sqrt{3}}(i+3j-k)\\
\overline{M_{-\frac{1}{2}}} & = & i+2j\\
\overline{M_{-\frac{1}{3}}} & = & \frac{1}{\sqrt{3}}(i+3j+k)\\
\overline{M_{-\frac{1}{4}}} & = & \frac{1}{\sqrt{6}}(3j+k)
\end{eqnarray*}
where
$$
i = \begin{pmatrix}
\sqrt{3}&0\\
0&-\sqrt{3}
\end{pmatrix},
$$
$$j = \begin{pmatrix}
0&1\\
-1&0
\end{pmatrix}
$$
and $k=i.j$. This projective monodromy group is contained in the image of a representation of the quaternion algebra $\left(\frac{3,-1}{\QQ[\sqrt{2},\sqrt{3}]}\right)$.

We now list the monodromy groups for each of our example families. For each example, we recall the differential equation and show its Riemann symbol detailing the eigenvalues at the singular points. We then show the monodromy group and write down the projective monodromy group in terms of elements in a quaternion algebra.
\begin{Eqn}{I, II, III}
$_2F_1( \frac{1}{12}, \frac{1}{12}, \frac{2}{3})$

The monodromy group is generated by the transformations
\begin{eqnarray*}
	M_0 & = & \frac{1+\sqrt{-3}}{2}\begin{pmatrix}1&1\\-1&0\end{pmatrix}\\
	M_1 & = & \sqrt{-1}\begin{pmatrix}0&1\\-1&0\end{pmatrix}\\
	M_\infty & = & \frac{\sqrt{3}+\sqrt{-1}}{2}\begin{pmatrix}1&1\\0&1\end{pmatrix}
\end{eqnarray*}
satisfying $M_0.M_1.M_\infty = id$.

Thus, the projective monodromy group is generated by $\overline{M}_1 = \begin{pmatrix}0&1\\-1&0\end{pmatrix}$ and $\overline{M}_\infty = \begin{pmatrix}1&1\\0&1\end{pmatrix}$ and is equal to $\textnormal{PSl}(2,\ZZ)$.
\end{Eqn}
 
\begin{Eqn}{IV}
$$
p(\lambda)\frac{d^2}{d\lambda^2} + \frac{1}{2}p^\prime(\lambda)\frac{d}{d\lambda} + \frac{3}{16}\lambda^2+\frac{3}{2}\lambda+\frac{3}{2}.
$$

with $p(\lambda) = (\lambda+1)(\lambda+2)(\lambda^2+16\lambda+16)$.

$$
\left\{\newcommand\T{\rule{0pt}{2.3ex}}
\newcommand\B{\rule[-1.0ex]{0pt}{0pt}}
\centering
\begin{tabular}{c c c c c}
	\T\B $-1$ & $-\alpha$ & $-2$ & $-\beta$ & $\infty$ \\
	\T\B $0$  & $0$       & $0$  & $0$   &   $-\frac{5}{4}$\\
	\T\B $\frac{1}{2}$ & $\frac{1}{2}$ & $\frac{1}{2}$ & $\frac{1}{2}$ & $\frac{1}{4}$
\end{tabular}\right\}
$$

The monodromy group is generated by the transformations

\begin{eqnarray*}
M_{-1} & = & \frac{\sqrt{-1}}{\sqrt{6}}\begin{pmatrix}3\sqrt{3}&6+\sqrt{3}\\-6+\sqrt{3}&-3\sqrt{3}\end{pmatrix}\\
M_{-\alpha} & = & \sqrt{-1}\begin{pmatrix}\sqrt{3}&2\\-2&-\sqrt{3}\end{pmatrix}\\
M_{-2} & = & \sqrt{-1}\begin{pmatrix}0&1\\-1&0\end{pmatrix}\\
M_{-\beta} & = & \sqrt{-1}\begin{pmatrix}0&2+\sqrt{3}\\-2+\sqrt{3}&0\end{pmatrix}\\
M_\infty & = & -\frac{1}{\sqrt{6}}\begin{pmatrix}\sqrt{3}&6+3\sqrt{3}\\-6+3\sqrt{3}&-\sqrt{3}\end{pmatrix}
\end{eqnarray*}

satisfying the relationship $M_{-1}.M_{-\alpha}.M_{-2}.M_{-\beta}.M_\infty = id$.

The projective monodromy group is generated by
\begin{eqnarray*}
	\overline{M_{-1}} & = & \frac{1}{\sqrt{6}}(3i+6j+k)\\
	\overline{M_{-\alpha}} & = & i + 2j\\
	\overline{M_{-2}} & = & j\\
	\overline{M_{-\beta}} & = & 2j+k\\
	\overline{M_\infty} & = & \frac{1}{\sqrt6}(i+6j+3k)
\end{eqnarray*}
in $\left(\frac{3,-1}{\QQ[\sqrt{6}]}\right)$.
\end{Eqn}

\begin{Eqn}{V}
$$
p(\lambda)\frac{d^2}{d\lambda^2}
+
\frac{1}{2}p^\prime(\lambda)\frac{d}{d\lambda}
+
\left(\frac{1}{18}\lambda+\frac{43}{144}\right)
$$
where $p(\lambda) = (\lambda+5)(\lambda+8)(\lambda+\frac{40}{9})$

After the substitution $\lambda = \frac{\mu-157}{27}$, this differential equation is of Lam\'e type with $p(\mu)=4(\mu+59)(\mu-22)(\mu-37)$, $n=-\frac{1}{3}$, $B = \frac{95}{36}$.

$$
\left\{\newcommand\T{\rule{0pt}{2.3ex}}
\newcommand\B{\rule[-1.0ex]{0pt}{0pt}}
\centering
\begin{tabular}{c c c c}
	\T\B $-59$ & $22$ & $37$ & $\infty$ \\
	\T\B $0$   & $0$  & $0$  & $-\frac{5}{3}$\\
	\T\B $\frac{1}{2}$ & $\frac{1}{2}$ & $\frac{1}{2}$ & $\frac{1}{6}$
\end{tabular}\right\}
$$

The monodromy group is generated by
\begin{eqnarray*}
M_{-59} & = & -\frac{\sqrt{-1}}{2\sqrt{2}}\begin{pmatrix}\sqrt{10} & -2\sqrt{3}+\sqrt{30}\\-2\sqrt{3}-\sqrt{30}& -\sqrt{10}\end{pmatrix}\\
M_{22} & = & -\frac{\sqrt{-1}}{2}\begin{pmatrix}\sqrt{2}&\sqrt{6}\\-\sqrt{6}&-\sqrt{2}\end{pmatrix}\\
M_{37} & = & -\sqrt{-1}\begin{pmatrix}0&1\\-1&0\end{pmatrix}\\
M_\infty & = & -\frac{\sqrt{-1}}{2}\begin{pmatrix}\sqrt{3} & 3-\sqrt{10}\\ 3+\sqrt{10} & \sqrt{3}\end{pmatrix}
\end{eqnarray*}
satisfying $M_{-59}.M_{22}.M_{37}.M_\infty = id$.

Hence the projective monodromy group is generated by
\begin{eqnarray*}
\overline{M_{-59}} & = & \frac{1}{2\sqrt{2}}(i-2j+k)\\
\overline{M_{22}} & = & \frac{1}{2\sqrt{5}}(i+k)\\
\overline{M_{37}} & = & \frac{1}{\sqrt{30}}k\\
\overline{M_\infty} & = & \frac{1}{2\sqrt{3}}(3+3j-k)
\end{eqnarray*}
in $\left(\frac{10,3}{\QQ[\sqrt{2},\sqrt{3},\sqrt{5}]}\right)$.
\end{Eqn}

\begin{Eqn}{VI, VII, VIII, IX}
$_2F_1(\frac{1}{8}, \frac{1}{8}, \frac{3}{4})$. 

The monodromy group is generated by the transformations
\begin{eqnarray*}
	M_0 & = & \frac{1+\sqrt{-1}}{\sqrt{2}}\begin{pmatrix}\sqrt{2}&1\\-1&0\end{pmatrix}\\
	M_1 & = & \sqrt{-1}\begin{pmatrix}0&1\\-1&0\end{pmatrix}\\
	M_\infty & = & \frac{1+\sqrt{-1}}{\sqrt{2}}\begin{pmatrix}1&\sqrt{2}\\0&1\end{pmatrix}
\end{eqnarray*}
satisfying $M_0.M_1.M_\infty = id$.

Thus, the projective monodromy group is generated by $\overline{M}_1 = \begin{pmatrix}0&1\\-1&0\end{pmatrix}$ and $\overline{M}_\infty = \begin{pmatrix}1&\sqrt{2}\\0&1\end{pmatrix}$ and up to a choice of basis is equal to the Fricke modular group of level 2, labeled $\Gamma_0(2)^+$ (see \cite{Dol}).
\end{Eqn}

\begin{Eqn}{X}
$_2F_1(\frac{1}{4}, \frac{1}{4}, 1)$

The monodromy group is generated by the transformations
\begin{eqnarray*}
	M_0 & = & \begin{pmatrix}2&1\\-1&0\end{pmatrix}\\
	M_1 & = & \sqrt{-1}\begin{pmatrix}0&1\\-1&0\end{pmatrix}\\
	M_\infty & = & \sqrt{-1}\begin{pmatrix}1&2\\0&1\end{pmatrix}
\end{eqnarray*}
satisfying $M_0.M_1.M_\infty = id$.

Thus, the projective monodromy group is generated by $\overline{M}_1 = \begin{pmatrix}0&1\\-1&0\end{pmatrix}$ and $\overline{M}_\infty = \begin{pmatrix}1&2\\0&1\end{pmatrix}$ which, up to a suitable choice of basis, is equal to $\Gamma_0(2)$.
\end{Eqn}

\begin{Eqn}{XI}
$$
p(\lambda)\frac{d^2}{d\lambda^2} +
\frac{1}{2}p^\prime(\lambda)\frac{d}{d\lambda} +
\left(\frac{3}{16}\lambda^2+\frac{3}{16}\lambda+\frac{1}{24}\right)
$$
where $p(\lambda)=(\lambda+1)(\lambda+\frac{1}{2})(\lambda+\frac{1}{3})(\lambda+\frac{1}{4})$

$$
\left\{\newcommand\T{\rule{0pt}{2.3ex}}
\newcommand\B{\rule[-1.0ex]{0pt}{0pt}}
\centering
\begin{tabular}{c c c c c}
	\T\B $-1$ & $-\frac{1}{2}$ & $-\frac{1}{3}$ & $-\frac{1}{4}$ & $\infty$ \\
	\T\B $0$  & $0$            & $0$            & $0$            & $-\frac{11}{4}$\\
	\T\B $\frac{1}{2}$ & $\frac{1}{2}$ & $\frac{1}{2}$ & $\frac{1}{2}$ & $\frac{3}{4}$
\end{tabular}\right\}
$$

The monodromy group is generated by the transformations
\begin{eqnarray*}
	M_\infty & = & \frac{1}{\sqrt{2}}\begin{pmatrix}0&-1+\sqrt{3}\\-1-\sqrt{3}&0\end{pmatrix}\\
	M_{-1} & = & -\sqrt{-1}\begin{pmatrix}1&-1+\sqrt{3}\\-1-\sqrt{3}&-1\end{pmatrix}\\
	M_{-\frac{1}{2}} & = & -\sqrt{-1}\begin{pmatrix}\sqrt{3}&2\\-2&-\sqrt{3}\end{pmatrix}\\
	M_{-\frac{1}{3}} & = & -\sqrt{-1}\begin{pmatrix}1&1+\sqrt{3}\\1-\sqrt{3}&-1\end{pmatrix}\\
	M_{-\frac{1}{4}} & = & -\frac{\sqrt{-1}}{\sqrt{2}}\begin{pmatrix}0&1+\sqrt{3}\\1-\sqrt{3}&0\end{pmatrix}\\
\end{eqnarray*}
satisfying $M_\infty M_{-1}M_{-\frac{1}{2}}M_{-\frac{1}{3}}M_{-\frac{1}{4}} = id$. The projective monodromy group is generated by 
\begin{eqnarray*}
\overline{M_\infty} & = & \frac{1}{\sqrt{6}}(3j-k)\\
\overline{M_{-1}} & = & \frac{1}{\sqrt{3}}(i+3j-k)\\
\overline{M_{-\frac{1}{2}}} & = & i+2j\\
\overline{M_{-\frac{1}{3}}} & = & \frac{1}{\sqrt{3}}(i+3j+k)\\
\overline{M_{-\frac{1}{4}}} & = & \frac{1}{\sqrt{6}}(3j+k)
\end{eqnarray*}
in $\left(\frac{3,-1}{\QQ[\sqrt{2},\sqrt{3}]}\right)$.
\end{Eqn}

\begin{Eqn}{XII}
$$
(4\lambda^2-1)(3\lambda^2-1)\frac{d^2}{d\lambda^2} + \lambda(24\lambda^2-7)\frac{d}{d\lambda} + \frac{1}{4}(9\lambda^2-2).
$$
After substituting $\mu = \lambda^2$, we get
$$
p(\mu)\frac{d^2}{d\mu^2} + \frac{1}{2}p^\prime(\mu)\frac{d}{d\mu} + \frac{3}{64}\mu-\frac{1}{96}
$$
with $p(\mu) = \mu\left(\mu-\frac{1}{4}\right)\left(\mu-\frac{1}{3}\right)$.

After the substitution $\mu = \frac{\nu+7}{36}$, this is the Lam\'{e} differential equation with $p(\nu) = 4(\nu+7)(\nu-2)(\nu-5)$, $n = -\frac{1}{4}$, and $B = \frac{3}{16}$.

$$
\left\{\newcommand\T{\rule{0pt}{2.3ex}}
\newcommand\B{\rule[-1.0ex]{0pt}{0pt}}
\centering
\begin{tabular}{c c c c}
	\T\B $0$ & $\frac{1}{4}$ & $\frac{1}{3}$ & $\infty$ \\
	\T\B $0$ & $0$  & $0$  & $-\frac{13}{8}$\\
	\T\B $\frac{1}{2}$ & $\frac{1}{2}$ & $\frac{1}{2}$ & $\frac{1}{8}$
\end{tabular}\right\}
$$

The monodromy group is generated by the transformations
\begin{eqnarray*}
	M_0 & = & \sqrt{-1}\begin{pmatrix}1&-1+\sqrt{3}\\-1-\sqrt{3}&-1\end{pmatrix}\\
	M_{1/4} & = & \frac{\sqrt{-1}}{\sqrt{2}}\begin{pmatrix}0&-1+\sqrt{3}\\-1-\sqrt{3}&0\end{pmatrix}\\
	M_{1/3} & = & \sqrt{-1}\begin{pmatrix}0&1\\-1&0\end{pmatrix}\\
	M_\infty &= & \frac{\sqrt{-1}}{\sqrt{2}}\begin{pmatrix}1+\sqrt{3}&2\\-2&1-\sqrt{3}\end{pmatrix}
\end{eqnarray*}
satisfying $M_0.M_{1/4}.M_{1/3}.M_\infty = id$.

Thus, the projective monodromy group is generated by
\begin{eqnarray*}
\overline{M}_0 & = & \frac{1}{\sqrt{3}}(i+3j-k),\\
\overline{M}_{1/4} & = & \frac{1}{\sqrt{6}}(3j-k),\\
\overline{M}_{1/3} & = & j,\\
\overline{M}_\infty & = & \frac{1}{\sqrt{2}}(1+i+2j)
\end{eqnarray*}
in $\left(\frac{3,-1}{\QQ[\sqrt{2},\sqrt{3}]}\right)$.
\end{Eqn}

\begin{Eqn}{XIII}
$$
p(\lambda)\frac{d^2}{d\lambda^2} +
\frac{1}{2}p^\prime(\lambda)\frac{d}{d\lambda} +
\left(\frac{3}{16}\lambda^2+\frac{3}{4}\lambda-\frac{2}{3}\right)
$$
where
$$
p(\lambda) = \lambda(\lambda-3)(\lambda+5)(\lambda+\frac{16}{3}).
$$

$$
\left\{\newcommand\T{\rule{0pt}{2.3ex}}
\newcommand\B{\rule[-1.0ex]{0pt}{0pt}}
\centering
\begin{tabular}{c c c c c}
	\T\B $-\frac{16}{3}$ & $-5$ & $0$ & $3$ & $\infty$ \\
	\T\B $0$  & $0$            & $0$            & $0$            & $-\frac{5}{4}$\\
	\T\B $\frac{1}{2}$ & $\frac{1}{2}$ & $\frac{1}{2}$ & $\frac{1}{2}$ & $\frac{1}{4}$
\end{tabular}\right\}
$$

The monodromy group is generated by the transformations
\begin{eqnarray*}
M_\infty & = & \begin{pmatrix}
-2&-\sqrt{5}\\
\sqrt{5}&2
\end{pmatrix}\\
M_{-\frac{16}{3}} & = & \sqrt{-1}\begin{pmatrix}
0&-1\\
1&0
\end{pmatrix}\\
M_{-5} & = & \frac{\sqrt{-1}}{\sqrt{2}}\begin{pmatrix}
0&1-\sqrt{3}\\
1+\sqrt{3}&0
\end{pmatrix}\\
M_0 & = & \sqrt{-1}\begin{pmatrix}
-\sqrt{5}&-3+\sqrt{3}\\
3+\sqrt{3}&\sqrt{5}
\end{pmatrix}\\
M_3 & = & \frac{\sqrt{-1}}{\sqrt{2}}\begin{pmatrix}
-2\sqrt{5} & -5+\sqrt{3}\\
5+\sqrt{3}&2\sqrt{5}
\end{pmatrix}
\end{eqnarray*}
satisfying $M_\infty.M_{-\frac{16}{3}}.M_{-5}.M_0.M_3 = id$.

\begin{eqnarray*}
	\overline{M}_\infty & = & \frac{1}{\sqrt{3}}(2i+k)\\
	\overline{M}_{-\frac{16}{3}} & = & \frac{1}{\sqrt{15}}k\\
	\overline{M}_{-5} & = & \frac{1}{\sqrt{10}}(j-k)\\
	\overline{M}_0 & = & \frac{1}{\sqrt{15}}(5i-3j+3k)\\
	\overline{M}_3 & = & \frac{1}{\sqrt{30}}(10i-3j+5k)
\end{eqnarray*}
in $\left(\frac{3,5}{\QQ[\sqrt{2},\sqrt{3},\sqrt{5}]}\right)$.
\end{Eqn}

\begin{Eqn}{XIV}
$$
(\lambda^2-1)(\lambda^2+80)\frac{d^2}{d\lambda^2} + 3\lambda(\lambda^2+26)\frac{d}{d\lambda} + \frac{3(\lambda^4+85\lambda^2+160)}{4(\lambda^2+80)}.
$$
After substituting $\mu = \lambda^2$, we get
$$
p(\mu)\frac{d^2}{d\mu^2} + \frac{1}{2}p^\prime(\mu)\frac{d}{d\mu} + \frac{3}{4}(\mu^2+85\mu+160)
$$
where $p(\mu) = 4\mu(\mu-1)(\mu+80)^2$.

$$
\left\{\newcommand\T{\rule{0pt}{2.3ex}}
\newcommand\B{\rule[-1.0ex]{0pt}{0pt}}
\centering
\begin{tabular}{c c c c c}
	\T\B $-80$ & $0$ & $1$ & $\infty$ \\
	\T\B $-\frac{1}{12}$  & $0$            & $0$          & $-\frac{5}{4}$\\
	\T\B $\frac{1}{12}$ & $\frac{1}{2}$ & $\frac{1}{2}$ & $\frac{1}{4}$
\end{tabular}\right\}
$$

The monodromy group is generated by the transformations
\begin{eqnarray*}
	M_{-80} & = &
	\frac{1}{2}\begin{pmatrix}
		8+\sqrt{3}&4\sqrt{5}+\sqrt{15}\\
		-4\sqrt{5}+\sqrt{15}&-8+\sqrt{3}
	\end{pmatrix}\\
	M_0 & = & \sqrt{-1}\begin{pmatrix}
		2&\sqrt{5}\\
		-\sqrt{5}&-2
	\end{pmatrix}\\
	M_1 & = & \frac{\sqrt{-1}}{2}\begin{pmatrix}
		\sqrt{5}&3\\
		-3&-\sqrt{5}
	\end{pmatrix}\\
	M_\infty & = & \begin{pmatrix}
		-\sqrt{5}&-3-\sqrt{3}\\
		3-\sqrt{3}&\sqrt{5}
	\end{pmatrix}
\end{eqnarray*}
satisfying $M_{-80}.M_0.M_1.M_\infty = id$. The projective monodromy group is generated by
\begin{eqnarray*}
\overline{M_{-80}} & = & \frac{1}{2\sqrt{3}}(3+8i+3j+4k)\\
\overline{M_0} & = & \frac{1}{\sqrt{3}}(2i+k)\\
\overline{M_1} & = & \frac{1}{2\sqrt{15}}(5i+3k)\\
\overline{M_\infty} & = & -\frac{1}{\sqrt{15}}(5i+3j+3k)
\end{eqnarray*}
in $\left(\frac{3,5}{\QQ[\sqrt{3},\sqrt{5}]}\right)$.

\end{Eqn}

\section{Fundamental Domains}

In this section, we construct the fundamental domains of the projective monodromy groups from the last section. We focus on the cases with compact fundamental domain. As we have noticed, in these examples, the differential equations are of generalised Lam\'e type and the projective monodromy groups are generated by elements of order 2. We make use of a classical result, Poincar\'e's polygon theorem, to find the fundamental domains for these examples.

This theorem is discussed in detail in \cite{Iversen} and is introduced by first letting $P$ denote a compact hyperbolic polygon. By a side-pairing transformation, we mean a hyperbolic isometry $\sigma$, other than the identity, that identifies two sides $s$ and $s^\prime$ of $P$ in such a way that $P\cap\sigma(P) = \sigma(s) = s^\prime$. If $P$ has sides $\{s_1,\ldots,s_k\}$, we let $\sigma_i$ be the transformation pairing $s_i$ to some other side.

If $v_i$ is a vertex of $P$ on side $i$, then $\sigma_i(v_i)$ is a vertex on side $s_j$. Similarly, $\sigma_j(\sigma_i(v_i))$ is a vertex on side $s_k$ and so on. Since there are finitely many vertices, we eventually come back to the start and find
$$
\left(\prod_{\alpha\in A}\sigma_\alpha\right)v_i = v_i$$
for some $A\subset\{1\ldots,k\}$. We call $\gamma_{v_i} := \prod_{\alpha\in A}\sigma_\alpha$ a vertex transformation associated to $v_i$.

Letting $\textnormal{angle}(v_\alpha)$ denote the interior angle of $P$ at the vertex $v_\alpha$, then we define the angle sum
$$
\textnormal{sum}(\gamma_{v_i}) = \sum_{\alpha\in A}\textnormal{angle}(v_\alpha)
$$
of the vertex transformation $\gamma_{v_i}=\prod_{\alpha\in A}\sigma_\alpha$. We can now state the theorem.
\begin{theorem}[Poincar\'e's Polygon Theorem]
If each vertex transformation satisfies
$$
\textnormal{order}(\gamma_v).\textnormal{sum}(\gamma_v) = 2 \pi,
$$
then the side pairing transformations $\sigma_1,\ldots,\sigma_k$ generate a Fuchsian group with fundamental domain $P$.
\end{theorem}
This theorem is usually used to find examples of Fuchsian groups corresponding to a given polygon. We shall construct a polygon from a set of generators of a projective monodromy group and use Poincar\'e's theorem to prove that it is a fundamental domain. Let $\Gamma$ be a Fuchsian group generated by order 2 rotations $\sigma_1,\ldots,\sigma_k$. Assume that the convex hull of the fixed points of the generators has these fixed points only at vertices and label the generators so that the index increases in an anticlockwise direction. Let $\sigma_0 = (\prod_{i=1}^k\sigma_i)^{-1}$ and assume that $\sigma_0$ has a fixed point $v\in\mathcal{H}$. We define a polygon $P_\Gamma$ as the polygon with vertices $v$, $\sigma_1(v)$, $\sigma_2\sigma_1(v)$, $\ldots$, $\sigma_k\sigma_{k-1}\ldots\sigma_1(v)=\sigma_0^{-1}(v)=v$. Our assumptions assure that this construction makes sense.

The polygon $P_\Gamma$ has $k$ sides, each of whose midpoints are a fixed point for one of the generators. Let $s_i$ denote the side between $\sigma_{i-1}\ldots\sigma_1(v)$ and $\sigma_i\ldots\sigma_1(v)$. Then the midpoint of $s_i$ is the fixed point of $\sigma_i$. The side-pairing transformations are the involutions $\sigma_i$ $i=1\ldots,k$ identifying $s_i$ with itself.

In our special case, where $\Gamma$ is generated by elements of order 2, it is particularly easy to check the conditions of Poincar\'e's polygon theorem. To show that $P_\Gamma$ is a fundamental domain for $\Gamma$, we need only check the condition
$$
\textnormal{sum}(\gamma_v) = \frac{2\pi}{\textnormal{order}(\gamma_v)}
$$
for the single vertex transformation $\gamma_v = \sigma_0$ since in this situation all vertex transformations are conjugate to $\sigma_0$ or $\sigma_0^{-1}$.

\begin{figure}[ht]
	\centering
		\includegraphics{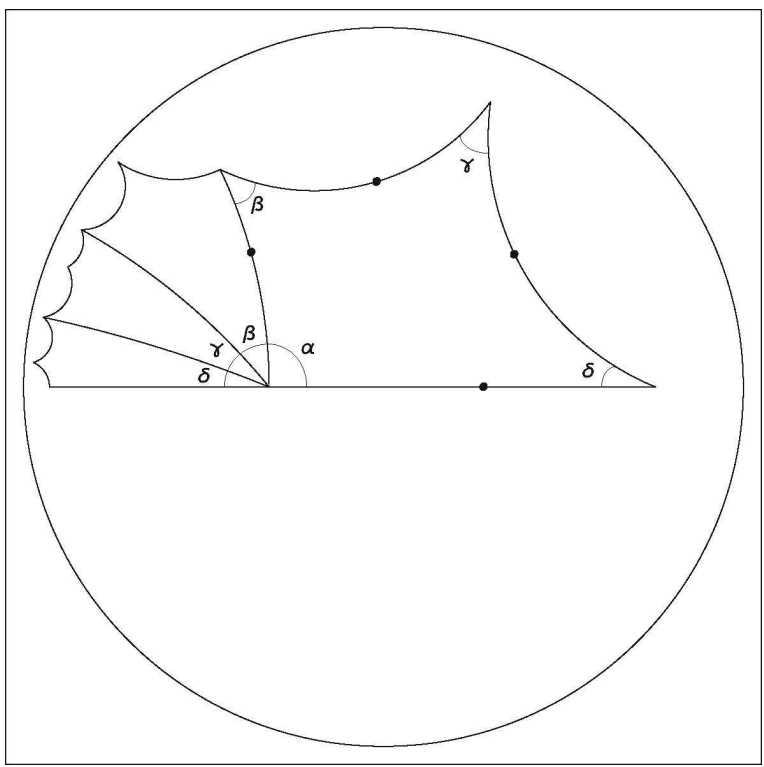}
	\caption{Poincar\'e's Condition}
	\label{fig:Polygon}
\end{figure}

To show that this condition holds, notice that $\frac{2\pi}{\textnormal{order}(\sigma_0)}$ is just the angle subtended by the rotation $\sigma_0$ at its fixed point $v$. On the other hand, $\textnormal{sum}(\sigma_0)$ is the sum of the interior angles of $P_\Gamma$. To show that these values are the same, write $\tau_m = \sigma_1\ldots\sigma_m$ for $m =1,\ldots,k$. The expression
\begin{eqnarray*}
\sigma_0^{-1} & = & \sigma_1\sigma_2\ldots\sigma_k\\
& = & (\tau_{k-1}\sigma_k\tau_{k-1}^{-1})(\tau_{k-2}\sigma_{k-1}\tau_{k-2}^{-1})\ldots(\tau_1\sigma_2\tau_1^{-1})\sigma_1
\end{eqnarray*}
shows that it is possible to decompose the vertex transformation $\sigma_0^{-1}$ (or $\sigma_0$) into a composition of rotations conjugate to the generators (or their inverses). This decomposition corresponds to the rearrangement of the angles in figure~\ref{fig:Polygon}. The marked points in this diagram show the fixed points of the generators of $\Gamma$.

Since the conditions of Poincar\'e's polygon theorem are satisfied, the polygon $P_\Gamma$ defined above is a fundamental domain for $\Gamma$.

It is customary to draw non--compact fundamental domains in the upper half--plane model of $\mathcal{H}^2$. However, compact fundamental domains do not have cusps and so there is no special boundary point to put at $\infty$ in the upper half--plane. For this reason, we choose to draw our compact fundamental domains in the Poincar\'e disk model. This better displays the symmetries of the tiling. For consistency, we also draw the non--compact fundamental domains in the disk model.

\begin{figure}
	\centering
		\includegraphics{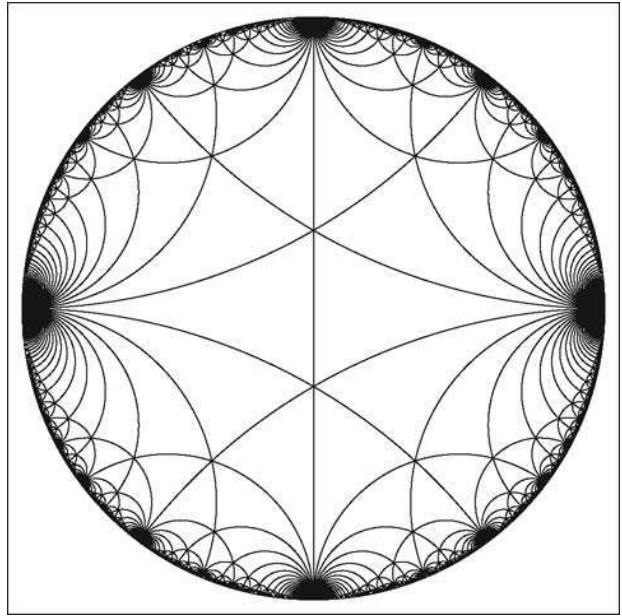}
	\caption{The Fundamental Domain for Examples I, II and III}
	\label{fig:123}
\end{figure}

\begin{figure}
	\centering
		\includegraphics{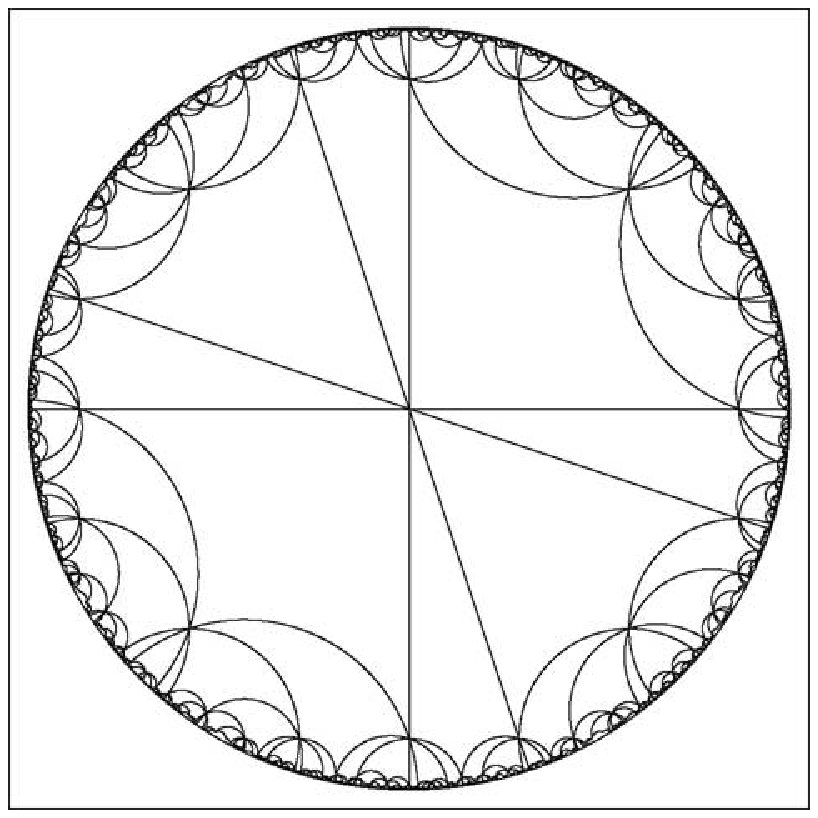}
	\caption{The Fundamental Domain for Example IV}
	\label{fig:4}
\end{figure}

\begin{figure}
	\centering
		\includegraphics{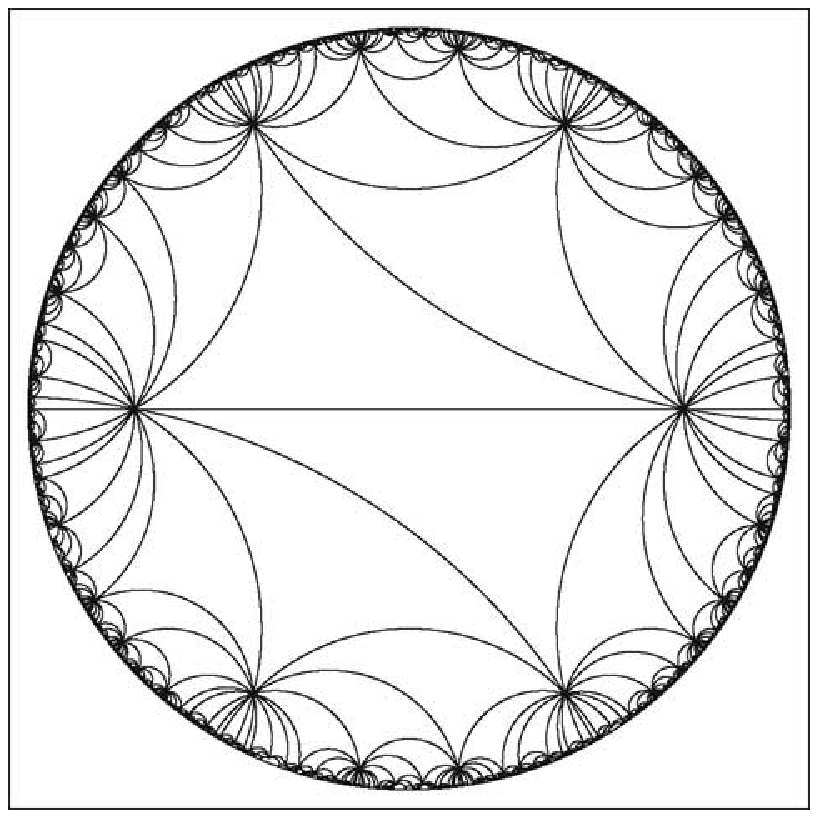}
	\caption{The Fundamental Domain for Example V}
	\label{fig:5}
\end{figure}

\begin{figure}
	\centering
		\includegraphics{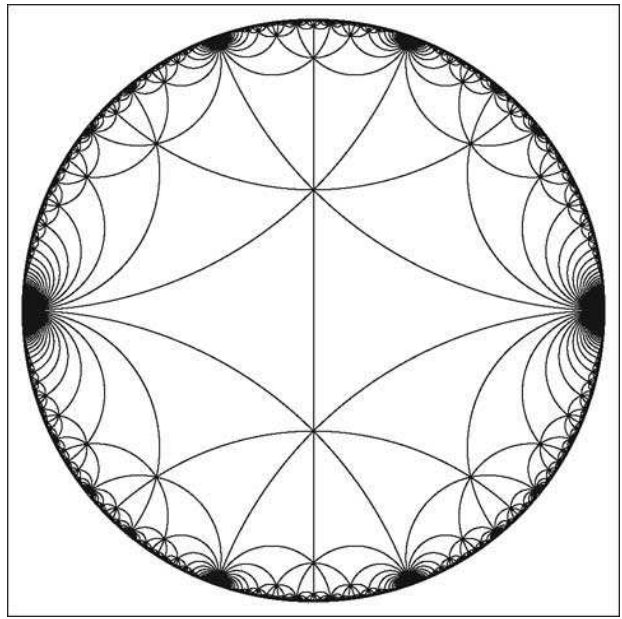}
	\caption{The Fundamental Domain for Examples VI, VII, VIII and IX}
	\label{fig:6789}
\end{figure}

\begin{figure}
	\centering
		\includegraphics{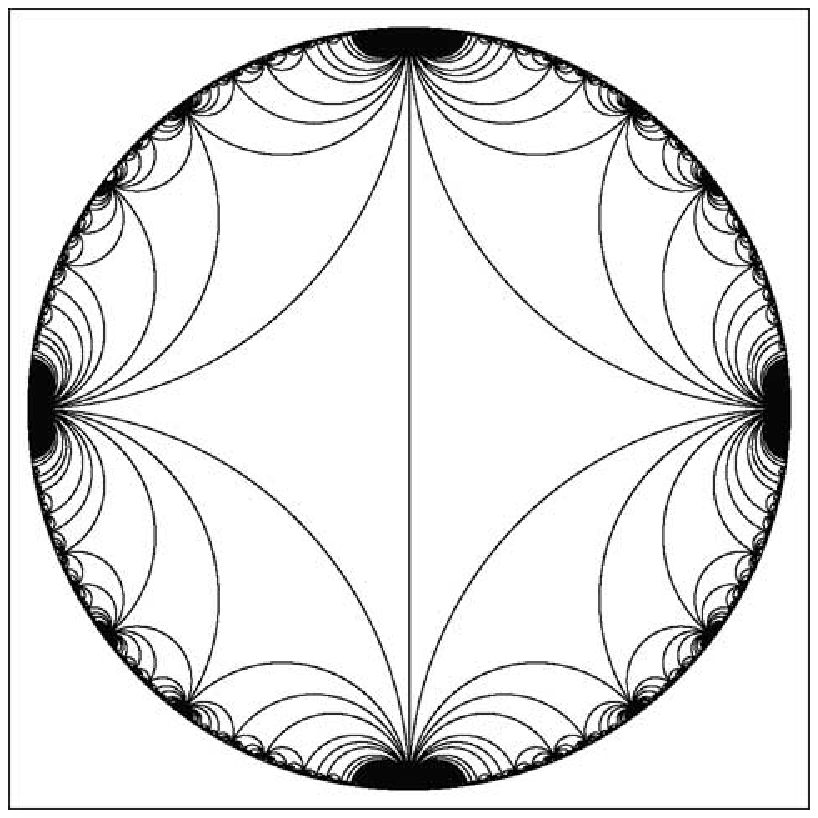}
	\caption{The Fundamental Domain for Example X}
	\label{fig:10}
\end{figure}

\begin{figure}
	\centering
		\includegraphics{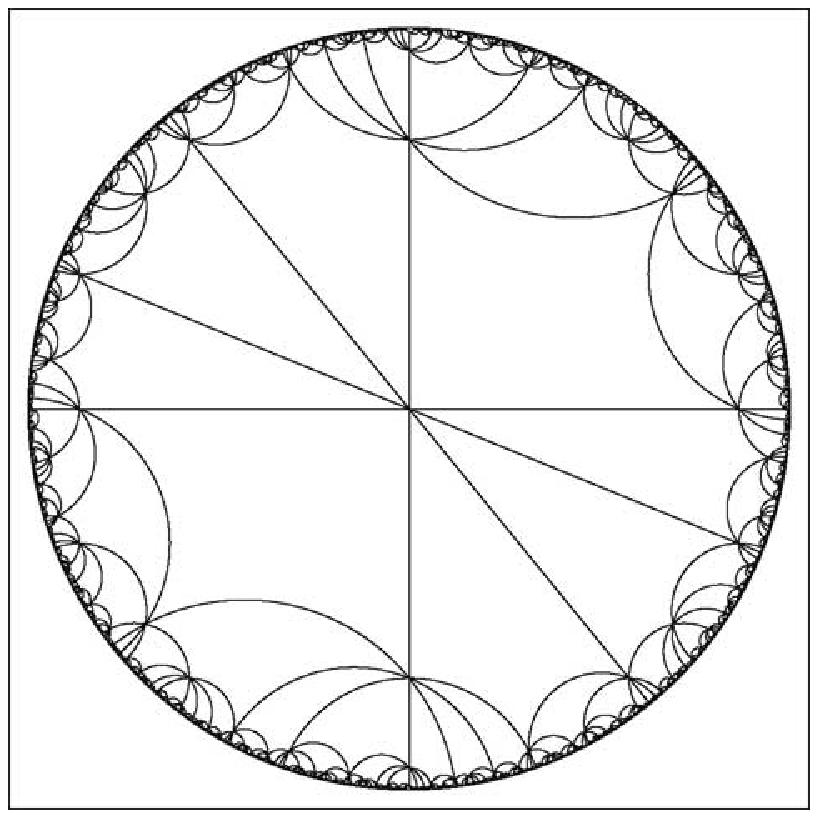}
	\caption{The Fundamental Domain for Example XI}
	\label{fig:11}
\end{figure}

\begin{figure}
	\centering
		\includegraphics{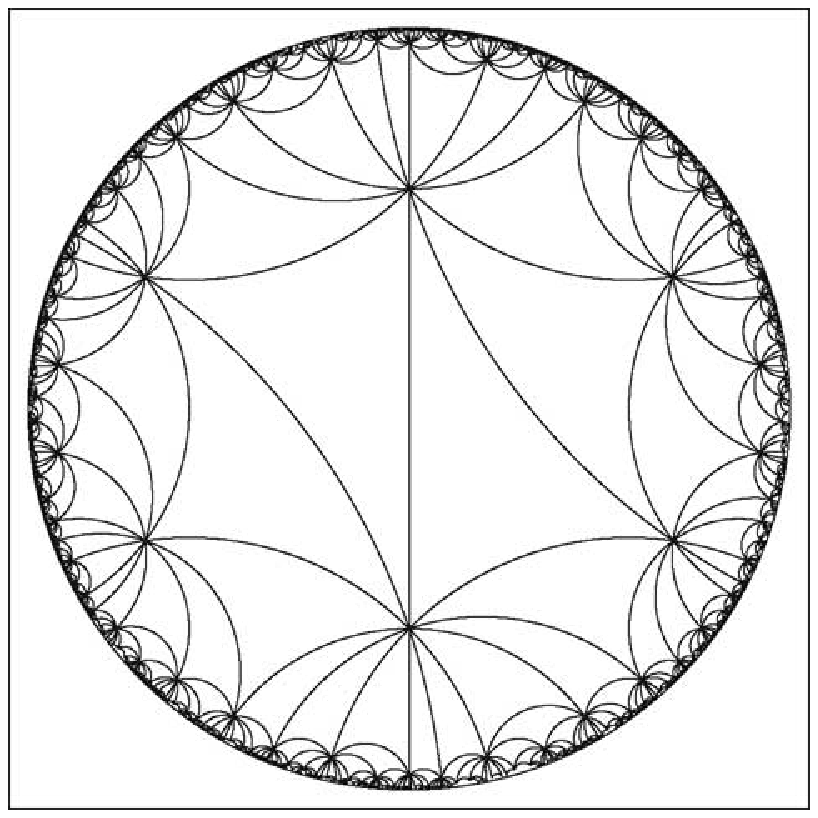}
	\caption{The Fundamental Domain for Example XII}
	\label{fig:12}
\end{figure}

\begin{figure}
	\centering
		\includegraphics{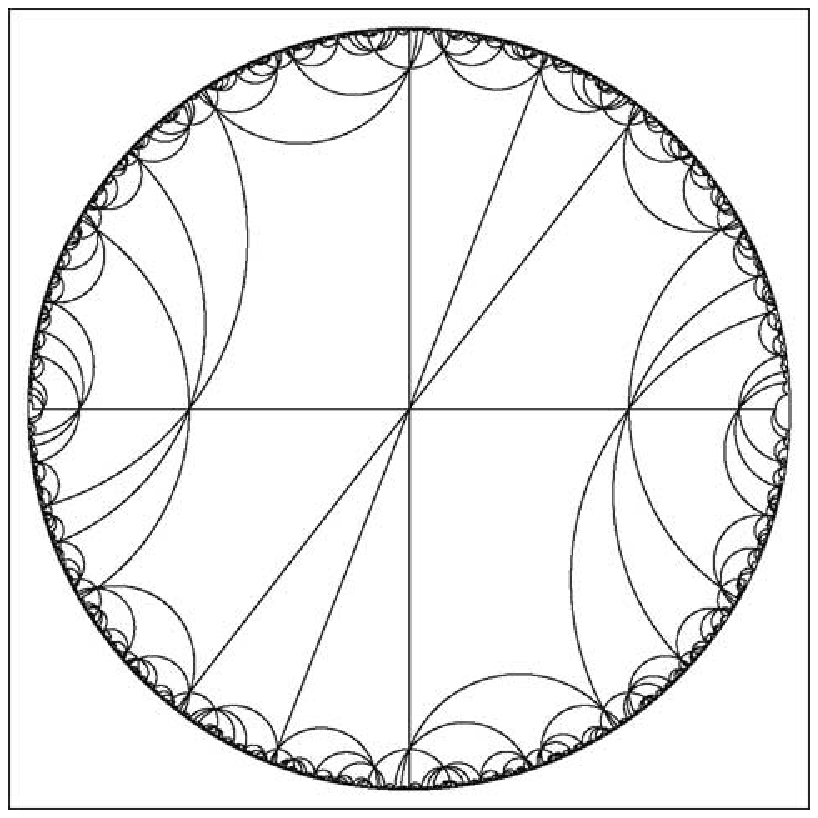}
	\caption{The Fundamental Domain for Example XIII}
	\label{fig:13}
\end{figure}

\begin{figure}
	\centering
		\includegraphics{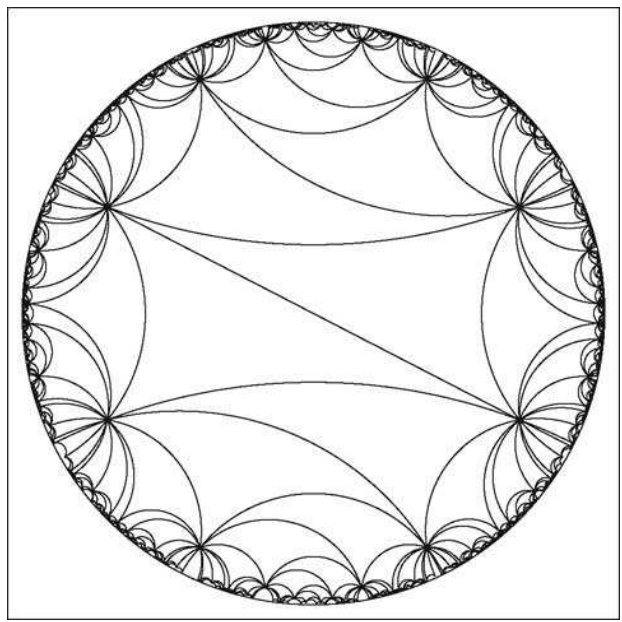}
	\caption{The Fundamental Domain for Example XIV}
	\label{fig:14}
\end{figure}

Figures \ref{fig:123}---\ref{fig:14} show the fundamental domains for the monodromy groups that have occurred in our examples. These were drawn in Mathematica using a hyperbolic geometry package, \cite{hypGeomPackage}.


	\newpage
	
\chapter{Quotients of Symmetric Surfaces}\label{chapter:quots}

We turn our attention to the following situation. Suppose that $X$ is an algebraic surface with a finite group of automorphisms $G\subset\textnormal{Aut}(X)$ with the property that $X/G$ has at worst DuVal singularities. Let $Y$ denote the minimal resolution of $X/G$. When is $Y$ a K3 surface?

A simple example is the case when $X$ is itself a K3 surface and $G$ is a finite group of symplectic automorphisms of $X$.

Suppose that $X$ is a smooth hypersurface in $\PP^3$, and $G\subset\textnormal{Sl}(4,\CC)$ is a finite group inducing an action of $\overline{G}\subset\PP\textnormal{Sl}(4,\CC)$ on $X$. Then $X/\overline{G}$ will have at worst DuVal singularities since it has only isolated fixed points on $X$ and at the fixed points the action is analytically equivalent to a finite group $H\subset\textnormal{Sl}(2,\CC)$ acting on $\CC^2$. In \cite{BS}, Barth and Sarti consider three families of hypersurfaces in $\PP^3$ each invariant under some finite subgroup of $\PP\textnormal{SO}(4)\subset\PP\textnormal{Sl}(4,\CC)$. After a detailed examination of the group action and the base locus of the families, it is determined that the minimal resolution of the quotient surfaces are K3 surfaces. Furthermore, the configurations of exceptional curves from the resolution are determined and it is observed that the K3 surfaces have generic Picard number 19.

We shall describe an alternative method for establishing whether a quotient surface is K3 and use it to find another example of this kind.

\section{Sarti's Examples}
As Lie groups, neither $SO(3)$ nor $SO(4)$ are simply connected, but have the special unitary group $SU(2)$ and $SU(2)\times{SU(2)}$ as their respective universal (double) covers. Hence, we obtain the following 2--to--1 homomorphisms.
$$\begin{array}{c c c}
                             &                              &   SO(3)\times{SO}(3)           \\
                             &          \nearrow            &                                \\
    SU(2)\times{SU(2)}       &                              &                                \\
                             &          \searrow            &                                \\
                             &                              &   SO(4)\phantom{\times SO(3)}  \\

\end{array}$$

Via pull-back and push-forward, two subgroups $A,\,B\subset{SO(3)}$ provide a group labeled $AB\subset{SO(4)}$ satisfying $AB / \{\pm{I}\} \cong A\times{B}$. Using this construction, it is possible to classify all finite subgroups of $SO(4)$ (see for example \cite{CS}).

We shall make use of a classical result, known as Molien's theorem, to calculate the Hilbert series of a finite matrix group, $G$. By definition, the Hilbert series is the power series whose $n$th coefficient is the dimension of the vector space of homogeneous $G$--invariant polynomials of degree $n$.
\begin{theorem}[Molien's Theorem]
If $G\subset\textnormal{Gl}(n,\CC)$ is a finite subgroup, then
$$
P(G,t):=\sum_{k=0}^\infty{(\dim{\CC[x_1,\ldots,x_n]^G})\,t^k}=\frac{1}{|G|}\sum_{A\in G}\frac{1}{\det(I_n - tA)}.
$$
where $I_n$ is the identity matrix.\qed
\end{theorem}

Since the invariant ring is always finitely generated, it is well--known that the Hilbert series of a finite matrix group is rational. The Hilbert series has an expression of the form
$$
P(G,t)=\frac{p(t)}{\prod{_{i=1}^n(1-t^{a_i})}}
$$
encoding the fact that the graded ring of invariant polynomials is generated by homogeneous elements of degree $a_i$. In our examples, the generators of the ring of invariants will often satisfy one or two polynomial relations and it will be easy to find this rational expression. The calculation of the Hilbert series can be backed up by an explicit determination of the generators and their relations using Magma. The calculation of the Hilbert series of a finite matrix group is conveniently included as a single function in Magma and we include some sample code in figure~\ref{magma} to demonstrate this.

This code should be self explanatory, but we give a quick summary of the main points. The first few inputs specify a matrix group $\textnormal{Sym}(4)\cong{G}\subset\textnormal{SO}(4)$ and define the graded ring of invariants. A set of generators is calculated using the \verb|FundamentalInvariants()| command and their relations by the \verb|Relations()| command. We see in this case, that the invariant ring is generated by elements of degree $1, 2, 4, 6$ and $9$ with one relation of degree $18$. The final line confirms that the Hilbert series is
$$
P(G,t) = \frac{1-t^{18}}{(1-t)(1-t^2)(1-t^4)(1-t^6)(1-t^9)}.
$$

\begin{figure}
\begin{quote}
\small
\begin{verbatim}
> Q := Rationals();
> G := MatrixGroup< 4,Q | A,B >
  where A is
  [1,0,0,0,
   0,0,1,0,
   0,0,0,1,
   0,1,0,0]
  where B is
  [1, 0,0,0,
   0,-1,0,0,
   0, 0,0,1,
   0, 0,1,0];
> Order(G);   // In fact, G is Sym(4)
24
> R := InvariantRing(G);
> FundamentalInvariants(R);
[
    x1,
    x2^2 + x3^2 + x4^2,
    x2^4 + x3^4 + x4^4,
    x2^6 + x3^6 + x4^6,
    x2^5*x3^3*x4 - x2^5*x3*x4^3 - x2^3*x3^5*x4
        + x2^3*x3*x4^5 + x2*x3^5*x4^3 - x2*x3^3*x4^5
]
> Relations(R);
[
    -1/36*f2^9 + 1/3*f2^7*f3 - 5/18*f2^6*f4 - 4/3*f2^5*f3^2
        + 13/6*f2^4*f3*f4 + 11/6*f2^3*f3^3 - 17/18*f2^3*f4^2
        - 25/6*f2^2*f3^2*f4 - 1/4*f2*f3^4 + 7/2*f2*f3*f4^2
        + 1/6*f3^3*f4 - f4^3 - h1^2
]
> H<t> := HilbertSeries(R);
> PS<v> := PowerSeriesRing(Q);
> PS ! ((1-t)*(1-t^2)*(1-t^4)*(1-t^6)*(1-t^9)*H);
1 - v^18
\end{verbatim}
\normalsize
\end{quote}
\caption{Example Magma Session}
\label{magma}
\end{figure}

In \cite{BS}, Barth and Sarti describe three families of surfaces in $\PP^3$ invariant under some finite subgroups of $\PP\textnormal{SO}(4)$. The three groups $\mathfrak{A}_4, \mathfrak{S}_4, \mathfrak{A}_5\subset\textnormal{SO}(3)$ of rotations of a tetrahedron, octahedron and icosahedron produce groups $\mathfrak{A}_4\times{\mathfrak{A}_4},\ \mathfrak{S}_4\times{\mathfrak{S}_4},\ \mathfrak{A}_5\times{\mathfrak{A}_5}$ with a natural action on $\PP^3$. Since $\textnormal{SO}(4)$ acts in a distance preserving manner, any subgroup of $\textnormal{SO}(4)$ leaves the polynomial $q:=x_0^2+x_1^2+x_2^2+x_3^2$ invariant. By Molien's theorem, we calculate the Hilbert series for these group actions:
\begin{eqnarray*}
P(\mathfrak{A}_4\times \mathfrak{A}_4,\,t)&=&\frac{1-2t^{24}+t^{48}}{(1-t^2)(1-t^6)(1-t^8)(1-t^{12})^3} \\
P(\mathfrak{S}_4\times \mathfrak{S}_4,\,t)&=&\frac{1-t^{36}-t^{48}+t^{84}}{(1-t^2)(1-t^8)(1-t^{12})(1-t^{18})(1-t^{24})^2} \\
P(\mathfrak{A}_5\times \mathfrak{A}_5,\,t)&=&\frac{1-t^{120}}{(1-t^2)(1-t^{12})(1-t^{20})(1-t^{30})(1-t^{60})}.
\end{eqnarray*}\label{threegroups}
In particular, in addition to the degree 2 invariant, $q$, we also have invariants of degree 6, 8 and 12 for $G=\mathfrak{A}_4\times{\mathfrak{A}_4},\ \mathfrak{S}_4\times{\mathfrak{S}_4}$, and $\mathfrak{A}_5\times{\mathfrak{A}_5}$ respectively. Barth and Sarti \cite{BS} consider the $G$ invariant pencils of surfaces $X_\lambda\colon q^{n/2}+\lambda s_n=0$ $(n = 6, 8$ and $12)$ and show that after minimally resolving singularities, the quotients $X_\lambda/G$ are pencils of K3 surfaces. Moreover, the Picard lattices of generic members of these pencils are found explicitly in terms of exceptional curves of the resolution and are found to be of rank 19 in each case.

\section{Weighted Projective Spaces and G.I.T.}\label{GIT}

In order to proceed, we need to take a look at some easy geometric invariant theory quotients and weighted projective spaces. All the results on weighted projective spaces we use can be found in \cite{Fletcher}.

\begin{definition}
A weighted projective space $\PP(a_0,a_1,\ldots,a_n)$ is said to be \emph{well-formed} if no $n$ of the $n+1$ weights share a common factor (and in particular, all $n+1$ weights have no common factor).
\end{definition}

Writing $R := \CC[x_0,x_1,\ldots,x_n]$ in its homogeneous parts, $R=\oplus_{m=0}^\infty R_m$, the $i$--th truncation is the graded ring $R^{[i]}:=\oplus_{m=0}^\infty R_{im}$.
\begin{lemma}
If $R$ is a graded ring, then for any truncation, $R^{[i]}$, we have
$$\Proj(R)\cong\Proj\left(R^{[i]}\right).$$
\end{lemma}
\begin{corollary}
Any weighted projective space is isomorphic to one with well-formed weights.
\end{corollary}
\begin{proof}If $n$ of the $n+1$ weights share a common factor $d$, then passing to the $d$--th truncation provides the isomorphism
$$\PP(da_0,\ldots,a_i,\ldots,da_n) = \Proj(R) \cong \Proj(R^{[d]}) = \PP(da_0,\ldots,da_i,\ldots,da_n)$$
so that all the weights may be assumed to be divisible by $d$. However, if all the weights are divisible by $d$, then directly from the definition of weighted projective space, we find
$$\PP(da_0,\ldots,da_n) \cong \PP(a_0,\ldots,a_n).$$
\end{proof}
\begin{definition}
A subvariety $X$ in $\PP(a_0,\ldots,a_n)$ is said to be \emph{well-formed} if the weighted projective space is well-formed and if $X\cap\PP_\textnormal{sing}$ has codimension at least 2 in $X$. Here, $\PP_\textnormal{sing}$ denotes the subvariety of singular points of $\PP(a_0,\ldots,a_n)$.
\end{definition}
There is a convenient characterisation of well-formedness for complete intersections discussed in \cite{Fletcher}. In particular, for a hypersurface of degree $d$ in a well--formed space $\PP(a_0,\ldots,a_n)$, we simply require that for any $i$ and $j$, the highest common factor of $a_0,\ldots,\widehat{a}_i,\ldots,\widehat{a}_j,\ldots,a_n$ divides $d$. For a complete intersection of degree $d_1$, $d_2$, we require that for any distinct $i$, $j$ and $k$, $\textnormal{hcf}(a_0,\stackrel{i,j,k}{\widehat{\ldots}},a_n)$ and $\textnormal{hcf}(a_0,\stackrel{i,j}{\widehat{\ldots}},a_n)$ divide either of $d_1$ or $d_2$.
\begin{definition}
A subvariety $X$ in an $n$ dimensional weighted projective space is said to be \emph{quasismooth} if the cone over $X$ in $\CC^{n+1}$ is smooth away from 0.
\end{definition}
\begin{theorem}[Adjunction Formula]\label{adjunction}
If $X$ is a well--formed, quasismooth complete intersection of degree $d_1,\ldots,d_m$ in the weighted projective space $\PP(a_0,\ldots,a_n)$, then the canonical sheaf, $\omega_X$, is given by
$$
\omega_X \simeq \Oh_X\left(\sum d_i - \sum a_j\right)
$$
\begin{flushright}
$\square$
\end{flushright}
\end{theorem}
If $R$ is a graded ring (with all the usual conventions), and $G$ is a finite
group acting on $R$, then the inclusion
$$R^G \hookrightarrow R$$
induces the quotient morphism
$$\textnormal{Spec}R \twoheadrightarrow{\textnormal{Spec}R^G} =
(\textnormal{Spec}R) / G.$$
For clarity, we consider the case $R = \CC[x_0,\ldots,x_n]$ with $G$ a finite
subgroup of $\textnormal{Gl}(n,\CC)$. Since the orbit of
$0\in\textnormal{Spec}(R) = \Aff^n$ is $\{0\}$, there is no problem in projectivising the quotient. Hence, we can construct the quotient $\PP^n / G$ as
$$\PP^n/G = \Proj(\CC[x_0,\ldots,x_n])^G = \Spec(R)\setminus{\{0\}} / {\CC^*}.$$
In practice, to compute a quotient using this method, the first step is to
determine the graded ring of invariants, $\CC[x_0,\ldots,x_n]^G$. (This is
achieved in Magma by inputting $G$ as a matrix group and using the \verb|InvariantRing()|
command). The quotient morphism $\pi : \PP^n \rightarrow \PP^n / G$ is then given
by
$$\pi : (x_0,\ldots,x_n) \mapsto (s_0,\ldots,s_m),$$
where $s_0,\ldots,s_m$ are homogeneous polynomials generating the invariant
ring. Since
$$(s_0(\underline{x}),\ldots,s_m(\underline{x})) =  (s_0(\lambda\underline{x}),\ldots,s_m(\lambda\underline{x})) = (\lambda^{\alpha_0}s_0(\underline{x}),\ldots,\lambda^{\alpha_m}s_m(\underline{x})),$$
where $\alpha_i = \deg s_i$, this quotient naturally lies in a
weighted projective space:
$$\PP^n / G \subset \PP(\alpha_0,\ldots,\alpha_m).$$
Although this construction depends on an initial choice of generators for $R^G$, different choices give isomorphic quotients related by weighted projective transformations.

\section{Quotient K3 Surfaces}

These results on geometric invariant quotients and weighted projective spaces justify the following procedure to find examples of K3 surfaces occurring as quotients of symmetric surfaces.
\begin{algorithm}A procedure to find examples of K3 surfaces occurring as quotients of invariant surfaces in $\PP^3$.
\begin{enumerate}
    \item Let $G\subset\textnormal{Sl}(4,\CC)$ be a finite subgroup and consider its induced action on $\PP^3$.
     Compute the graded ring, $R^G$, of polynomial invariants of this group action. Say $R^G = \CC[a_1,\ldots,a_n] / (r_1,\ldots,r_m)$ where each $a_i$ is a generating invariant of homogeneous degree $\alpha_i$ and the $r_j$ are the polynomial relations between the generators. We may now express $\PP^3 / G$ as
    $$\PP^3 / G :\ (r_1 = \ldots = r_m = 0) \subset \PP(\alpha_1,\ldots,\alpha_n).$$
    \item Consider each homogeneous piece, $R^G_{(d)}$, of $R^G$. This is a (projectivised) finite dimensional vector space with basis $\{p_i,\ldots,p_k\}$, say, defining the $k-1$ parameter family
    $$\lambda_1p_1+\ldots+\lambda_kp_k=0$$
    of $G$--invariant surfaces in $\PP^3$. We require this family to be generically nonsingular, or at least for the quotient to have only DuVal singularities. Suppose this is the case (if not, then it may be true for some subfamily).
     The quotient of these $G$--invariant surfaces by the induced action of $G$ is then
    $$(\lambda_1p_1+\ldots+\lambda_kp_k = r_1 = \ldots = r_m = 0) \subset \PP(\alpha_1,\ldots,\alpha_n).$$
    \item Check whether this is a complete intersection. If so, then it is isomorphic to a well--formed complete intersection. Reduce to the well--formed expression and check whether it defines a family of K3 surfaces by applying the adjunction formula.
\end{enumerate}
\end{algorithm}

\begin{example}We describe the procedure for the group $\mathfrak{A}_4\times \mathfrak{A}_4\subset\textnormal{SO}(4)$.
\begin{enumerate}
    \item Consider the group $\mathfrak{A}_4 \times \mathfrak{A}_4\subset\PP\textnormal{SO}(4)$. As the Hilbert series on page \pageref{threegroups} suggests, this group action has $6$ fundamental invariants of degrees $2, 6, 8, 12, 12$, and $12$ with two relations between these, both of degree $24$. Hence $\PP^3 / (\mathfrak{A}_4\times \mathfrak{A}_4)$ is of the form
    $$\PP^3 / (\mathfrak{A}_4\times \mathfrak{A}_4):\ (f_{24} = g_{24} = 0) \subset \PP(2,6,8,12,12,12).$$
We don't actually need to specify what $f_{24}$ or $g_{24}$ are, only their degrees are important here.
    \item There are two invariants of degree $6$; the generating invariant of degree $6$, $s_6$, and the cube of the degree $2$ invariant, $q^3 = (x_0^2+x_1^2+x_2^2+x_3^2)^3.$ The resulting one--parameter family of surfaces
    $$(\lambda_1s_3 + \lambda_2\,q^3 = 0) \subset \PP^3$$
    is generically nonsingular and leads to the quotient
    $$(\lambda_1s_3 + \lambda_2\,q^3 = f_{24} = g_{24} = 0) \subset \PP(2,6,8,12,12,12).$$
    \item This complete intersection is not well-formed since, for example, each weight is divisible by $2$. Also, away from the degenerate surface at $\lambda_1 = 0$, we may write $\lambda = \frac{\lambda_2}{\lambda_1}$ and substitute $s_6$ with $-\lambda\,q^3$ wherever it occurs to write our quotient as
    $$(F_{24,\lambda} = G_{24,\lambda} = 0)\subset\PP(2,8,12,12,12)$$
    where $F_{24,\lambda}$ and $G_{24,\lambda}$ are obtained from $f_{24}$ and $g_{24}$ by making the above substitution.
    All the weights should be divided by their common factor, $2$, to get a complete intersection of degree $(12,12)$ in $\PP(1,4,6,6,6)$. The last four weights still have a common factor of $2$ and so
    $$\PP(1,4,6,6,6)\cong\PP(2,4,6,6,6).$$
    Dividing the weights by $2$ again expresses our quotient as a complete intersection of degree $(6,6)$ in $\PP(1,2,3,3,3)$.
    This is well-formed, and since $6+6 = 1+2+3+3+3$, this is a family of K3 surfaces with DuVal singularities.
\end{enumerate}

\end{example}
This procedure may be carried out in all the other examples from \cite{BS} and offers a significantly simplified way to prove that the quotients are K3. We aim to find some more examples of K3 surfaces that are quotients of symmetric surfaces. As mentioned earlier, the K3 surfaces with symplectic automorphisms, as found in chapter \ref{autChapter}, provide a source of such examples. To find more non--trivial examples like those of \cite{BS}, we take a look at the primitive finite subgroups of $\textnormal{Sl}(4,\CC)$ listed in \cite{HH}.

Following the enumeration of subgroups of $\textnormal{Sl}(4,\CC)$ used in \cite{HH}, only groups $\textnormal{II}_*$, $\textnormal{VI}_*$, $\textnormal{X}_*$, $\textnormal{XIV}_*$, $\textnormal{XVI}_*$, $\textnormal{XXVI}_*$, $\textnormal{XXIX}_*$ and $\textnormal{XXX}_*$ provide us with invariant rings that are complete intersection rings (required to apply the adjunction formula). Of these groups, $\textnormal{X}_*$, $\textnormal{XIV}_*$ and $\textnormal{XVI}_*$ correspond to the three groups $\mathfrak{A}_4\times\mathfrak{A}_4$, $\mathfrak{S}_4\times\mathfrak{S}_4$ and $\mathfrak{A}_5\times\mathfrak{A}_5$ considered in \cite{BS}.

Group $\textnormal{II}_*$ is isomorphic to $\mathfrak{A}_5$ and $\textnormal{XXVI}_*$ is the lift of $M_{20}$ from $\PP\textnormal{Sl}(4,\CC)$. These two groups occurred in chapter \ref{autChapter} and have invariant quartic K3 surfaces.

The groups $\textnormal{VI}_*$ and $\textnormal{XXX}_*$ have the following Hilbert series:
$$
P(\textnormal{VI}_*,t) = \frac{1-t^{120}}{(1-t^{12})(1-t^{18})(1-t^{24})(1-t^{30})(1-t^{40})}
$$
$$
P(\textnormal{XXX}_*,t) = \frac{1-t^{64}-t^{120}+t^{184}}{(1-t^8)(1-t^{24})^2(1-t^{32})(1-t^{40})(1-t^{60})}.
$$
It can be verified that neither of these groups yield new examples of K3 quotients.

The final group $G = $ $\textnormal{XXIX}_*$ \emph{does} provide us with a new family of K3 surfaces and we look at this in some extra detail. As an alternative to the generators given in \cite{HH}, the group may be generated as
$$
G = \left<\overline{F}_{384}, \overline{M}_{20} \right> \subset \textnormal{Sl}(4,\CC)
$$
where the groups $\overline{F}_{384}$ and $\overline{M}_{20}$ are the lift to $\textnormal{Sl}(4,\CC)$ of those appearing in section~\ref{p3}. Explicitly,  $\textnormal{XXIX}_*$ is generated by the matrices
$$
\begin{pmatrix}
i&0&0&0\\
0&-i&0&0\\
0&0&1&0\\
0&0&0&1
\end{pmatrix},\ \ 
\begin{pmatrix}
1&0&0&0\\
0&i&0&0\\
0&0&-i&0\\
0&0&0&1
\end{pmatrix},\ \ 
\begin{pmatrix}
1&0&0&0\\
0&1&0&0\\
0&0&i&0\\
0&0&0&-i
\end{pmatrix},\ \ 
\begin{pmatrix}
-i&0&0&0\\
0&1&0&0\\
0&0&1&0\\
0&0&0&i
\end{pmatrix}
$$
$$
\begin{pmatrix}
1&0&0&0\\
0&0&1&0\\
0&0&0&1\\
0&1&0&0
\end{pmatrix},\ \ 
\begin{pmatrix}
0&1&0&0\\
1&0&0&0\\
0&0&1&0\\
0&0&0&-1
\end{pmatrix},\ \ 
\frac{1}{2}\begin{pmatrix}
-i&-i&-1&1\\
-i&i&-1&-1\\
i&i&-1&1\\
i&-i&-1&-1
\end{pmatrix}
$$
and has order $5760\times{4}$. Its Hilbert series is
$$
P(\textnormal{XXIX}_*,t) = \frac{1-t^{120}}{(1-t^8)(1-t^{12})(1-t^{20})(1-t^{24})(1-t^{60})}
$$
indicating the fact that the ring of invariants is generated by invariants $s_n$ of degrees $n = 8, 12, 20, 24$ and $60$ satisfying a relation of degree $120$.

The family of surfaces
$$
X_\lambda : (s_{20} = \lambda s_8s_{12})\subset\PP^3
$$
is invariant under $G = \textnormal{XXIX}_*$. We choose not to write down the invariant polynomials here, although we are required to check that this family is generically nonsingular. To do this, we may use Magma to check that one of the surfaces in the family is nonsingular to conclude that the Zariski open set of base points above which the fibre is nonsingular is nonempty and hence dense.

The quotient family $X_\lambda/G$ is given by
$$
X_\lambda/G : (s_{20} - \lambda s_8s_{12} = f_{120} = 0 ) \subset \PP(8,12,20,24,60).
$$
After substituting $s_{20}$ with $\lambda s_8s_{12}$, we get
$$
X_\lambda/G : (F_{\lambda,120} = 0 ) \subset \PP(8,12,24,60)
$$
where $F_{\lambda,120}(s_8,s_{12},s_{24},s_{60}) = f_{120}(s_8,s_{12},\lambda s_8s_{12},s_{24},s_{60})$. This reduces to the well--formed expression
$$
(F_{\lambda,10} = 0 ) \subset \PP(2,1,2,5)
$$
which describes a family of K3 surfaces with DuVal singularities because $10 = 2+1+2+5$.

	\newpage

\appendix
	
\chapter{Macaulay2 Code to Calculate the Picard--Fuchs Differential Equation}\label{macaulay}
\begin{verbatim}
-----------------------------------------
-- Start Macaulay2 and load this file: --
-- load "...path.../PicardF.m2"        --
-----------------------------------------

R = QQ[a,x,y,z,t];               -- Set up the polynomial ring

------------------------ BEGIN INPUTS ------------------------

Q = value read "Q = "          -- Defining equation of the family
Coeff = {value read "c3 = "}   -- First coefficient of the P-F eqn
Deg = value read "Degree = "   -- Expected degree of the P-F eqn

Exponent = Deg;
q = (diff Q)_(0,0);            -- The derivative with respect to a
if (q==0) then (<<"ERROR: constant family?"<<endl; end)  --Check

----------------- BEGIN FIRST REDUCTION STEP -----------------

I = ideal(
   (diff Q)_(0,1), (diff Q)_(0,2),
   (diff Q)_(0,3), (diff Q)_(0,4)
);
gbI = gb(I, ChangeMatrix=>true); 
if ( not( (Coeff#0)*q^Exponent % gbI == 0) ) then
   (<<"ERROR - BAD COEFFICIENT"<<endl; end)  -- CHECK
Reduced = ( (Coeff#0)*q^Exponent // gbI);
          -- Express c_k*q^k in terms of the partial derivatives of Q
Numerator =  (diff Reduced)_(0,1) + (diff Reduced)_(1,2)
           + (diff Reduced)_(2,3) + (diff Reduced)_(3,4);

--------------------- BEGIN REDUCTION LOOP --------------------

Continue = true;
while (Exponent >1 and Continue) do(
     
   Exponent = Exponent - 1;
   I = ideal(
      (diff Q)_(0,1), (diff Q)_(0,2),
      (diff Q)_(0,3), (diff Q)_(0,4),
      (-q)^Exponent
   );
   gbI = gb(I, ChangeMatrix=>true);
   if (not (Numerator) % gbI == 0) then (Continue = false; end);
   Reduced = (Numerator // gbI);
   Coeff = append(Coeff, Reduced_(4,0));
   Numerator =  (diff Reduced)_(0,1) + (diff Reduced)_(1,2)
              + (diff Reduced)_(2,3) + (diff Reduced)_(3,4);      

)

if (not Continue) then (<< "ERROR" << endl; end)

Coeff = append(Coeff, Numerator);      -- Add the final coefficient

-------------------- BEGIN OUTPUT -----------------------------

<<endl;
<< "The coefficients are:"<<endl;
i=0;
while(i < Deg+1) do (
   << "c" << Deg-i << " = " << Coeff#i << "." << endl;
   i=i+1;
);
<< " " <<endl;
\end{verbatim}

	\newpage 
	
\chapter{Maple Procedure to Calculate the Monodromy Representation}\label{maple}
Inputs:

Matrix $A = \sum_{i=1}^k\frac{R_i}{y-\alpha_i}$

Array VertexList $= [v_0,v_1,\ldots,v_m]$. A list of points defining a piecewise linear closed loop based at $v_0$. This is the path along which the monodromy is to be calculated. The paths must not intersect any of the singular points, $\alpha_1,\ldots,\alpha_m$.

\ 
\\
Outputs:

The trace of the monodromy around the closed loop and the square of the trace. As it stands, the trace is calculated to an accuracy of at least 3 decimal places in a typical example. This can be controlled by using the \emph{abserr} and \emph{relerr} options in Maple's dsolve() procedure, but is left as it is since we are only required to find the integer $Tr(M)^2$.

\begin{verbatim}
with(linalg);
Monodromy := proc(A, VertexList)

   local a, b, c, d, p0, p1, B, a11, a12, a21, a22;
   local eqn, icf, icg, solf, solg, i, noOfLines, M;
	
   noOfLines := nops(VertexList);
   a := 1; b := 0;
   c := 0; d := 1;
	
   for i from 1 to noOfLines do
      p0 := VertexList[i];
      p1 := VertexList[i mod noOfLines + 1];
      B := evalm((p1-p0)*A((p1-p0)*y+p0));
      a11 := B[1,1]; a12 := B[1,2];
      a21 := B[2,1]; a22 := B[2,2];
      eqn := {
         diff(f(y),y)=a11*f(y)+a12*g(y),
         diff(g(y),y)=a21*f(y)+a22*g(y)
      };
      icf := {f(0) = a, g(0)=c};
      solf := dsolve(eqn union icf, {f(y),g(y)}, type = numeric);
      a := rhs(solf(1)[2]);
      c := rhs(solf(1)[3]);
      icg := {f(0) = b, g(0)=d};
      solg := dsolve(eqn union icg, {f(y),g(y)}, type = numeric);
      b := rhs(solg(1)[2]);
      d := rhs(solg(1)[3]);
   end do;

   M := matrix(2,2,[a,b,c,d]);
   print(Tr(M) = trace(M));
   print(Tr(M)^2 = trace(M)^2);

end proc;
\end{verbatim}

To use this function, we need to specify the matrix $A=\sum\frac{R_i}{z-\alpha_i}$ and a path. In Maple, the composition of two paths defined by the lists of vertices \ttfamily{V1}\rm\ and \ttfamily{V2}\rm\ is the path defined by the concatenation of lists: \ttfamily{[op(V1),op(V2)]}\rm.

\begin{verbatim}

R0 := matrix(2,2,[ 1/6     , 0     ,-97111/233280 ,-5/3     ] );
Ra := matrix(2,2,[ 347/120 , 405/32,-99589/182250 ,-287/120 ] );
Rb := matrix(2,2,[-427/120 ,-15    , 207949/216000, 487/120 ] );
Rc := matrix(2,2,[ 1/2     , 75/32 , 0            , 0       ] );

A := z ->  R0/z + Ra/(z-9/64) + Rb/(z-1/7) + Rc/(z-5/32);

V0 := [127/896+I, -1+I, -1-I, 3/64-I, 3/64+I];
Va := [127/896+I, 3/64+I, 3/64-I, 127/896-I];
Vb := [127/896+I, 127/896-I, 67/448-I, 67/448+I];
Vc := [127/896+I, 67/448+I, 67/448-I, 1-I, 1+I];

Monodromy(A, [op(V0),op(Va)] );

                 Tr(M) = 3.162258490 + 0.000012081 I
                      2
                 Tr(M)  = 9.999878757 + 0.00007640648964 I 
\end{verbatim}
It is clear in this case that since $\textnormal{trace}(M)^2$ is approximating an integer, then $\textnormal{trace}(M) = \sqrt{10}$.

	\newpage

\typeout{Bibliography}
\providecommand{\bysame}{\leavevmode\hbox to3em{\hrulefill}\thinspace}
\providecommand{\MR}{\relax\ifhmode\unskip\space\fi MR }
\providecommand{\MRhref}[2]{%
  \href{http://www.ams.org/mathscinet-getitem?mr=#1}{#2}
}
\providecommand{\href}[2]{#2}

\addcontentsline{toc}{chapter}{Bibliography} 

\end{document}